\documentclass{article}
\usepackage[english]{babel}
\usepackage[style=numeric-comp, sortcites=true, doi=false, url=false, eprint=false]{biblatex}
\usepackage[letterpaper,top=2cm,bottom=2cm,left=3cm,right=3cm,marginparwidth=1.75cm]{geometry}
\usepackage{authblk}
\usepackage{amsmath}
\usepackage{mathrsfs}
\usepackage{graphicx}
\usepackage[colorlinks=true, allcolors=blue]{hyperref}
\usepackage{subfigure}
\usepackage{multirow}
\usepackage{booktabs}
\usepackage{amssymb}
\usepackage{bm}
\usepackage{enumitem}
\usepackage{makecell}
\usepackage{authblk}
\DeclareMathOperator{\sinc}{sinc}
\usepackage{mathtools}

\newcommand{\SSS}{\mathbb{S}}
\newcommand{\RR}{\mathbb{R}}
\newcommand{\ZZ}{\mathbb{Z}}
\newcommand{\Q}{\mathcal{Q}}
\newcommand{\ve}{\varepsilon}
\newcommand{\G}{\mathcal{G}}
\newcommand{\Tau}{\mathcal{T}}
\newcommand{\LL}{\mathcal{L}}
\DeclareMathOperator{\tr}{Trace}
\DeclareMathOperator{\diverg}{div}

\date{}

\title{An IMEX multi-level domain-decomposition method for the Boltzmann equation on unstructured meshes}

\author[1,2]{Domenico Caparello\thanks{Domenico.CAPARELLO@univ-cotedazur.fr, domenico.caparello@unife.it}}
\author[2,3]{Lorenzo Pareschi\thanks{l.pareschi@hw.ac.uk}}
\author[1]{Thomas Rey\thanks{Thomas.REY@univ-cotedazur.fr}}

\affil[1]{Université Côte d’Azur, CNRS, LJAD, Parc Valrose, F-06108 Nice, France}
\affil[2]{Department of Mathematics and Computer Science University of Ferrara}
\affil[3]{Maxwell Institute for Mathematical Sciences, Heriot-Watt University, Edinburgh, UK}

\begin{document}
\maketitle

\begin{abstract}
We introduce a dynamically adaptive multi-level method for the multiscale Boltzmann equation on unstructured meshes, coupling the compressible Euler, Navier–Stokes, and full Boltzmann descriptions within a single asymptotic-preserving IMEX framework. Unlike standard domain-decomposition approaches, in which independently evolved solvers exchange data only through interface conditions, the three models are integrated here as components of a unified multiregime scheme. They share the same globally stiffly accurate IMEX Runge–Kutta discretization and are synchronized at every internal stage through kinetic-to-fluid moment projections and fluid-to-kinetic lifting. This stage-wise coupling is essential for maintaining temporal consistency across model interfaces, while the asymptotic-preserving structure ensures compatibility with the fluid limits as the Knudsen number becomes small.
The active description is selected dynamically in each cell through moment-realizability indicators derived from successive Chapman–Enskog approximations. The resulting hierarchy uses Euler in near-equilibrium regions, Navier–Stokes where first-order nonequilibrium corrections are sufficient, and the full Boltzmann equation where a kinetic description is required. Employing Navier–Stokes as the intermediate level avoids evolving a velocity distribution in mildly nonequilibrium regions and is particularly advantageous on unstructured meshes, where kinetic transport represents the dominant computational cost. The spatial discretization combines finite-volume WENO reconstruction on triangular meshes with a conservative fast spectral treatment of the Boltzmann collision operator.
Numerical experiments on several test cases show that the new method closely reproduces the full Boltzmann reference while significantly reducing the computational cost. These results demonstrate that tightly coupled AP model hierarchies provide an effective alternative to conventional interface-based domain decomposition for deterministic kinetic simulations in complex geometries.
\end{abstract}

\smallskip
\noindent\textbf{Keywords:} Boltzmann equation; kinetic--fluid coupling;
dynamic domain decomposition; asymptotic-preserving schemes;
IMEX Runge--Kutta methods; unstructured meshes.

\tableofcontents

\section{Introduction}
Fluid flows in transitional regimes arise in many engineering and physical applications \cite{dim-par-2014, cer-1988}, including atmospheric reentry of spacecraft, micro-electro-mechanical systems (MEMS), and industrial heat transfer processes. In such regimes, classical macroscopic models such as the Euler equations are often inadequate due to shocks, kinetic boundary layers, and other non-equilibrium effects. A reliable description of these phenomena therefore requires kinetic models that resolve the microscopic behaviour of the gas. 
Moreover, to overcome the geometric constraints of structured grids, this work employs unstructured meshes, which are well suited to complex geometries. Such meshes permit a non-uniform distribution of elements, allowing increased spatial resolution in regions where it is most needed and providing a flexible balance between accuracy and computational cost.

Kinetic equations describe the evolution of collective statistical systems by evolving a particle distribution function in a high-dimensional phase space \cite{cer-1988}, typically six-dimensional for a three-dimensional configuration domain. Due to the intrinsic physical constraints of the model, the choice of numerical discretization is critical. Numerical discretizations should account for the fundamental structural properties  of kinetic equations \cite{dim-par-2014}, including positivity of the distribution  function, conservation of mass, momentum, and energy, and entropy dissipation. 
Preserving all these properties simultaneously remains a major challenge, especially for high-order methods.
Furthermore, in regimes where microscopic interactions dominate and collision processes become extremely frequent, the governing equations exhibit a high degree of stiffness. To handle these multiscale transitions efficiently without resolving the microscopic scales globally, asymptotic-preserving (AP) schemes are required \cite{fil-jin-2010,fil-jin-2011,jin-1999,jin-par-2001,DiPa13, cap-par-rey-2025-2, jin-li-2013}.

Dynamic model selection strategies based on the local fluid regime \cite{DiPa08, tiw-1998, dim-mie-ris-2014, Klar98, dim-par-2014, deg-dim-mie-2010, deg-jin-mie-2005, deg-liu-mie-2006, deg-dim-mie-2007, deg-dim-2012, cap-par-rey-2025, cap-ten-2025, fil-rey-2014} provide an effective approach to reducing the curse of dimensionality associated with the high-dimensional kinetic description. Rather than employing a single, computationally expensive solver across the entire domain, this approach allows for the adaptive partitioning of the flow into distinct regions. By dynamically selecting the most appropriate mathematical description, ranging from macroscopic fluid models in near-equilibrium zones to full kinetic solvers where non-equilibrium effects dominate, it becomes possible to concentrate high-dimensional resolution only where it is strictly required. In this way the computational cost is reduced, while the correct asymptotic limits are still preserved.

The framework developed in this work is inspired by these developments and by the foundational works of Levermore, Morokoff, and Nadiga \cite{lev-mor-nad-1998}, Filbet and Rey \cite{fil-rey-2014}, and Caparello et al. \cite{cap-par-rey-2025, cap-ten-2025}. Specifically, we adopt a decomposition strategy consistent with the frameworks established in \cite{fil-rey-2014, cap-par-rey-2025, cap-ten-2025}. In \cite{cap-par-rey-2025}, the computational domain is partitioned into Euler, ES-BGK, and full Boltzmann regimes on Cartesian meshes, where the use of the ES-BGK model as an intermediate step is justified by its reduced computational cost compared to the full Boltzmann collision operator. However, as we will demonstrate later in this work (see Section~\ref{secNumEff}), and as previously highlighted by Boscheri and Dimarco \cite{bos-dim-2021, bos-dim-2020}, the computational bottleneck on unstructured meshes shifts significantly. In such configurations, the predominant cost of the kinetic level arises from the transport term rather than the collision operator. Consequently, the ES-BGK model offers no clear computational advantage over the full Boltzmann solver in an unstructured framework. For this reason, we choose to employ the Navier--Stokes equations as the intermediate regime instead of the ES-BGK model. This substitution is further motivated by physical consistency: the ES-BGK model is specifically designed to recover the correct hydrodynamic limit up to the Navier--Stokes level \cite{fil-jin-2011, cap-par-rey-2025-2}, making the latter the natural macroscopic counterpart for our domain decomposition strategy.

While a similar Euler, Navier--Stokes, and BGK decomposition has been explored in a 1D setting \cite{xio-qiu-2017}, spatial adaptivity has also been considered in other multidimensional kinetic contexts. For instance, in \cite{cai-li-2010} the authors implement a two-level Euler with fast kinetic solver and Boltzmann-BGK coupling on two-dimensional $h$-adaptive meshes, while in \cite{wu-li-che-koo-tse-2006} unstructured meshes have been employed for a coupled Direct Simulation Monte Carlo (DSMC) and Navier--Stokes method. The present work differs from these approaches in several respects. The main novelty is that the Euler, Navier--Stokes, and Boltzmann models are not advanced as independent solvers coupled only through interface conditions. Instead, they are integrated within a unified AP IMplicit--EXplicit (AP-IMEX) Runge--Kutta framework, sharing the same temporal stages and exchanging kinetic and macroscopic representations at every stage \cite{bos-par-rus-2017, bos-par-rus-2024, cap-ten-2025}. This stage-wise synchronization preserves high-order temporal matching across static model interfaces and, together with the AP and globally stiffly accurate structure of the time discretization, ensures consistency with the fluid limit. A further distinctive feature is the combination of this tightly coupled three-level hierarchy with a deterministic spectral collision solver on general unstructured triangular meshes.
First, in the kinetic regions we solve the full Boltzmann equation instead of a simplified relaxation model, so that no phenomenological approximation of the collision operator is introduced. Second, the resulting hierarchy is physically more consistent. 
In simpler models where the kinetic level fails to recover the correct Prandtl number, its placement above the Navier--Stokes level may lead to a hierarchy that is not fully consistent across the different descriptions. In contrast, our strategy utilizes Navier--Stokes as a natural bridge between the equilibrium Euler regime and the fully non-equilibrium Boltzmann description. Finally, while 1D implementations focus almost exclusively on reducing collisional complexity, they do not address one of the primary challenges of modern CFD: unstructured spatial discretizations. By confining the full kinetic solver to small regions of the unstructured mesh, our strategy reduces the cost of the transport step, which is the most expensive part of the kinetic solver (see Section~\ref{secNumEff} and \cite{bos-dim-2020, bos-dim-2021}).

Overall, this hierarchical strategy offers two main advantages \cite{cap-par-rey-2025}:
\begin{itemize}
\item \textbf{Higher accuracy in multiscale regimes than classical fluid solvers:} The method automatically and dynamically applies the appropriate physical model in each region, ensuring the accurate resolution of shocks, boundary layers, and other kinetic effects where standard fluid models fail.

\item \textbf{Higher computational efficiency than full Boltzmann solvers:} The method is faster than full kinetic simulations, with comparable accuracy. This is achieved by confining the expensive kinetic computations to a minimal region of the domain while relying on highly efficient fluid solvers elsewhere.
\end{itemize}

One of the main numerical difficulties in solving the Boltzmann equation occurs in stiff regimes, where collisions are strong enough to induce stiffness but insufficient to drive the system to equilibrium \cite{dim-par-2014}. In such cases, reduced macroscopic models cannot be applied, and explicit time integration schemes are severely limited by stability constraints. A robust solution is offered by AP schemes, which remain stable and accurate across a wide range of Knudsen numbers, including the stiff limit, without the need for restrictive mesh or time-step refinement \cite{cor-per-1991, jin-1995, caf-shi-rus-1995, jin-par-tos-1998, cap-par-rey-2025-2, jin-par-2001, klar-1998, gos-tos-2002}.

Our numerical method incorporates these AP strategies to ensure robust behaviour across regimes, coupled with fast spectral solvers for the Boltzmann operator \cite{MoPa:2006,fil-mou-par-2006} to maximize performance. Fast spectral AP techniques provide an efficient deterministic approach for solving the full Boltzmann collision operator. Enhancing these techniques with a dynamic domain decomposition on unstructured meshes is one of the new aspects of this work. While other strategies for multiscale kinetic flows exist, such as Discrete Unified Gas Kinetic Schemes (DUGKS) \cite{zha-kun-2021,zha-wan-xu-2014,guo-xu-wan-2013, wu-shi-cha-wan-2016} and neural-network-based solvers 
\cite{tia-mar-2021,gan-tap-esm-myo-2024,lee-sch-xia-kru-fra-2024}, our approach retains the full Boltzmann collision operator in the kinetic regions. Its computational cost is reduced through a dynamic three-level hierarchy in which the Euler, Navier--Stokes, and Boltzmann models are tightly coupled, stage by stage, within a common AP-IMEX discretization.

Here, we develop and implement this hierarchical strategy in a fully two-dimensional spatial domain with two-dimensional velocity dependence on unstructured meshes, utilizing high-order accurate finite volume methods in space and time. In the remainder of the paper, we describe in detail the moment realizability-based switching criterion, the unified stage-wise IMEX coupling 
across the three models, and a Weighted Essentially Non-Oscillatory (WENO) reconstruction on unstructured meshes. Through benchmark problems involving shock waves and flows around complex 
obstacles, we show that the hybrid solutions closely reproduce the corresponding full Boltzmann reference while reducing the computational cost. These results illustrate the potential of the proposed framework for deterministic multiscale kinetic simulations in nontrivial geometries.

The remainder of this paper is organized as follows. We begin by introducing the kinetic and fluid models in Section~\ref{secPrelim}. Section~\ref{sec:UnifiedFramework} then presents the unified AP-IMEX scheme applied across all solvers. The criteria for the hybrid decomposition are established in Section~\ref{secRegimeIndic}, while the specific details of the numerical discretizations are outlined in Section~\ref{secNumSch}. To validate the proposed framework, Section~\ref{secNumSim} is dedicated to numerical experiments. This is followed by an analysis of the computational efficiency in Section~\ref{secNumEff}. Finally, concluding remarks and perspectives for future work are summarized in Section~\ref{secCon}.

\section{The Boltzmann equation}\label{secPrelim}
For a given non-negative initial condition $f_0$, we will study a particle distribution function $f = f(t,x,v)$, for $t \geq 0$, $x \in \Omega \subset \mathbb{R}^{d_x}$ and $v \in \mathbb{R}^{3}$, solution to the initial-boundary value problem
\begin{equation} 
    \label{eqCollision}
    \left\{ \begin{aligned}
      & \frac{\partial f}{\partial t} + v \cdot \nabla_x f \,=\, \frac{1}{\ve}\;\Q(f), 
      \\
    &  \;
      \\
      & f(0, x, v) = f_{0}(x,v),
    \end{aligned} \right.
\end{equation}
where the collision operator $\Q$ is a  Boltzmann-like operator, and $\ve > 0$ is the Knudsen number. This last quantity is the ratio between the mean free path of particles and the physical length scale of the system, and it gives an indication of the rarefaction of the gas: if $\ve \ll 1$ the particle system is subject to many collisions and then the gas is \emph{dense} and in fluid regime; vice-versa, if $\ve \sim 1$ the number of collisions is small and the gas can be considered \emph{rarefied} \cite{fil-rey-2014}.\\
The open set $\Omega$ is a bounded Lipschitz-continuous domain of $\RR^{d_x}$, which means  that the model \eqref{eqCollision} has to be supplemented with proper boundary conditions described later.
		
The Boltzmann operator satisfies the following properties:
\begin{enumerate}[label=\textbf{(H$\bm{_\arabic{*}}$)}, ref=\textbf{(H$\bm{_\arabic{*}}$)}]
    \item \label{hypConservations}  It preserves mass, momentum and kinetic energy, which can be respectively expressed as:
    \begin{equation*}
        \int_{\RR^{3}} \Q(f)(v) \, dv = 0, \quad  \int_{\RR^{3}} \Q(f)(v) \, v \, dv = 0, \quad \int_{\RR^{3}} \Q(f)(v) \, |v|^2 \, dv = 0;
    \end{equation*}
			 
    \item \label{hypEntropy} It dissipates the Boltzmann entropy (H-theorem)
    \begin{equation*}
        \int_{\RR^{3}} \Q(f)(v) \, \log(f)(v) \, dv \, \leq \, 0;
    \end{equation*}
        
    \item \label{hypEquilib} Its equilibria are given by Maxwellian distributions
    \begin{equation*}
        \Q(f) \, = \, 0 \quad \Leftrightarrow \quad f = \mathcal M_{\rho, \bm u, T} := \frac{\rho}{(2 \pi T)^{3/2}} \exp \left ( - \frac{|v-\bm u|^2}{2 T} \right ),
    \end{equation*}
    where the \emph{density}, \emph{velocity} and	\emph{temperature} of the gas $\rho$, $\bm u$ and $T$  are computed from the distribution function $f$ as  
    \begin{equation*}
        \rho = \int_{v\in{\RR}^{3}}f(v)\,dv, \quad \bm{u} = \frac{1}{\rho}\int_{v\in{\RR}^{3}}v f(v) \, dv, 
        \quad T = \frac{1}{3\rho} \int_{v\in{\RR}^{3}}\vert \bm{u} - v \vert^2 f(v) \,dv.
    \end{equation*}
  
\end{enumerate}  

According to assumptions \ref{hypEntropy}-\ref{hypEquilib}, it can be proved that, when $\ve \to 0$, the distribution $f^\ve$ converges (at least formally) to a Maxwellian distribution, whose moments are solution to the compressible Euler system
\begin{equation}
    \label{eqHydroClosedEuler}
    \left\{ \begin{aligned}
      & \partial_t \rho  + \diverg_x  (\rho \, \bm{u} ) \,=\, 0, 
      \\
      \,
      \\
      & \partial_t(\rho \, \bm{u} ) + \diverg_x  \left(\rho \, \bm{u} \otimes \bm{u} \,+\, \rho \, T  \,\bm{I}\right) \, =\, \bm{0}_{\RR^{3}}, 
      \\
      \,
      \\
      & \partial_t E + \diverg_x \left ( \bm{u} \left ( E +\rho \, T\right )  \right ) \,=\, 0,
    \end{aligned} \right.
\end{equation}
where $E$ represents the total energy 
    \begin{equation*}
        E= \int_{v\in{\RR}^{3}}\frac{1}{2}\vert v \vert^2 f(v) \,dv=\frac{1}{2}\rho |\bm{u}|^2+\frac{3}{2}\rho T.
    \end{equation*}
This corresponds to what is typically known as the \emph{hydrodynamic limit}.\\
At the next order in $\ve$, the Chapman--Enskog expansion (Appendix~\ref{app:ChapmanEnskog}) yields the compressible Navier--Stokes equations
\begin{equation}
\label{eqHydroNS}
    \left\{ \begin{aligned}
      & \partial_t \rho + \diverg_x (\rho\,\bm u) = 0,\\
      & \partial_t(\rho\,\bm u) + \diverg_x\left(\rho\,\bm u\otimes\bm u
        + \rho T\,\bm{I} - \ve\,\mu\,\bm D(\bm u)\right) = \bm 0_{\RR^3},\\
      & \partial_t E + \diverg_x\left((E+\rho T)\,\bm u
        - \ve\,\mu\,\bm D(\bm u)\,\bm u - \ve\,\kappa\,\nabla_x T\right) = 0,
    \end{aligned} \right.
\end{equation}
where $\bm D(\bm u)=\nabla_x\bm u+(\nabla_x\bm u)^\intercal
-\frac{2}{3}(\diverg_x\bm u)\,\bm{I}$ is the traceless deformation tensor, and the viscosity $\mu$ and thermal conductivity $\kappa$ are specified in Section~\ref{secRegimeIndic}.

In the following we present two operators that satisfy these hypotheses: the Boltzmann operator, which is the kinetic model of our hierarchy, and the ES-BGK operator, which is used only as the penalization operator in the time discretization of the Boltzmann equation (Section~\ref{sub:AP_Penalization_GSA}).
    
\subsection{The Boltzmann operator}
    \label{subBoltzOp}
    If we consider only binary elastic collisions between particles, then a dilute gas can be described by the kinetic equation with the Boltzmann operator, and in a velocity space of dimension $d_v$ the collision operator is

    \begin{align}
      \label{eqQBoltz}
      \Q (f,f)(v) & \,=\, \int_{\RR^{d_v}} \int_{\mathbb{S}^{d_v-1}}  B(|v-v_*|,\cos \theta) \, \left[ f'_* f' \,-\, f_* f \right] \, d\sigma \, dv_*,
        \end{align}
    where we used the notation $f = f(v)$, $f_* = f(v_*)$, $f ^{'} = f(v')$, $f_* ^{'} = f(v_* ^{'})$.   
        The pair $(v,v_*)$ are the pre-collisional velocities, while the quantities $(v',v'_*)$ are the post-collisional velocities. These two pairs are related by the following relations
    \begin{equation*}
      v' \,=\, \frac{v+v_*}{2} \,+\, \frac{|v-v_*|}{2} \,\sigma, \qquad v'_* \,=\, \frac{v+v_*}{2} \,-\, \frac{|v-v_*|}{2} \,\sigma.
    \end{equation*}
    The \emph{collision kernel} $B$ is a non-negative function which by physical arguments of invariance only depends on $|v-v_*|$ and	$\cos \theta = {\widehat u} \cdot \sigma$, where ${\widehat u} =(v-v_*)/|v-v_*|$. 
                            
    It can be proved that the Boltzmann collision operator has the fundamental properties \ref{hypConservations}, \ref{hypEntropy} and \ref{hypEquilib}.
    In the rest of the paper we will concentrate on the case of the so-called Maxwell pseudo-molecules \cite{dim-par-2014, cer-1988}, in which the collision rate is independent with respect to the relative velocities between particles:
    \begin{equation}
        B(|v-v_*|, \cos{(\theta))}=b_0 (\cos{(\theta)}).
    \end{equation}
          
    It is important to remark that, for this last class of models, it is possible to split the collision operator as follows
    \begin{equation}
        \Q(f,f) = \Q^+(f,f) - L(f)f,
    \end{equation}
    where 
    \begin{equation}
        \Q^+(f,f) = \int_{\RR^{d_v}}\int_{\SSS^{d_v-1}}B(|v-v_*|, \cos{(\theta)})f' f'_* d\sigma dv_*,
    \end{equation}
    and
    \begin{equation}
        L(f) = \int_{\RR^{d_v}}\int_{\SSS^{d_v-1}}B(|v-v_*|, \cos{(\theta)})f_* d\sigma dv_*.
    \end{equation}        
    This structure is particularly useful for applying fast spectral methods to compute the Boltzmann operator \cite{dim-par-2014}.
    The computation of Boltzmann-like collision operators using deterministic numerical methods requires restricting the problem to a bounded velocity space. A common approach is to periodize both the particle distribution function and the collision operator, effectively introducing some non-physical binary collisions, and then use Fourier series to compute the truncated operator. While this procedure may lead to the loss of certain local invariants, mass can be preserved through careful periodization. This forms the basis of spectral methods \cite{rey-2023}, which we discuss in more detail in Section \ref{secSpeMet}.

\subsection{The ES-BGK operator}
    \label{subBGK}
    The BGK operator \cite{bha-gro-kro-1954}, and its Ellipsoidal Statistical (ES-BGK) extension \cite{and-let-per-per-2000}, are two other well-known collision operators which satisfy the properties \ref{hypConservations}-\ref{hypEntropy}-\ref{hypEquilib}. 
    The bilinear Boltzmann collision operator $\Q$ is replaced by a nonlinear relaxation operator, which matches the same hydrodynamic limit as the Boltzmann operator: at the Euler level for the BGK operator, and at the Navier--Stokes level for the ES-BGK operator.

    In the classical BGK model \cite{bha-gro-kro-1954} the collision operator is a relaxation operator towards the local Maxwellian
    \begin{equation}
        \Q_{BGK}=\nu_{BGK} \left(\mathcal{M}_{\rho, \bm{u}, T} - f\right),
    \end{equation}
    where $\nu_{BGK}$ depends on macroscopic quantities. This simple model conserves mass, momentum and total energy, and has the correct Euler \eqref{eqHydroClosedEuler} limit when $\ve \to 0$. Unfortunately the viscosity $\mu$ and the thermal conductivity $\kappa$ at Navier--Stokes level are not satisfactory. In particular, for the BGK model the Prandtl number defined by \cite{fil-jin-2011}
    \begin{equation}
        \text{Pr}=\frac{\gamma}{\gamma-1}\frac{\mu}{\kappa},
    \end{equation}
    is equal to 1, whereas for most gases $\text{Pr}<1$. In this last equation, the quantity $\gamma=\left(d_v+2\right)/d_v$ is the ratio between specific heat at constant pressure and at constant volume, and $d_v$ is the number of dimensions of the velocity space.

    A possible way to overcome this limitation is to replace the Maxwellian $\mathcal{M}$, with an anisotropic Gaussian $\mathcal{G}$: this is the so called ES-BGK operator \cite{and-let-per-per-2000, cer-1988, hol-1965, fil-jin-2011}. In order to introduce this new quantity we need some additional notation. We begin by presenting the \emph{stress tensor}
    \begin{equation*}
      \rho_f\,\Theta_f(t,x)\, =\,  \int_{{\RR}^{3}}  (v-\bm{u}_f)\otimes (v-\bm{u}_f)\,f(t,x,v) \,dv.
    \end{equation*} 
    Therefore the \emph{translational temperature} is related to the stress tensor as $T_f = \tr(\Theta_f)/d_v$. We finally introduce the corrected tensor
    \begin{equation*}
      \Tau_f(t,x) \,\,=\,\, \left[(1-\beta) \, T_f \,\bm{I} \,\,+\,\,\beta \,\Theta_f\right](t,x),  
    \end{equation*} 
    which can be viewed as a linear combination of the initial stress tensor $\Theta_f$ and of the isotropic stress tensor $T_f \,\bm{I}$. 
    In three dimensions in the velocity space, the parameter $-\frac{1}{2}\leq \beta < 1$ is used to modify the value of the Prandtl number through the formula 
    \begin{equation*}
    \frac{2}{3} \,\leq\, {\rm Pr}  \,=\, \frac{1}{1-\beta} \,<\, +\infty\quad {\rm for } \quad \beta\in \Big[-\frac{1}{2}\,,\,1\Big).
    \end{equation*}
    Please note that the classical BGK operator, obtained for $\beta = 0$, has a Prandtl number equal to $1$.
    
    To define the ES-BGK operator, we introduce a corrected Gaussian $\G[f]$ defined by
    \begin{equation}
    \G[f]= \frac{\rho_f}{\sqrt{{\rm det}(2\pi\,\Tau_f)}}\,\exp\left(-\frac{(v-\bm{u}_f)^\intercal\,\Tau_f^{-1}\,(v-\bm{u}_f)}{2}\right),
    \end{equation}
    and the corresponding collision operator is now
    \begin{equation} 
      \label{eqOpES-BGK}
      \Q_{ESBGK}(f) =  \nu_{ESBGK} \left( \G[f] \,-\, f \right),
    \end{equation}
    where $\nu_{ESBGK}$ is the collision frequency, and it depends on macroscopic quantities.\\
    Unlike the BGK model, the ES-BGK model is able to reproduce transport coefficients, viscosity and thermal conductivity, in the Chapman--Enskog expansion \cite{bar-gol-lev-1991, bou-gol-pul-2000, fil-jin-2011, cap-par-rey-2025-2}, recovering the Navier--Stokes equations for the density $\rho$, momentum $\rho \bm{u}$ and temperature $T$ with the correct Prandtl number.\\
    This is one of the main advantages of the ES-BGK operator over the BGK operator. For a more comprehensive discussion on this topic, encompassing both single- and multi-species frameworks, we refer the reader to \cite{cap-par-rey-2025-2}.\\
    The details regarding the BGK and ES-BGK models provided in this section serve as a framework for the subsequent development of our numerical methodology. Specifically, these models are used in the construction of an Asymptotic Preserving (AP) Implicit-Explicit (IMEX) scheme for the full Boltzmann solver. A comprehensive discussion of this numerical framework and its implementation will be presented in Sections \ref{secNumSch} and \ref{sub:AP_Penalization_GSA}.

\section{The multi-level AP-IMEX framework}
\label{sec:UnifiedFramework}
 
 In this section, we present the high-order temporal discretization, the coupling between the different models, and the regime indicators that drive the dynamic domain decomposition.
 The formulation is designed to integrate a hierarchy of physical models within a single scheme, while retaining both the Asymptotic-Preserving (AP) property in the fluid limit and a high-order temporal matching at static interfaces between different physical regimes.

\subsection{Common IMEX time integration across the model hierarchy}
\label{sub:UnifiedIMEX}

To keep the framework general and flexible, the temporal discretization is designed to be common to a generic hierarchy of models $\{\mathscr{M}_0, \mathscr{M}_1, \dots, \mathscr{M}_Q\}$, characterized by increasing physical fidelity and computational cost.

To this end, let us consider a generic system of partial differential equations where the right-hand side can be additively split into two terms: a non-stiff operator $F(U)$, which can be safely treated explicitly, and a stiff operator $G(U)/\ve$, which requires an implicit treatment to avoid severe stability restrictions on the time step $\Delta t$. Here, $U$ represents a generalized state vector whose meaning depends on the active model. At time level $n$, each spatial cell is associated with explicit and implicit operators $F^{l_n}$ and $G^{l_n}/\ve$, respectively, corresponding to the locally active model $\mathscr{M}_{l_n}$ with $l_n \in \{0, \dots, Q\}$. Please note that, in general, the locally active model $\mathscr{M}_{l_n}$ varies from cell to cell. In this work $Q=2$, and $\mathscr{M}_0$, $\mathscr{M}_1$ and $\mathscr{M}_2$ denote the Euler, Navier--Stokes and Boltzmann models, respectively.

Regardless of the local nature of $U$ and the specific model active in the spatial cell, an IMEX Runge--Kutta scheme with $\nu$ stages can be written in the following general form \cite{bos-par-rus-2024}:
\begin{equation}
\label{eq:IMEX_system}
\begin{cases}
U_i \displaystyle = U^n + \Delta t \sum_{j=1}^{i-1} \tilde{a}_{ij} F^{l_n}(t^n + \tilde{c}_j \Delta t, U_j) + \frac{1}{\ve} \Delta t \sum_{j=1}^i a_{ij} G^{l_n}(t^n + c_j \Delta t, U_j),\\
U^{n+1} \displaystyle = U^n + \Delta t \sum_{i=1}^\nu \tilde{\omega}_i F^{l_n}(t^n + \tilde{c}_i \Delta t, U_i) + \frac{1}{\ve} \Delta t \sum_{i=1}^\nu \omega_i G^{l_n}(t^n + c_i \Delta t, U_i).
\end{cases}
\end{equation}
Here the subscripts $i$ and $j$ refer to the Runge--Kutta stages. The scheme is uniquely characterized by the $\nu \times \nu$ matrices $\tilde{A}=(\tilde{a}_{ij})$ (with $\tilde{a}_{ij}=0$ for $j \geq i$) and $A=(a_{ij})$, ensuring the scheme is explicitly evaluated for $F$ and implicitly evaluated for $G$. Along with the coefficient vectors $\tilde{c}$, $c$, $\tilde{\omega}$, and $\omega$, the method is conventionally represented by a double \emph{tableau} in Butcher notation:
\begin{center}
\begin{tabular}{c|c}
    $\tilde{c}$ & $\tilde{A}$\\                
    \hline
    \\
                & $\tilde{\omega}^\intercal$\\
\end{tabular}, \  \   
\begin{tabular}{c|c}
    $c$ & $A$\\                
    \hline
    \\
        & $\omega^\intercal$\\
\end{tabular}.
\end{center}
From an algorithmic and computational perspective, the strict Diagonally Implicit (DIRK) structure of the matrix $A$ ($a_{ij}=0$ for $j > i$) and the explicit structure of $\tilde{A}$ ($\tilde{a}_{ij}=0$ for $j \geq i$) are exploited to solve the system sequentially stage-by-stage. 
Following the formal splitting introduced in \cite{par-rus-2001}, let us denote the operator evaluations at the intermediate stages as $F_j = F^{l_n}(t^n + \tilde{c}_j \Delta t, U_j)$ and $G_j = G^{l_n}(t^n + c_j \Delta t, U_j)$. The computation of each internal stage $i = 1, \dots, \nu$ is practically decoupled into an explicit evaluation followed by an implicit step:
\begin{enumerate}
    \item \textbf{Explicit step:} First, an intermediate explicit state $U_i^*$ is computed using solely the explicit fluxes from the previously known stages:
    \begin{equation}
        U_i^* = U^n + \Delta t \sum_{j=1}^{i-1} \tilde{a}_{ij} F_j.
    \end{equation}
    \item \textbf{Implicit step:} Subsequently, the complete stage solution $U_i$ is evaluated by adding the history of the implicit operators and solving for the current implicit contribution:
    \begin{equation}
    \label{eq:ImplicitCorrector}
        U_i = U_i^* + \frac{\Delta t}{\ve} \sum_{j=1}^{i-1} a_{ij} G_j + \frac{\Delta t}{\ve} a_{ii} G^{l_n}(t^n + c_i \Delta t, U_i).
    \end{equation}
\end{enumerate}
Finally, after all $\nu$ stages have been computed sequentially, the numerical solution is advanced to the next time level using the standard quadrature update:
\begin{equation}
    U^{n+1} = U^n + \Delta t \sum_{i=1}^\nu \tilde{\omega}_i F_i + \frac{\Delta t}{\ve} \sum_{i=1}^\nu \omega_i G_i.
\end{equation}
This splitting is convenient from a computational point of view: instead of solving one large nonlinear system, at each stage only the diagonal term $a_{ii} G_i$ has to be inverted, see \eqref{eq:ImplicitCorrector}. The stiff operator can thus be inverted locally at each stage, making the high-order fully coupled IMEX integration practically affordable. Note that in many cases the implicit step based on a DIRK solver for $G$ can be solved explicitly. More details on this specific topic are discussed in Section \ref{sub:AP_Penalization_GSA}.

For the macroscopic models of the hierarchy (Euler and Navier--Stokes), the state vector is $U=(\rho,\rho\bm u,E)^\intercal$ and no stiff relaxation term is present: setting $G^{l_n}=0$, the explicit  operator $F^{l_n}$ contains the full set of macroscopic fluxes, including the viscous and heat fluxes of the Navier--Stokes equations. The IMEX scheme \eqref{eq:IMEX_system} then reduces to the fully explicit Runge--Kutta method defined by $(\tilde A,\tilde\omega)$, with the same stages as the kinetic solver.
We stress that the same double Butcher tableau (i.e., the same set of temporal coefficients) is used for all the models $\mathscr{M}_{l_n}$ of the hierarchy. Sharing the same temporal scheme across the whole computational domain is a key element of our coupling strategy.

As detailed in Section~\ref{sub:coupling}, synchronizing the internal stages of all the models preserves the high-order accuracy in time at the static interfaces between cells governed by different models \cite{cap-ten-2025}.

\subsection{Asymptotic-Preserving penalization and GSA property}
\label{sub:AP_Penalization_GSA}

The time discretization must handle the stiff limit ($\ve \to 0$) without time-step restrictions depending on $\ve$, i.e. it must be Asymptotic-Preserving (AP). In the context of kinetic theory, the state variable is the distribution function $U = f(t,x,v)$, and the evolution equation involves a convective transport operator $F(f)$ and the full spectral Boltzmann collision operator $G(f) = \mathcal{Q}(f)$. 

Spatially, we employ a Finite Volume (FV) framework. Denoting the spatial cell as $P_i$, the exact evolution equation for the cell average $\bar{f}_i(t)$ reads:
\begin{equation}
\bar{f}_i' = \bar{F}_i(f) + \frac{1}{\ve}\bar{G}_i(f).
\end{equation}
To achieve third-order spatial accuracy, we employ a third-order WENO reconstruction \cite{bal-gar-flo-bos-2020} alongside Lax--Friedrichs numerical fluxes \cite{hu-shu-1999} for the transport part. However, preserving this third-order accuracy for the highly stiff collision operator requires a multi-point quadrature rule within each cell. An implicit treatment of the Boltzmann operator is already prohibitively expensive at a single point; with the multi-point quadrature needed for third-order accuracy, such a system would have to be solved at every quadrature node of each cell, which makes it even more expensive.

To overcome this computational bottleneck while strictly retaining the high-order spatial accuracy, we rely on a penalization technique \cite{fil-jin-2010, bos-dim-2021, bos-par-rus-2024}. We introduce a simplified, cell-centered low-order operator, specifically an ES-BGK-type relaxation operator $\mathcal{P}_{ESBGK}(f)$ (defined in Equation \eqref{eqOpES-BGK}), which is treated implicitly, while the high-order accuracy is recovered through an explicit correction. In all simulations we take $\nu_{ESBGK}=6$ and $\beta=-1$, which in a two-dimensional velocity space gives ${\rm Pr}=1/(1-\beta)=1/2$, i.e. the Prandtl number of the Boltzmann operator for Maxwell molecules. Denoting the high-order multi-point approximation of the Boltzmann operator as $\bar{G}^{HO}$, and the low-order cell-centered ES-BGK operator as $\bar{\mathcal{P}}_{ESBGK}^{LO}$, the semi-discrete system is reformulated by adding and subtracting the penalization term:
\begin{equation}
\bar{f}_i' = \bar{F}_i(f) + \frac{1}{\varepsilon} \left( \bar{G}^{HO}(f) - \bar{\mathcal{P}}_{ESBGK}^{LO}(f) \right) + \frac{1}{\varepsilon} \bar{\mathcal{P}}_{ESBGK}^{LO}(f).
\end{equation} 
When the IMEX Runge--Kutta scheme \eqref{eq:IMEX_system} is applied to this reformulated system, the explicit operator is $F=\bar{F}_i+\ve^{-1}\big(\bar{G}^{HO}-\bar{\mathcal{P}}_{ESBGK}^{LO}\big)$ and the stiff implicit operator is $G/\ve$ with $G=\bar{\mathcal{P}}_{ESBGK}^{LO}$. The cell-centered ES-BGK operator is therefore treated implicitly, but it can be computed explicitly thanks to the strategy proposed in \cite{fil-jin-2011, cap-par-rey-2025-2}.
Since it is evaluated only at the cell center, its local inversion is analytically exact and computationally inexpensive. Conversely, the stiff non-equilibrium residual $\left( \bar{G}^{HO}(f) - \bar{\mathcal{P}}_{ESBGK}^{LO}(f) \right)$ is shifted to the explicit part of the scheme.

This strategy removes the algebraic bottleneck, but the stiff residual scaled by $1/\varepsilon$ is now treated explicitly, which may compromise the stability and the AP property of the scheme as $\varepsilon \to 0$.
To guarantee that the explicit intermediate errors are filtered out, we use \textit{Globally Stiffly Accurate} (GSA) schemes \cite{bos-par-rus-2013}. An IMEX scheme is GSA if its coefficients satisfy:
\begin{equation}
\omega_i = a_{\nu i} \quad \text{and} \quad \tilde{\omega}_i = \tilde{a}_{\nu i} \qquad \forall \,\, i=1,\dots,\nu.
\end{equation}
This structural algebraic constraint ensures that the updated solution at the new time step identically coincides with the final internal stage, i.e., $f^{n+1} \equiv f_\nu$. 
Rather than computing the final update as a linear combination of the stages, the GSA property ensures that the scheme ends with an implicit stage. 

In this work we use the GSA BPR(3,5,3) scheme \cite{bos-par-rus-2013} (also referred to as BPR(3,4,3) in \cite{bos-dim-2021}), whose explicit and implicit tableaus are given by:
\begin{center}
\begin{tabular}{c|c c c c c}
    $0$   & $0$   & $0$   & $0$   & $0$ & $0$\\ 
    $1$   & $1$   & $0$   & $0$   & $0$ & $0$\\
    $2/3$ & $4/9$ & $2/9$ & $0$   & $0$ & $0$\\
    $1$   & $1/4$ & $0$   & $3/4$ & $0$ & $0$\\
    $1$   & $1/4$ & $0$   & $3/4$ & $0$ & $0$\\
    \hline
          & $1/4$ & $0$   & $3/4$ & $0$ & $0$\\
\end{tabular}, \  \  \  \
\begin{tabular}{c|c c c c c}
    $0$   & $0$    & $0$    & $0$   & $0$    & $0$\\ 
    $1$   & $1/2$  & $1/2$  & $0$   & $0$    & $0$\\
    $2/3$ & $5/18$ & $-1/9$ & $1/2$ & $0$    & $0$\\
    $1$   & $1/2$  & $0$    & $0$   & $1/2$  & $0$\\
    $1$   & $1/4$  & $0$    & $3/4$ & $-1/2$ & $1/2$\\
    \hline
          & $1/4$ & $0$   & $3/4$ & $-1/2$ & $1/2$\\
\end{tabular}.
\end{center}

\subsection{Interface coupling}
\label{sub:coupling}

Once the unified temporal structure is established, the numerical solution must be coupled across neighbouring cells belonging to different physical regimes. In what follows we implement the strategy proposed in \cite{cap-ten-2025}. Let $\mathcal{T}$ denote the computational mesh and let
\begin{equation}
    \chi_i \in \{\mathrm{K},\mathrm{E},\mathrm{NS}\},
\end{equation}
be the regime indicator associated with cell $i$, where $\mathrm{K}$ denotes a kinetic cell, $\mathrm{E}$ an Euler cell, and $\mathrm{NS}$ a Navier--Stokes cell. Furthermore, let $\mathcal{S}_i$ be the numerical stencil associated with cell $i$, i.e. the set of all cells involved in the update of $P_i$, including the reconstruction stencils of its face neighbours. The kinetic solver evolves the distribution function $f$, whereas the Euler and Navier--Stokes solvers evolve the vector of conservative variables $U=(\rho,\rho \bm u,E)^\intercal$. Consequently, all stencil-dependent quantities must be expressed in a representation consistent with the governing model of the central cell.

\paragraph{Coupling for kinetic cells}
Consider a cell $i$ such that $\chi_i=\mathrm{K}$. Since the kinetic discretization requires a distribution function in every cell belonging to $\mathcal{S}_i$, we define the coupled distribution $\widetilde f_j$ for $j\in\mathcal S_i$ as:
\begin{equation}
\label{eq:kinetic_coupling}
\widetilde f_j=
    \begin{cases}
        f_j,
        & \chi_j=\mathrm{K},
        \\[2mm]
        \mathcal M(U_j),
        & \chi_j\in\{\mathrm{E},\mathrm{NS}\},
    \end{cases}
\end{equation}
where $\mathcal M(U_j)$ denotes the Maxwellian distribution associated with the conservative variables $U_j$. The kinetic operator in cell $i$ is then evaluated using the set $\{\widetilde f_j\}_{j\in\mathcal S_i}$.

\paragraph{Coupling for macroscopic cells}
Consider a cell $i$ such that $\chi_i\in\{\mathrm{E},\mathrm{NS}\}$. The macroscopic discretization requires conservative variables throughout the stencil. We define the coupled state $\widetilde U_j$ for $j\in\mathcal S_i$ as:
\begin{equation}
\widetilde U_j=
\begin{cases}
    U_j,
    & \chi_j\in\{\mathrm{E},\mathrm{NS}\},
    \\[3mm]
    \displaystyle \int_{\mathbb{R}^d} m(v)\,f_j(v)\,dv,
    & \chi_j=\mathrm{K},
\end{cases}
\label{eq:macro_coupling}
\end{equation}
where $m(v) = (1,\, v,\, |v|^2/2)^\intercal$ is the vector of collision invariants. The macroscopic operator in cell $i$ is then computed using the set $\{\widetilde U_j\}_{j\in\mathcal S_i}$.

\paragraph{Dynamic adaptation and stage-wise synchronization}
The bidirectional coupling detailed above defines a local projection between kinetic and macroscopic representations: kinetic-to-macroscopic transfer is achieved through moment projection, whereas macroscopic-to-kinetic lifting is performed via Maxwellian reconstruction. This procedure is applied at two stages of the temporal evolution:
\begin{enumerate}
\item \textbf{Initialization at the start of the time step:} At the beginning of each time step $t^{n+1}$, the domain decomposition map is dynamically updated from the previous state $\chi_i^n$ to the new partition $\chi_i^{n+1}$ according to local physical criteria. Before evaluating any spatial fluxes or RK stages, the solution state is immediately synchronized across newly transitioned cells and stencil halos using the underlying fields available from step $t^n$:
\begin{itemize}
\item \textbf{Local remapping:} Cells undergoing a regime transition update their primary internal variables. If $\chi_i^n = \mathrm{K} \to \chi_i^{n+1} \in \{\mathrm{E}, \mathrm{NS}\}$, the newly required conservative moments $U_i^n$ are extracted via quadrature of $f_i^n$. If $\chi_i^n \in \{\mathrm{E}, \mathrm{NS}\} \to \chi_i^{n+1} = \mathrm{K}$, the discrete distribution function is instantiated as the equilibrium state $f_i^n = \mathcal{M}(U_i^n)$. Transitions between macroscopic regimes ($\mathrm{E} \leftrightarrow \mathrm{NS}$) retain the primary conservative state $U_i^n$ directly.
\item \textbf{Stencil preparation:} To evaluate spatial transport, data required from neighbouring stencil cells $j \in \mathcal{S}_i$ is gathered based on their \textit{previous} regime state at $t^n$. If cell $i$ is now kinetic ($\chi_i^{n+1} = \mathrm{K}$), it requires distribution functions throughout its stencil; if neighbour $j$ was macroscopic ($\chi_j^n \in \{\mathrm{E}, \mathrm{NS}\}$), its kinetic distribution function is reconstructed via Maxwellian $\mathcal{M}(U_j^n)$, otherwise the existing $f_j^n$ is retained. Conversely, if cell $i$ is now macroscopic ($\chi_i^{n+1} \in \{\mathrm{E}, \mathrm{NS}\}$), it requires macroscopic moments; if neighbour $j$ was kinetic ($\chi_j^n = \mathrm{K}$), its moments are recovered via quadrature, otherwise the existing $U_j^n$ is directly used.
\end{itemize}
\item \textbf{Synchronization within the time step:} 
To preserve the high-order temporal accuracy of the IMEX scheme across physical boundaries throughout the time step, this bidirectional coupling is repeated after each implicit update of the Runge--Kutta integrator, using the regime partition $\chi_i^{n+1}$ valid for that step. Exchanging interface information only at the end of a full time step (step-wise coupling) would introduce an artificial temporal lag at the regime interfaces, degrading the local boundary accuracy to first order \cite{cap-ten-2025}. Stage-by-stage synchronization guarantees that interface states in every numerical stencil remain synchronized with the intermediate RK stage solutions.
\end{enumerate}
We stress that this stage-by-stage synchronization guarantees the high-order temporal accuracy of the IMEX scheme only for static interfaces \cite{cap-ten-2025}. In our framework the regime partition evolves in time; nevertheless, the dynamic initialization and the stage-wise coupling avoid the first-order temporal lag that a step-wise exchange would introduce at the regime boundaries.

\paragraph{Lifting operator}
In principle, alternative kinetic reconstructions could be employed when lifting macroscopic Navier--Stokes states to the kinetic description, such as higher-order non-equilibrium closures based on Chapman--Enskog type expansions or other kinds of closure approaches \cite{lev-1996, cli-mcd-2009}. However, these higher-order reconstructions do not, in general, preserve the non-negativity of the resulting distribution function. The occurrence of negative values in the distribution function near interfaces can lead to numerical instabilities in a coupled kinetic--fluid framework. For this reason, and to improve the robustness of the coupling, we restrict the macroscopic-to-kinetic lifting operator to the Maxwellian distribution, which is non-negative and uniquely determined by the local conservative variables, thus ensuring the physical admissibility of the lifted distribution.

\subsection{Regime indicators}
\label{secRegimeIndic}

 Several decomposition methods rely on the continuum breakdown parameter introduced by Boyd, Chen, and Candler \cite{boy-che-can-1995}, namely a local Knudsen number based on the gradient length of the macroscopic quantities, to switch between hydrodynamic and kinetic descriptions. Large Knudsen numbers correspond to rarefied regimes where a kinetic description becomes necessary, whereas hydrodynamic models are employed in near-equilibrium regions.

Several authors have developed methods of this type, incorporating spatial derivatives directly into the numerical scheme \cite{kol-ars-ari-fro-zab-2007, deg-dim-2012}, or relying on indicators based on the viscous stress and heat flux derived via moment expansions \cite{tiw-1998, tiw-kla-har-2009, ala-pup-2011, ala-pup-2012}. Since the literature on this topic is very large, we only cite here the works most closely related to the present approach: \cite{DiPa08, dim-mie-ris-2014, Klar98, deg-dim-mie-2010, deg-jin-mie-2005, deg-liu-mie-2006, deg-dim-mie-2007, fil-rey-2014, cap-par-rey-2025, cap-ten-2025, xio-qiu-2017}. Additionally, hybrid schemes have recently been extended to Vlasov--BGK and Vlasov--Poisson--BGK equations \cite{lai-2023, lai-rey-2023, dim-mie-ris-2014}, highlighting the growing interest in domain decomposition to speed up kinetic solvers.

The present approach adopts a decomposition strategy consistent with the frameworks established in \cite{fil-rey-2014, cap-par-rey-2025, cap-ten-2025}, which builds upon classical and recent hybrid solvers \cite{tiw-2000, fil-xio-2018, li-son-wan-2021, xio-qiu-2017}. Specifically, while on Cartesian meshes a three-regime partition (Euler, ES-BGK, and full Boltzmann) is computationally justified \cite{cap-par-rey-2025}, the computational bottleneck on unstructured meshes shifts significantly. As highlighted by Boscheri and Dimarco \cite{bos-dim-2021, bos-dim-2020}, the predominant cost on unstructured grids arises from the transport term rather than the collision operator. Consequently, the ES-BGK model offers no clear computational advantage over the full Boltzmann solver. 

For this reason, we employ the compressible Navier--Stokes equations as our intermediate regime instead of the ES-BGK model. This choice is physically consistent, as the ES-BGK model is designed to recover the correct hydrodynamic limit up to the Navier--Stokes level \cite{fil-jin-2011, cap-par-rey-2025-2}, making the latter the natural macroscopic counterpart. We also implement the high-order in time coupling technique developed in \cite{cap-ten-2025}, discussed and numerically validated in Section \ref{secNumSim}.

We defer the formal Chapman--Enskog expansion and the analytical derivations of the hydrodynamic equations to Appendix~\ref{app:ChapmanEnskog}. 

   The transport coefficients are the Chapman--Enskog values for Maxwell molecules used in \cite{fil-jin-2011}: $\mu=\frac{\sqrt2}{3\pi}\frac{T}{A_2(5)}$ with $A_2(5)\simeq0.436$, and $\kappa=4\mu$, corresponding to the Prandtl number ${\rm Pr}=1/2$ of the two-dimensional Boltzmann operator \cite{fil-jin-2011}. These coefficients are used consistently in all the solvers. 

\subsubsection{The moment realizability matrix and regime indicators}
\label{subMomentRealizMain}

To construct our breakdown criteria, we introduce the moment realizability matrix $\bm{M}$, which captures the internal structure of the distribution function $f$ through its collisional invariants \cite{lev-mor-nad-1998, fil-rey-2014, cap-par-rey-2025, cap-ten-2025}. Heuristically, this matrix encodes whether the computed macroscopic moments lie within the physically admissible set generated by a positive distribution function $f$. When $\bm{M}$ is no longer positive semi-definite, the macroscopic closure becomes unreliable, indicating hydrodynamic breakdown.

Following the formulation proposed in \cite{cap-ten-2025} (with theoretical details in \cite{lev-mor-nad-1998, fil-rey-2014}), the matrix $\bm{M}$ is defined as:
\begin{equation}
\label{eq:MRM_Main}
    \bm{M} = \frac{1}{\rho}\int (\bm{m} \otimes \bm{m}) f \, dv = 
    \begin{pmatrix}
     1  & \bm{0}^\intercal            & -1\\
     \bm{0}  & \bm{I}+\bm{\bar{A}}    & \bm{\bar{B}}\\
     -1 & \bm{\bar{B}}^\intercal & \bar{C}
    \end{pmatrix},
\end{equation}
where $\bm{m} = \left(1, \frac{v-\bm{u}}{\sqrt{T}}, \frac{|v-\bm{u}|^2}{2T}-\frac{5}{2}\right)^\intercal$, $\bm{V}=(v-\bm{u})/\sqrt{T}$, and the moment-based quantities $\bm{\bar{A}}$, $\bm{\bar{B}}$, and $\bar{C}$ are defined as:
\begin{equation}
    \label{eqDefABV}
    \left\{\begin{aligned}
    & \bm{\bar{A}} := \frac{1}{\rho} \int_{\RR^{3}} \bm{A}(\bm{V}) f(v) \, dv, && \bm{A}(\bm{V}) = \bm{V} \otimes \bm{V} - \frac{|\bm{V}|^2}{3} \bm{I}, \\
    & \bm{\bar{B}} := \frac{1}{\rho}\int_{\RR^{3}} \bm{B}(\bm{V}) f(v) \, dv, && \bm{B}(\bm{V}) = \frac12 \left [ |\bm{V}|^2 - 5\right ] \bm{V}, \\
    & \bar{C} = \frac{1}{\rho}\int\left(\frac{|v-\bm u|^2}{2T}-\frac{5}{2}\right)^2 f(v) \, dv.
    \end{aligned}\right.
\end{equation}

By exploiting Hölder's inequality (see Appendix~\ref{app:Realizability} for the mathematical proof), it can be shown that $\bar{C} > 1$. This inequality allows us to compute the Schur complement of $\bm{M}$ in two successive steps. The positive semi-definiteness of $\bm{M}$ is mathematically equivalent to the positive semi-definiteness of the reduced matrix $\mathcal{S}$:
\begin{equation}\label{eqmatS_Main}
    \mathcal{S} = \bm{I}+\bm{\bar{A}} - \displaystyle\frac{\bm{\bar{B}}\bm{\bar{B}}^\intercal}{\bar{C}-1}.
\end{equation}
Evaluating $\mathcal{S}$ at successive orders of the Chapman--Enskog expansion yields three distinct indicators ($\mathcal{S}_0$, $\mathcal{S}_1$, and $\mathcal{S}_2$), which serve as the foundation for our model-switching criteria.

\subsubsection{Breakdown criteria}
\label{subBreakdownMain}

Our domain decomposition solver evaluates the reduced matrices $\mathcal{S}_\ell$ in each spatial cell to dynamically determine whether a transition between description layers is required. If the description at the $\ell$-th order layer is no longer sufficiently accurate compared to the $(\ell+1)$-th order layer, the cell transitions to the higher description level.

Evaluating the terms $\bm{\bar{A}}$, $\bm{\bar{B}}$, and $\bar{C}$ (with $\bar{C}_{\text{Euler}} = \bar{C}_{\text{NS}} = 5/2$) at the zeroth and first orders yields:
\begin{equation}
    \mathcal{S}_0 = \bm{I},
\end{equation}
and
\begin{equation}
    \mathcal{S}_1 = \bm{I} - \varepsilon \bm{D}(\bm u)\frac{\mu}{\rho T} - \frac{2}{3}\left(\varepsilon \frac{\kappa}{\rho T^{\frac{3}{2}}}\right)^2 \left(\nabla_x T\right) \left(\nabla_x T\right)^\intercal,
\end{equation}
where $\bm{D}(\bm u)$ is the traceless deformation tensor of $\bm{u}$. 

The moment realizability matrix~\eqref{eq:MRM_Main} and the indicator
$\mathcal S_1$ are derived here, as in \cite{cap-par-rey-2025, fil-rey-2014}, for a three-dimensional velocity space. In the simulations of
Section~\ref{secNumSim}, which employ a two-dimensional velocity space,
$\mathcal S_1$ is evaluated with the same expression, using the two-dimensional macroscopic fields ($\partial_z\equiv 0$, $u_z=0$), and its eigenvalues are computed on the full $3\times3$ matrix. With respect to a derivation carried out directly in a two-dimensional velocity space, this choice only modifies the numerical constants appearing in $\bm D(\bm u)$ and in the heat-flux term (e.g.\ $2/3$ instead of $1$). These constants do not affect the qualitative behaviour of the indicator and are absorbed in the empirical calibration of the threshold $\eta_1$. The surrogate indicator $\lambda_2$ introduced below is a scalar quantity and is used without modification.

\paragraph{Euler to Navier--Stokes Transition (Layer 0 $\to$ Layer 1)}
The Euler description is considered inadequate when the discrepancy between the eigenvalues $\lambda_{1}$ of $\mathcal{S}_{1}$ and those of $\mathcal{S}_{0}$ ($\lambda_{0}=1$) becomes too large. Thus, the switch from Euler (Layer 0) to Navier--Stokes (Layer 1) occurs if:
\begin{equation}
\label{crit:CompEuler_Main}
    \exists \lambda_{1} \in\, {\rm Sp}(\mathcal{S}_{1}) \qquad \text{such that} \,\, \left| \lambda_{1}-1\right| > \eta_1,
\end{equation}
where $\eta_1 > 0$ is a small, problem-dependent threshold. Conversely, the switch back from Layer 1 to Layer 0 is allowed only when:
\begin{equation}
    |\lambda_{1}(t, x)-1|\leq\eta_1 \qquad \forall \,\, \lambda_{1} \in\, {\rm Sp}(\mathcal{S}_{1}).
\end{equation}

\paragraph{Navier--Stokes to Boltzmann Transition (Layer 1 $\to$ Layer 2)}
In principle, the breakdown of the first-order (Navier--Stokes) description could be monitored by directly evaluating the Schur complement $\mathcal{S}_2$ derived from the second-order (Burnett) Chapman--Enskog expansion (detailed in Appendix~\ref{app:ChapmanEnskog}). However, as noted in \cite{fil-xio-2018}, the explicit evaluation of the $O(\varepsilon^2)$ Burnett corrections $\bm{\bar{A}}_{Burnett}$ and $\bm{\bar{B}}_{Burnett}$ involves many nonlinear terms. Computing $\mathcal{S}_2$ exactly in each cell at every time step would therefore introduce a significant computational overhead.

To circumvent this, following \cite{fil-xio-2018}, we adopt a computationally efficient scalar surrogate indicator $\lambda_2$. By assuming macroscopic quantities $\rho, T \sim O(1)$ and transport coefficients $\mu/\sqrt{T} \sim O(1)$, $\lambda_2$ isolates and measures the dominant spatial derivative scales of the $O(\varepsilon^2)$ non-equilibrium stresses and heat fluxes:
\begin{equation}
    \lambda_2 := \varepsilon^2 \left[ \frac{|\nabla_x T|^2}{T} + |\nabla_x \bm{u}|^2 + \sqrt{\left( |\Delta_x \bm{u}|^2 + |\Delta_x \rho / \rho|^2 \right)(1+T^2)} \right].
\end{equation}
The transition from the Navier--Stokes description (Layer 1) to the full Boltzmann kinetic solver (Layer 2) is triggered when:
\begin{equation}
\label{criterioSecondoSwitch}
    |\lambda_2| > \eta_2,
\end{equation}
while the reverse transition (Layer 2 $\to$ Layer 1) is permitted only when $|\lambda_2| \leq \eta_2$.
This indicator-based switching strategy has already been validated in \cite{fil-xio-2018, cap-par-rey-2025}.

Both transitions towards lower levels rely on the same macroscopic
indicators, which are derived from the Chapman--Enskog expansion. The
kinetic-to-fluid transition therefore implicitly assumes that, in kinetic cells, the distribution function is close to its Chapman--Enskog form. This assumption may fail for strongly non-equilibrium initial data; in such cases the first time step is performed with the kinetic model in the whole domain (see \hyperref[sec:test3]{Test 3}). A kinetic-side criterion based on the distance between $f$ and its local Maxwellian \cite{fil-rey-2014, cap-par-rey-2025} could be added to further control this transition, at the price of an additional threshold; we did not pursue this option here, in order to keep the number of free parameters minimal.

    \section{Numerical schemes}\label{secNumSch}
    In this section we discuss the spatial and velocity  discretization of our scheme and the conservation of macroscopic quantities at the kinetic level. All the details regarding the temporal discretization and the coupling between different solvers are discussed in Section \ref{sec:UnifiedFramework}.
    
    \subsection{Velocity discretization}
    In this subsection, we introduce the velocity space discretization of the distribution function governed by the full Boltzmann equation. It is worth noting that this procedure applies exclusively to the kinetic model, and not to the macroscopic equations (i.e., the Euler and Navier--Stokes equations). To perform the numerical approximation, the unbounded velocity space is truncated to a finite domain, thereby neglecting the tails of the distribution function. Since the exact distribution function generally has non-compact support, this truncation inevitably leads to a loss of exact conservation for the macroscopic quantities. The same issue affects the evaluation of the Boltzmann collision integral. To overcome this lack of conservation, we implement the strategy proposed by \cite{gam-tha-2009}, which is discussed in Section \ref{secConMas}. This bounded space is then discretized using a finite set of points, referred to as discrete velocities, which represent the only allowable microscopic states for the particles:
    \begin{equation}
        \mathcal{V} = \{ v_k = k \Delta v + a, \quad k \in \mathbb{Z}^{d_v}, \quad -K^{(j)} \le k_j \le K^{(j)}, \quad j=1,\dots,d_v \},
    \end{equation}
    where $a \in \mathbb{R}^{d_v}$ is an arbitrary shift vector, $\Delta v$ is a constant velocity mesh size, and $K^{(j)}$ represents the bounds for the components of the multi-index $k$. In this discrete-ordinate framework, the continuous distribution function $f$ is replaced by a vector $f_K(t,x)$ of size $N$. Each component of this vector approximates the continuous distribution function evaluated at the corresponding discrete velocity node $v_k$: 
    \begin{equation}
        f_K(t, x) = (f_k(t,x))_{k=1}^N, \quad \text{with} \quad f_k(t,x) \approx f(t,x,v_k).
    \end{equation}
    Consequently, the discrete-ordinate kinetic model consists of the following system of $N$ coupled equations:
    \begin{equation}
        \partial_t f_k + v_k \cdot \nabla_{x} f_k = \frac{1}{\ve}\, \mathcal{Q}_k(f), \quad k=1, \dots, N.
    \end{equation}
    
    \subsection{Spatial discretization}
    The spatial discretization is based on a finite volume framework combined with a WENO reconstruction on triangular meshes. Specifically, we consider a centroid-based Delaunay-type domain tessellation composed of $N_P$ non-overlapping triangles $P_i$, for $i=1, \dots, N_P$. We denote by $x_i$ the position vector of the barycenter of each triangular cell. Let us define the data interpolation layer on this mesh. The cell averages of each component of the vector representing the distribution function at time $t^n$ are obtained via:
    \begin{equation}
        \bar{f}_{k,i}^n = \frac{1}{|P_i|} \int_{P_i} f_k(t^n, x) \, dx, \quad \forall k \in \{1, \dots, N\}, \quad i \in \{1, \dots, N_P\},
    \end{equation}
    where $|P_i|$ represents the volume (or area) of the cell $P_i$. Similarly, the cell averages for the macroscopic quantities evolved using the Euler or Navier--Stokes equations are defined as:
    \begin{equation}
    \bar{U}_{b,i}^n = \frac{1}{|P_i|} \int_{P_i} U_b(t^n, x) \, dx, \quad i \in \{1, \dots, N_P\},
    \end{equation}
    where $b \in \{\text{Euler}, \text{Navier--Stokes}\}$ identifies the specific macroscopic model, and $U_b(t, x) = [\rho, \allowbreak \rho u_x, \allowbreak \rho u_y, \allowbreak E]_b^\intercal(t, x)$ denotes the continuous vector of the macroscopic conservative variables. Our goal is now to construct a high-order, non-oscillatory polynomial representation of both the distribution function $f_k(t,x)$ and the macroscopic variables $U_b(t,x)$. These local spatial reconstructions will be denoted as $f^w_{k,i}(t^n,x)$ and $U_{b,i}^w(t^n, x)$, respectively. To achieve this, we implement the third-order WENO reconstruction on unstructured meshes described in Section 2 of \cite{bal-gar-flo-bos-2020}. The numerical fluxes at the cell interfaces (edges) are evaluated at high-order Gaussian quadrature points using a classical Lax--Friedrichs flux, according to the approach outlined in \cite{hu-shu-1999}.

    \subsection{Spectral discretization of the Boltzmann operator}\label{secSpeMet}
    The aim is to build a numerical method to compute the Boltzmann operator \eqref{eqQBoltz}, which depends only on the velocity variable $v$. We shall then neglect the time and space variables in this section.
    The starting point is the truncation of the integration domain of the Boltzmann integral \eqref{eqQBoltz}, then we suppose that the distribution function $f$ has compact support on the ball $B_0(R)$, of radius $R$ centred in the origin. We have the classical result from \cite{PaRu:SINUM:2000}, which states that if the distribution function $f$ is compactly supported on the ball $B_0(R)$ of radius $R$ centered in the origin, then $\text{Supp}_v(\Q(f,f)) \subset B_0 (\sqrt{2}R)$.
    It is possible to introduce a collision operator with cut-off in the following way   
    \begin{equation}
        \Q^R(f,g)(v)=\int \int_{B_0(2R)\times \SSS^{d-1}}\left[f(v')g(v_*')-f(v)g(v_*)\right]B(v-v_*, \sigma)dv_* d\sigma, \ \ \ \forall v \in \RR^d.
    \end{equation}
    In order to have a spectral approximation without superposition of periods, it is sufficient to have a distribution function $f(v)$ restricted on the cube $\left[-T, T\right]^d$, with $T\geq (2+\sqrt{2})R$. Then, we assume $f(v)=0$ on $\left[-T, T\right]^d \backslash B_0(R)$ and extend $f(v)$ to a periodic function on the set $\left[-T, T\right]^d$. We want to remark that the lower bound for $T$ can be improved: in fact, the choice $T=(3+\sqrt{2})R/2$ guarantees the absence of intersection between periods where the distribution function $f$ is different from 0. Unfortunately, the support of $f$ will increase in time, so we can just minimize the errors due to aliasing \cite{rey-2023}.\\                
    In the following we will assume $T=\pi$ and hence $R=2\pi/(3+\sqrt{2})$ in order to simplify the notation.
    Projecting the collision operator onto the space of trigonometric polynomials we obtain
    \begin{equation}\label{eqColOpTrigPol}
        \hat{Q}_k = \int_{\left[-\pi, \pi\right]^{d_v}} \Q^R(f_N, f_N)e^{-ik\cdot v}dv, \ \ \ k = -N, ..., N,
    \end{equation}
    where $f_N$ is the truncated Fourier series of the distribution function $f$,
    \begin{equation}
        \begin{aligned}                            
            f_N(v)=\sum_{k=-N}^{N} \hat{f}_k e^{ik\cdot v},\\
            \hat{f}_k = \frac{1}{(2\pi)^{d_v}}\int_{[-\pi,\pi]^{d_v}} f(v) e^{-ik\cdot v}dv,
        \end{aligned}
    \end{equation}
    where we have used just one index to denote the \emph{$d_v$}-dimension sums with respect to the vector $k=(k_1, ..., k_{d_v})\in \ZZ^{d_v}$, and so we have
    \begin{equation}
        \sum_{k=-N}^N \coloneq \sum_{k_1, ..., k_{d_v} = -N}^N.
    \end{equation}
    Now we can substitute the truncated Fourier series $f_N$ inside the expression \eqref{eqColOpTrigPol}, and after some computations we obtain \cite{PaRu:SINUM:2000}
    \begin{equation}\label{eqQk}
        \hat{Q}_k = \sum_{\substack{l,m=-N\\l+m=k}}^N \hat{f}_l \hat{f}_m \beta(l, m), \ \ \ k = -N, ..., N,
    \end{equation}
    where $\beta(l,m)=\hat{B}(l,m) - \hat{B}(m,m)$, and
    \begin{equation}
        \hat{B}(l,m)=\int_{B_0(2\lambda \pi)}\int_{\SSS^{d-1}} |q|\sigma(|q|, \cos{(\theta)})e^{-i(l\cdot q^+ + m \cdot q^-)}d\omega \ dq,
    \end{equation}
    with
    \begin{equation}
        \lambda = 2/(3+\sqrt{2}), \ \ \ q^+ = \frac{1}{2}(q+|q|\sigma), \ \ \ q^-=\frac{1}{2}(q-|q|\sigma), \ \ \ q = v - v_*.
    \end{equation}
    We would like to remark that the naive evaluation of \eqref{eqQk} requires $\mathcal{O}(n^2)$ operations, where $n=N^d$, which means that this spectral method is computationally expensive. The computation can be sped up by reducing the number of operations needed to evaluate the collision integral, introducing a FFT-based convolution approach in our spectral method. It is possible to find a new representation, called \emph{Carleman-like representation}, of the collision operator which conserves more symmetries when approximated on a bounded domain.
    Introducing the change of variables $x=r\sigma/2$ and $y = v_* - v - x$ the collision operator can be rewritten as \cite{MoPa:2006}
    \begin{equation}
        \Q(f,f)(v)=\int_{x\in \RR^d}\int_{y\in \RR^d}\tilde{B}(x,y)\delta(x\cdot y)\left[f(v+y)f(v+x)-f(v+x+y)f(v)\right]dx \ dy,
    \end{equation}
    where
    \begin{equation}\label{BetaTildeCarleman}
        \tilde{B}(x,y) = \tilde{B}(|x|, |y|) = 2^{d-1}B\left(\sqrt{|x|^2+|y|^2},\frac{|x|}{\sqrt{|x|^2+|y|^2}}  \right)(|x|^2+|y|^2)^{-\frac{d-2}{2}},
    \end{equation}
    where the last expression of $\tilde{B}(x,y)$ is valid on the manifold defined by $x\cdot y = 0$.
    Assuming the validity of the \emph{variable hard-sphere} (VHS) model, 
    \begin{equation}\label{eqVHS}
        B(|v-v_*|, \cos{(\theta)})=C_\gamma |v-v_*|^\gamma ,
    \end{equation}
    with $\gamma \in [0,1]$ and a constant $C_\gamma >0$, we obtain a constant $\tilde{B}(x,y)$ if $\gamma=0$ and $d=2$, or if $\gamma=1$ and $d=3$.
    
    Now let us discuss two possible truncation strategies \cite{rey-2023}. We start by considering a bounded domain $\mathcal{D}_T=\left[-T,T\right]^d \ (0<T<+\infty)$. A first approach is removing the collisions involving velocities outside of the box. This is the natural first stage for deriving conservative schemes based on the discretization of the velocity. 
    {A truncation on the modulus of $x$ and $y$ is not necessary in this case, since we are imposing them to be in the box.}
    So we obtain 
    \begin{multline}
        \Q^{tr}(f,f)(v) = \int \int_{\{x,\ y\in \ \RR^d \ | \ v+x,\ v+y,\ v+x+y\ \in \ \mathcal{D}_T\}}\tilde{B}(x,y)\delta(x\cdot y)\\
        \left[f(v+y)f(v+x)-f(v+x+y)f(v)\right]dx \ dy,
    \end{multline}
    where $v\in \mathcal{D}_T$.\\
    This last operator satisfies the following weak form
    \begin{multline}
        \int \Q^{tr}(f,f)(v) \zeta(v) dv =
        \frac{1}{4} \int \int \int_{\{v,\ x,\ y\in \ \RR^d \ | v,\ v+x,\ v+y,\ v+x+y\ \in \ \mathcal{D}_T\}}\tilde{B}(x,y)\delta(x\cdot y)\\
        f(v+x+y)f(v)\left[\zeta(v+y)+\zeta(v+x)-\zeta(v+x+y)-\zeta(v)\right] dv \ dx \ dy,
    \end{multline}
    from which we can deduce the conservation of mass, momentum, energy and the validity of the H-theorem for the entropy. In this derivation the collision kernel acquires artificial dependencies on $v, v_*, v', v_*'$, and the convolution-like structure is lost. This is the drawback of this kind of truncation.

    A second possible approach is truncating the integration with respect to variables $x$ and $y$ by setting them to be in $B_0(R)$. In order to obtain all possible collisions we consider $R=2S$ for a compactly supported function $f$ whose support is $B_0(S)$.\\
    Our goal is to implement the FFT algorithm to evaluate the quadrature approximation, so we will make use of periodic distribution functions, and then we have to properly take into account the aliasing effect due to periods superposition in the Fourier space. A geometrical argument shows that using the periodicity of the function it is enough to take $T\geq (3+\sqrt{2})S/2$ to prevent intersections of the region where $f$ is different from zero. The operator now can be written as 
    \begin{equation}\label{eqQRapSpe}
        \mathcal{Q}^R(f,f)(v)=
        \int_{x\in \mathcal{B}_0(R)}
        \int_{y\in \mathcal{B}_0(R)}
        \tilde{B}(x,y)\delta(x\cdot y)
        \left[f(v+y)f(v+x)-f(v+x+y)f(v)\right] \ dx \ dy,
    \end{equation}
    for $v\in \mathcal{D}_T$. The main interest of this representation is preserving the real collision kernel and its properties.

    The main idea to develop a more efficient algorithm is to use the previous representation, where taking $\zeta(v)=e^{ik\cdot v}$ we obtain the following spectral quadrature formula
    \begin{equation}
        \hat{Q}_k = \sum_{\substack{l,m=-N\\l+m=k}}^N \hat{f}_l \hat{f}_m \hat{\beta}_F(l, m), \ \ \ k = -N, ..., N,
    \end{equation}
    where $\hat{\beta}_F (l,m)=\hat{B}_F(l,m) - \hat{B}_F(m,m)$ are now given by
    \begin{equation}
        \hat{B}_F(l,m) =
        \int_{\mathcal{B}_0(R)}
        \int_{\mathcal{B}_0(R)}
        \tilde{B}(x,y)\delta(x\cdot y)
        e^{i(l\cdot x + m \cdot y)} \ dx \ dy.
    \end{equation}
    Finding a convolution structure in the previous equation  reduces the number of operations needed to evaluate it. In particular, we want to approximate each $\hat{\beta}_F(l,m)$ as follows
    \begin{equation}\label{BSpectral}
        \hat{\beta}_F(l,m) \simeq \sum_{p=1}^A \alpha_p(l)\alpha_p'(m),
    \end{equation}
    where $A$ represents the number of finite possible directions of collisions. In the end, we will obtain a sum of $A$ discrete convolutions, which can be computed in $\mathcal{O}(An\log_2{n})$ operations using the standard FFT technique.\\
    Now we make the decoupling assumption
    \begin{equation}
        \tilde{B}(x,y)=a(|x|)b(|y|),
    \end{equation}
    in order to obtain the convolution form. We would like to remark that this assumption is satisfied when $\tilde{B}$ is constant.\\
    This is the case of the Maxwellian molecules in dimension two, and hard-spheres in dimension three, the two most used cases. Indeed, putting the variable hard-sphere kernel inside the equation \eqref{BetaTildeCarleman} we obtain
    \begin{equation}
        \tilde{B}(x,y)=2^{d-1}C_\gamma (|x|^2 + |y|^2)^{-\frac{d-\gamma-2}{2}},
    \end{equation}
    and so we obtain that $\tilde{B}$ is constant if $d=2$, $\gamma=0$ (Maxwellian molecules in dimension two) or if $d=3$, $\gamma=1$ (hard-sphere molecules in dimension three).

    \paragraph{The 2D case}
    We start with Maxwellian molecules in dimension two, assuming, without loss of generality, $\tilde{B}=1$, since a different constant simply multiplies the operator. We move to polar coordinates $x=\rho e$ and $y=\rho' e'$ inside \eqref{BSpectral}, and we get
    \begin{equation}
        \hat{B}_F(l,m)=\frac{1}{4}\int_{\SSS^1}\int_{\SSS^1}
        \delta(e\cdot e')
        \left[\int_{-R}^Re^{i\rho(l\cdot e)}d\rho\right]
        \left[\int_{-R}^Re^{i\rho'(m\cdot e')}d\rho'\right]
        de \ de'.
    \end{equation}
    Now we introduce $\phi_R^2(s)=\int_{-R}^R e^{i\rho s} d\rho$, for $s\in \RR$, and we obtain the explicit formula
    \begin{equation}
        \phi_R^2(s)=2R\sinc{(Rs)},
    \end{equation}
    where $\sinc(x)=(\sin(x)/x)$.
    Using the parity properties of the $\sinc$ function inside the expression of $\hat{B}_F(l,m)$ we get 
    \begin{equation}
        \hat{B}_F(l,m)=
        \int_0^\pi \phi_R^2(l\cdot e_\theta) 
        \phi_R^2(m\cdot e_{\theta + \pi/2}) \ d\theta.
    \end{equation}
    Finally a regular discretization of $A$ equally spaced points $\theta_p=\pi p/A$ of the half circle, which is spectrally accurate \cite{can-hus-qua-zan-1991} because of the periodicity of the function, gives
    \begin{equation}
        \hat{B}_F(l,m)=\frac{\pi}{A}
        \sum_{p=1}^A \alpha_p(l)\alpha_p'(m),
    \end{equation}
    with
    \begin{equation}
        \alpha_p(l)=\phi_R^2(l\cdot e_{\theta_p}), \ \ \ \ \ \ \ \alpha_p'(m)=\phi_R^2(m\cdot e_{\theta_{p}+\pi/2}).
    \end{equation}

    \subsection{Conservation of macroscopic quantities}\label{secConMas}
    It is well known that the discretization of the velocity space on discrete grids introduces a loss of conservation \cite{mie-2000, dim-par-2014}. A natural solution to this problem is to consider a velocity support big enough such that the energy outside the domain is negligible. However, the truncation of the domain is not the only issue, as quadrature errors on the grid must also be taken into consideration. Therefore, it is highly desirable to construct a solver where conservation properties are exactly maintained during the simulation.\\    
    Several approaches have been proposed in the literature, exploiting either entropy minimization (as in \cite{mie-2000}) or $L^2$-minimization (as in \cite{gam-tha-2009, bos-cho-rus-2022}). The \emph{entropy approach} allows the construction of a conservative discrete Maxwellian, while $L^2$-minimization can be applied to more general distribution functions. In this paper, we implement the $L^2$-minimization approach proposed in \cite{gam-tha-2009}. Specifically, we apply this projection technique to the Anisotropic Maxwellian $\G[f]$ (see Equation \eqref{eqOpES-BGK}) involved in the AP-penalization technique (see Section \ref{sub:AP_Penalization_GSA}), ensuring that it exactly reproduces the desired macroscopic moments on the discrete grid. Furthermore, we apply the same correction to the full Boltzmann spectral collision operator (see Section \ref{secSpeMet}) in order to strictly enforce the exact conservation of macroscopic quantities at the collisional level. The same projection is also applied to the Maxwellian distributions used in the fluid-to-kinetic lifting (Section~\ref{sub:coupling}). Finally, it is worth noting that a similar weighted projection technique \cite{bos-cho-rus-2022} has been recently applied to a gas mixtures setting in \cite{cap-par-rey-2025-2}.

\section{Numerical simulations}\label{secNumSim}
    In this section we numerically test and investigate the implemented numerical scheme on complex geometries. The code is massively parallelized on GPU in order to speed up the code execution. The numerical simulations are executed on a shared workstation equipped with two \emph{Intel(R) Xeon(R) Gold 5418Y} CPUs, 128 GB of RAM and an \emph{NVIDIA RTX A1000} GPU, kindly provided by the Department of Mathematics and Computer Science of the University of Ferrara. Since other workloads may have been running on the same machine, the reported execution times may be affected by them and should be regarded as indicative. The implementation does not require specialized hardware, and thanks to the extensive exploitation of GPU parallelism and the highly optimized code design, the software can also be executed efficiently on conventional laptop-class systems. All simulations have been executed using the third-order in space and time solvers discussed above.\\
    The time step $\Delta t$ employed in the simulations is generally determined by the relation:
    \[
    \Delta t = 0.5 \frac{\Delta x_{\min}}{L},
    \]
    where $\Delta x_{\min}$ is the minimum edge length across all Delaunay triangles of the fluid cells in the spatial mesh, and $L$ is the bound of the velocity domain $\mathcal{V}=[-L,L]^2$.
    This criterion was adopted for all simulations with the following exceptions: in \hyperref[sec:test4]{Test 4} the safety factor is equal to 0.8 and in \hyperref[sec:test3]{Test 3} it is approximately 0.6, while in the temporal convergence analysis (\hyperref[sec:test1]{Test 1}) the time step was varied as described in Section~\ref{sec:test1}.

    \subsection{Validation}
    In this subsection some preliminary tests are executed in order to validate the implemented numerical scheme.
    \subsubsection{Test 1: Time convergence test}
    \label{sec:test1}
\begin{figure}[htbp]
        \centering
        \includegraphics[width=0.5\linewidth]{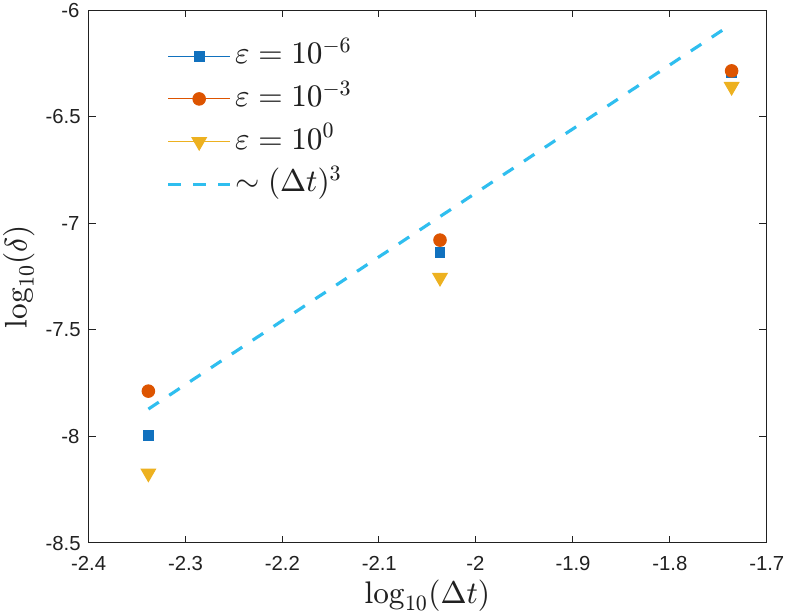}
        \caption{Temporal convergence rates of the BPR(3,5,3) scheme for the benchmark described in \hyperref[sec:test1]{Test 1}. The numerical error $\delta$ ($L^1$ norm) is plotted against the time step $\Delta t$ for $\varepsilon = 1, 10^{-3}, 10^{-6}$ under a frozen random regime decomposition, confirming the theoretical third-order accuracy across all solvers and highlighting the potential order reduction in the intermediate regime.}
        \label{Fig1}
    \end{figure}

    The primary objective of this benchmark is to demonstrate the third-order temporal convergence of the BPR(3,5,3) scheme and to test the coupling procedure under a fully multi-model setting. This test allows us to verify at the same time the temporal integration of the three solvers (Euler, Navier--Stokes, and full Boltzmann) and the accuracy of the interface coupling strategy detailed in Section \ref{sec:UnifiedFramework}. To this end, the regime indicator of each computational cell is randomly initialized to one of the three physical models (Euler, Navier--Stokes, or kinetic) and is kept fixed throughout the time evolution. The same random partition is used for all time steps $\Delta t$ and for the reference solution. As in \cite{cap-ten-2025}, keeping the decomposition fixed excludes the errors due to the regime adaptation, so that any loss of accuracy can only come from the coupling at the interfaces. We consider a physical domain $[0,2]\times[0,2]$ discretized with approximately 430 triangular cells. The Euler and Navier--Stokes solutions are initialized with the following macroscopic moments $W=(\rho, u_x, u_y, T)$:
    \begin{equation}
        W(x,y) = (1+0.1\exp{(-\lambda r^2)}\,,0\,,0\,,0.5),
    \end{equation}
    where $r^2=(x-x_c)^2+(y-y_c)^2$, $x_c=y_c=1$, and $\lambda=10$. The distribution function in the kinetic cells is initialized as a Maxwellian distribution matching the moments $W(x, y)$. The velocity space $[-7,7]^2$ is discretized using 32 Cartesian cells in both directions. The test is conducted across three distinct Knudsen numbers: $\varepsilon = 1$, $\varepsilon = 10^{-3}$, and $\varepsilon = 10^{-6}$. For each case, we systematically vary the time step $\Delta t$ and evaluate the numerical error at the final time $\delta$ for the density $\rho$ in the $L^1$ norm by comparing the numerical solution against a reference solution computed with the finest temporal resolution. The corresponding convergence rates are depicted in Figure \ref{Fig1}. These results confirm that third-order accuracy is maintained in both the highly rarefied regime ($\varepsilon = 1$) and the fluid limit ($\varepsilon = 10^{-6}$). This demonstrates that the scheme possesses \emph{Asymptotic Accuracy} (AA), as it recovers the third-order slope in the fluid limit ($\varepsilon \to 0$). However, the scheme does not exhibit \emph{Uniform Accuracy} (UA) across the entire range of Knudsen numbers. In the intermediate transition regime ($\varepsilon = 10^{-3}$), a small deviation from the optimal slope occurs, reducing the empirical convergence rate to approximately $2.5$. Since the same frozen decomposition and coupling preserve third-order convergence for $\varepsilon = 1$ and $\varepsilon = 10^{-6}$, this reduction is most likely not due to the interface coupling, but to the well-known order reduction phenomenon that affects IMEX Runge--Kutta schemes applied to kinetic equations in intermediate regimes.
    
    \subsubsection{Test 2: Space convergence test}\label{sec:test2}
    The objective of this second benchmark is to assess the spatial accuracy of the implemented numerical scheme. The interface coupling technique does not generally preserve high-order spatial convergence. Consequently, to correctly evaluate the spatial order of the scheme without the interference of the coupling error, this test is conducted exclusively in the kinetic regime, employing the full Boltzmann solver across the entire computational domain. We consider the same physical domain $[0,2]\times[0,2]$ and initial macroscopic moments $W=(\rho, u_x, u_y, T)$ as in \hyperref[sec:test1]{Test 1}:
    \begin{equation}
        W(x,y) = (1+0.1\exp{(-\lambda r^2)}\,,0\,,0\,,0.5),
    \end{equation}
    where $r^2=(x-x_c)^2+(y-y_c)^2$, $x_c=y_c=1$, and $\lambda=10$. The distribution function is correspondingly initialized as a Maxwellian, and the velocity space $[-7,7]^2$ is discretized using 32 Cartesian cells in both directions. To measure the spatial convergence rate, we systematically vary the spatial grid resolution. For each mesh size, the time step $\Delta t$ is suitably adjusted to strictly satisfy the hyperbolic Courant--Friedrichs--Lewy (CFL) stability condition at kinetic level. 

\begin{figure}[htbp]
        \centering
        \includegraphics[width=0.5\linewidth]{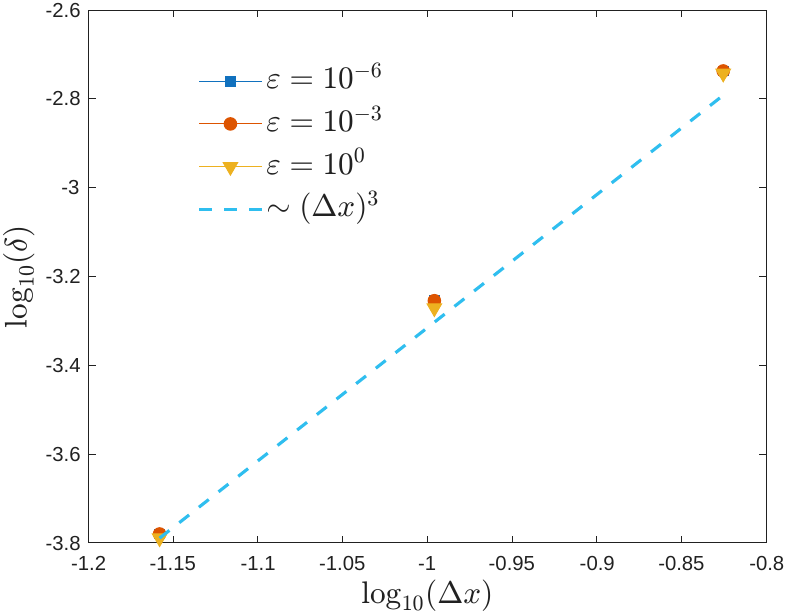}
        \caption{Spatial convergence rates for the benchmark described in \hyperref[sec:test2]{Test 2}. The numerical error $\delta$ ($L^1$ norm) is plotted as a function of the spatial mesh size for $\varepsilon = 1, 10^{-3}, 10^{-6}$ in a purely kinetic regime. The curves nearly overlap and show third-order spatial accuracy in all the regimes considered.}
        \label{Fig2}
    \end{figure}

    The test is performed across three distinct Knudsen numbers, specifically $\varepsilon = 1$, $\varepsilon = 10^{-3}$, and $\varepsilon = 10^{-6}$. The numerical error $\delta$, computed in the $L^1$ norm by comparing the obtained results for the density $\rho$ against a reference solution calculated on the finest spatial grid, confirms that the expected third-order spatial convergence is fully maintained in the purely kinetic regime. The corresponding convergence rates are depicted in Figure \ref{Fig2}, where the numerical error $\delta$ is plotted against the characteristic mesh size $\Delta x$, defined in this context as the square root of the mean cell area.

    Furthermore, it can be observed that the error curves corresponding to the different Knudsen numbers nearly overlap, showing no degradation of the spatial accuracy in any of the considered regimes. This uniform behaviour is due to the fact that the spatial discretization error dominates by several orders of magnitude over the temporal one, as can be seen by comparing the error scales of Figure \ref{Fig2} with those of Figure \ref{Fig1}. Consequently, the third-order spatial accuracy is not affected by the temporal order reduction. This shows that, for the time steps dictated by the CFL condition, the spatial accuracy of the purely kinetic scheme is not affected by the stiffness of the collision term in the considered regimes.

    \subsubsection{Test 3: Numerical error convergence}\label{sec:test3}
    In this test, we investigate the sensitivity of the coupled numerical scheme to the user-defined thresholds that govern the dynamic domain identification algorithm. The objective is to analyze the trade-off between computational efficiency and numerical accuracy under different physical regimes.

    To this end, we consider a physical domain $[0,2]\times[0,2]$ discretized with approximately 430 cells. The distribution function $f(x, y)$ is initialized as a sum of two Maxwellians $f_P(x,y)$ on a static and uniform background $f_B(x,y)$ according to
    \begin{equation}
        \begin{cases}
            f(x, y) = f_B(x,y)+f_P(x,y),\\
            f_B(x,y) = \mathcal{M}_{1, 0, 0, 0.2},\\
            f_P(x, y) = 0.5\left(\mathcal{M}_{\rho_p, 0.5, 0.5, 0.2}+\mathcal{M}_{\rho_p, -0.5, -0.5, 0.2}\right),\\          \rho_p(x,y)=0.1\exp{(-\lambda r^2)},
        \end{cases}
    \end{equation}
    where $\mathcal{M}_{\rho, u_x, u_y, T}$ is the Maxwellian initialized with the density $\rho$, velocities $u_x$, $u_y$, temperature $T$ and $r^2=(x-x_c)^2+(y-y_c)^2$, $x_c=y_c=1$, and $\lambda=20$.
    
    The macroscopic moments of the Navier--Stokes and Euler equations $W=(\rho, u_x, u_y, T)$ are initialized everywhere by computing the moments of the distribution function. The velocity space $[-4,4]^2$ is discretized using 32 Cartesian cells in both directions.
    
    At $t=0$ all cells are assigned to the kinetic model, and the first time step is performed with the full Boltzmann solver in the whole domain; the dynamic domain decomposition is activated from the second time step on. This is needed because the regime indicators rely only on macroscopic quantities and cannot detect the non-equilibrium part of the initial datum, which is concentrated in the central region where $f_P$ is not negligible (see Section~\ref{subBreakdownMain}).
    
    Starting from the second iteration, the test is conducted implementing the dynamic domain decomposition discussed in this paper. The test is executed for different values of the Knudsen number, and for each of them, we systematically vary the threshold parameters that dictate the assignment of each computational cell to the Euler, Navier--Stokes, or full Boltzmann solver. By tuning these thresholds, we effectively manipulate the dynamic spatial decomposition of the domain, which can be quantified by the ratio of cells assigned to the $i$-th regime (Euler, NS or Boltzmann) over the total number of cells ($N_i/N_{\text{tot}}$) throughout the simulation. For each configuration, we evaluate the numerical error by comparing the obtained results against a high-fidelity reference solution calculated using the full Boltzmann solver uniformly across the entire domain.

\begin{figure}
        \centering
        \includegraphics[width=1.0\linewidth]{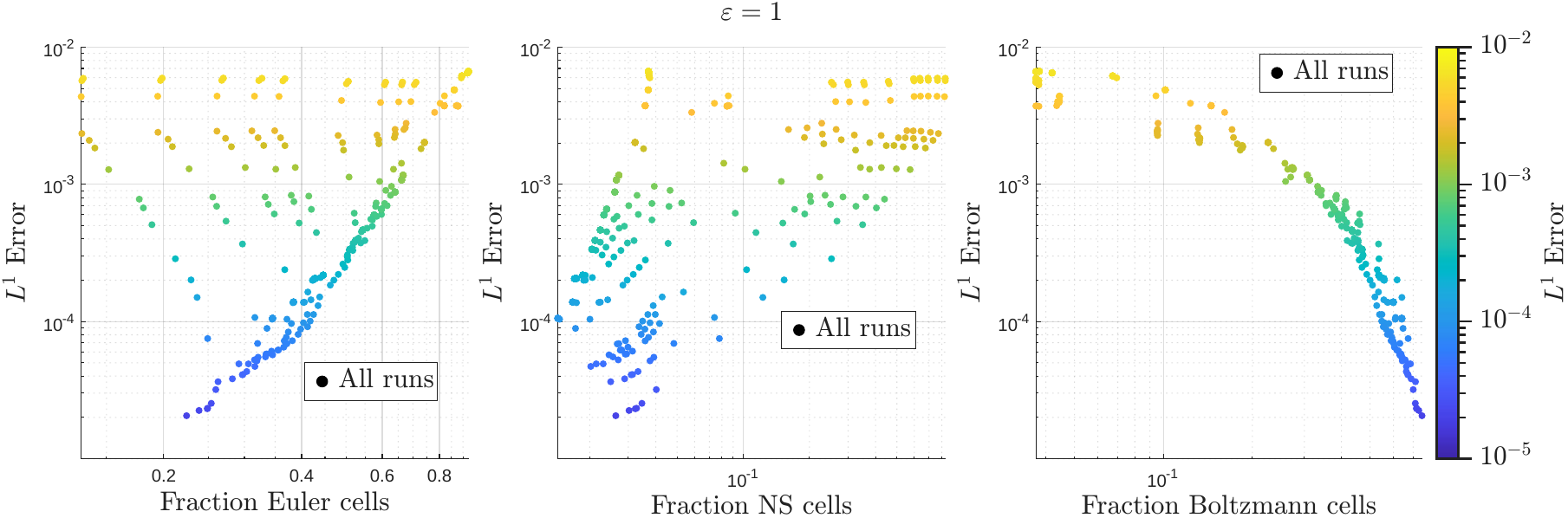}
        
        \smallskip
        
        \includegraphics[width=1.0\linewidth]{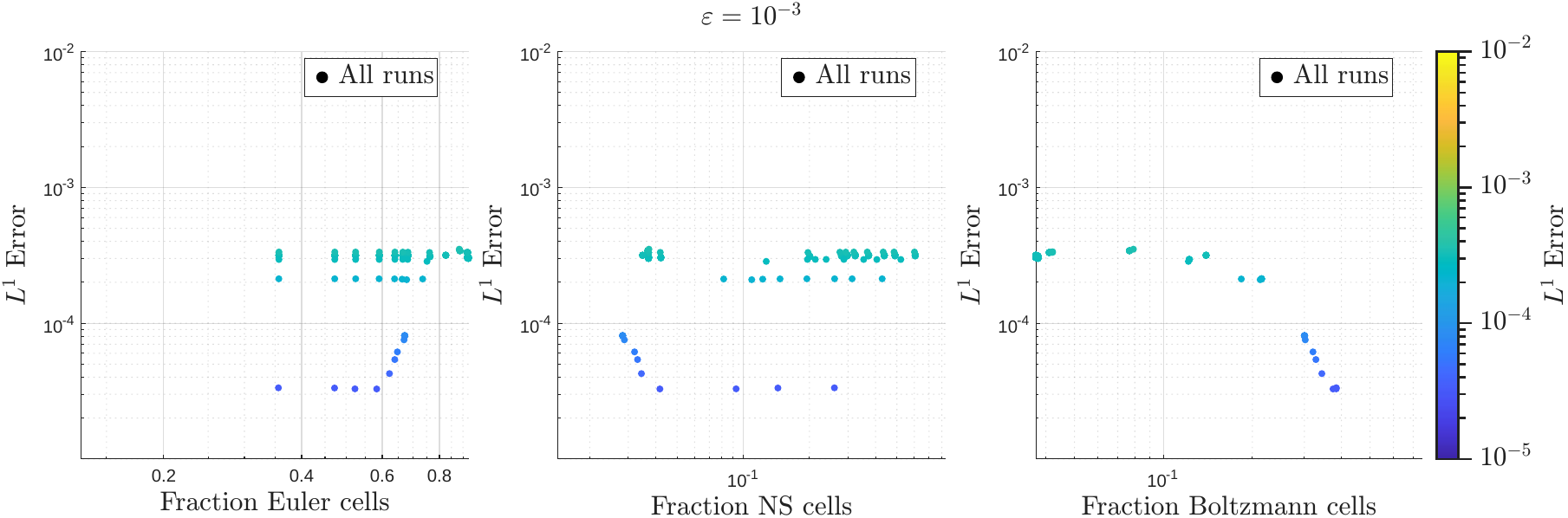}
        \caption{\hyperref[sec:test3]{Test 3}: One-dimensional projections of the global $L^1$ error norm for the density $\rho$ (evaluated at the final simulation time) as a function of the time-averaged cell fractions for the Euler (left subplots), Navier--Stokes (center subplots), and Boltzmann (right subplots) domains. The scatter points represent individual simulation runs; in each subplot the other domain fractions are not fixed. The figure highlights the  shift from the highly rarefied regime ($\varepsilon = 1$, top macro-panel) to the stiff fluid limit ($\varepsilon = 10^{-3}$, bottom macro-panel). The same axis limits, axis scalings and colorbar are used in both macro-panels to allow a direct comparison between the two regimes.}
        \label{Fig4}
    \end{figure}

    In Figure \ref{Fig4}, we evaluate the accuracy of the adaptive domain decomposition method by plotting the $L^1$ error of the density $\rho$ at the final time with respect to the reference full Boltzmann solution for two distinct physical regimes. Specifically, the upper panel illustrates the highly rarefied regime characterized by $\varepsilon = 1$, whereas the lower panel displays the stiff, collisional fluid regime with $\varepsilon = 10^{-3}$. Within each macro-panel, the data is organized into three distinct subplots to independently isolate the sensitivity of the error with respect to the fraction of Euler cells (left), Navier--Stokes cells (center), and Boltzmann cells (right). Since the number of Euler, Navier--Stokes, and Boltzmann cells changes at each iteration, these metrics are the time-averaged cell fractions computed over the entire simulation.

    In total, we present the results of 441 simulations for the rarefied regime and 289 simulations for the stiff regime. We use more configurations in the rarefied case because the points are more widely spread over the parameter space, so that more data are needed to identify the trend. In the stiff regime the points are more clustered, and fewer simulations are sufficient.

    In these visualizations, each scatter point represents a distinct threshold configuration (i.e., a different simulation). In each panel the error is plotted against one fraction only, while the other fractions are not fixed and change from point to point. The color of each scatter point is associated with a colorbar to better highlight the trend of the error. The same axis limits and colorbar are used for both regimes, to allow a direct comparison between them.

    In the highly rarefied regime ($\varepsilon = 1$, upper panel), the domain adaptation algorithm produces significant deviations from the reference solution when a dominant fraction of Euler or Navier--Stokes cells is employed. This behaviour is physically justified by the well-known failure of macroscopic hydrodynamic closures in capturing non-equilibrium kinetic phenomena. However, as the transition thresholds are systematically tightened (thereby reducing the Euler and Navier--Stokes regions), the error decreases by more than an order of magnitude, approaching the lowest errors observed in the fluid regime.

    Conversely, in the stiff regime ($\varepsilon = 10^{-3}$, lower panel), the overall error is consistently lower across the entire parameter space. In this collisional limit, the error shows an almost flat dependency over a wide range of Euler and Navier--Stokes cell fractions, indicating that the macroscopic hydrodynamic models can be widely utilized without a severe accuracy penalty. This confirms that the fluid closures reproduce the correct behaviour in the dense regime, which strongly reduces the accuracy penalty associated with fluid-kinetic transitions. As expected, for both values of $\varepsilon$, increasing the time-averaged fraction of Boltzmann cells drives the computed solution closer to the full Boltzmann reference.

    These trends show that, by calibrating the model transition thresholds, it is possible to control the trade-off between computational cost and numerical accuracy and to select the desired level of accuracy.

    \subsection{Numerical applications}
    In this section, a series of numerical experiments is presented to assess the performance of the proposed methodology. The test cases are designed to validate the coupling procedure, evaluate the capability of the scheme to handle complex geometries, and demonstrate the effectiveness of the dynamic domain decomposition strategy. In particular, the results highlight the ability of the proposed approach to accurately reproduce the behaviour of the full Boltzmann solution, outperforming both Euler and Navier--Stokes models in regimes where kinetic effects are significant. We also compare the computational cost of the hybrid method with that of the full Boltzmann solver.

    \subsubsection{Test 4: Sod shock tube test}\label{sec:test4}
    The first test case considers the classical Sod shock-tube problem in a strongly collisional regime corresponding to a very small Knudsen number $\varepsilon=10^{-5}$. Such a stiff regime does not constitute the most natural application scenario for a Boltzmann solver, since the flow is expected to be accurately described by the Euler equations and a purely macroscopic approach would be computationally more efficient. 
    Nevertheless, this test is useful for several checks. Each solver must reproduce the correct Sod solution, so the test validates the individual kinetic and continuum solvers. It also checks that the coupling does not introduce spurious artifacts at the kinetic--continuum interfaces, and it tests the AP property, since the solution must converge to the continuum limit without resolving the kinetic scales. Finally, it verifies that the adaptive algorithm reproduces the reference solution while reducing the kinetic region, and it provides a benchmark for the regime-detection criterion.    
    In particular, the adaptive procedure is expected to identify kinetic regions only in localized portions of the domain characterized by large gradients and stronger departures from thermodynamic equilibrium, while leaving the remaining cells in the continuum description. Although non-equilibrium effects are weak in this stiff regime, the ability of the algorithm to correctly localize the regions where they are expected to be most pronounced constitutes an important consistency check of the proposed adaptive methodology. Therefore, the objective of this test is not to demonstrate the necessity of a kinetic description for the Sod problem itself in the stiff limit, but rather to verify that the proposed adaptive framework correctly recovers the continuum limit while dynamically identifying the regions where a kinetic treatment would be most justified.

    The computational domain is the box $\mathcal{T} = [0, 1]\times[0, 0.01]$ discretized with a total number of cells approximately equal to 810. Periodic boundary conditions are implemented along the $y$-direction, and zero-gradient boundary conditions along the $x$-direction. The velocity space $\mathcal{V}=[-8,8]^2$ is discretized using 32 regular Cartesian control volumes in both directions. The macroscopic moments $W=(\rho, u_x, u_y, T)$ associated to each cell are
    \begin{equation}
        \begin{aligned}
            W_L &= (1,\,0,\,0,\,1),      & \qquad x<0.5,\\
            W_R &= (0.125,\,0,\,0,\,0.8), & \qquad x\ge 0.5.
        \end{aligned}
    \end{equation}
    The macroscopic variables evolved by the Euler and Navier--Stokes solvers are initialized according to the Riemann data given above. The kinetic distribution function is instead initialized throughout the computational domain as the Maxwellian distribution associated with the corresponding local macroscopic moments. Regarding the regime classification, all cells are initially assigned to the kinetic regime. This choice avoids introducing any a priori assumption on the spatial distribution of the flow regimes and provides a uniform initial state for the adaptive procedure. Starting from the first iteration, the dynamic domain decomposition algorithm automatically identifies the most appropriate physical description in each cell according to the regime-detection criterion. The solution is then evolved using the adaptive coupling strategy described in the previous section, allowing the computational domain to dynamically partition into kinetic, Euler, and Navier--Stokes regions as dictated by the local flow conditions.

\begin{figure}
        \centering
        \includegraphics[width=0.32\linewidth]{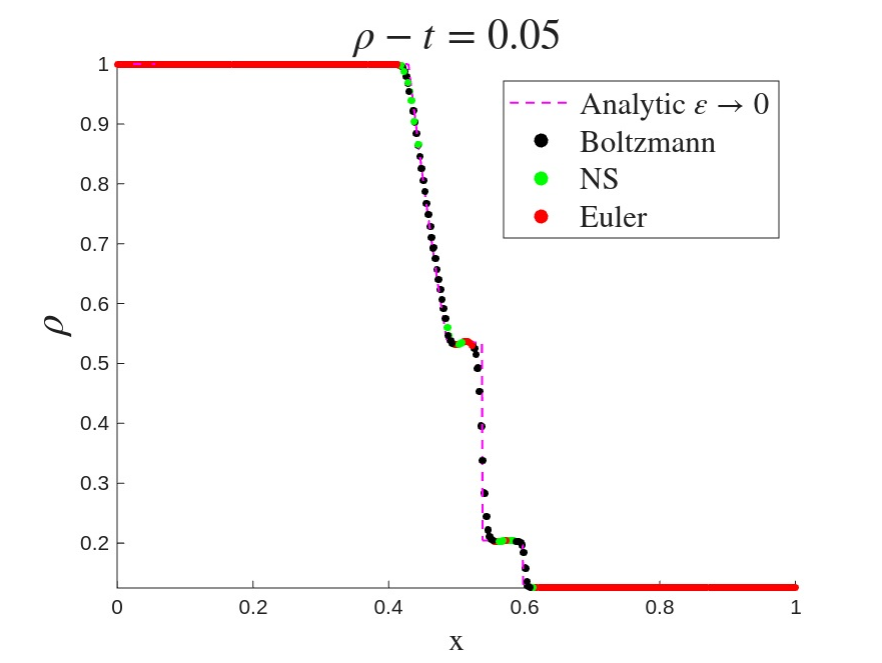}
        \includegraphics[width=0.32\linewidth]{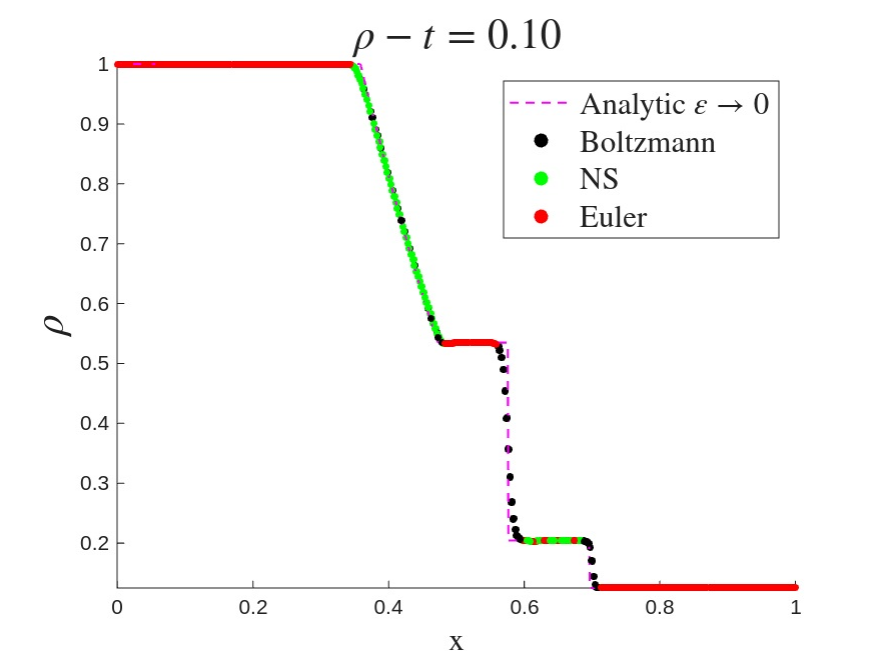}
        \includegraphics[width=0.32\linewidth]{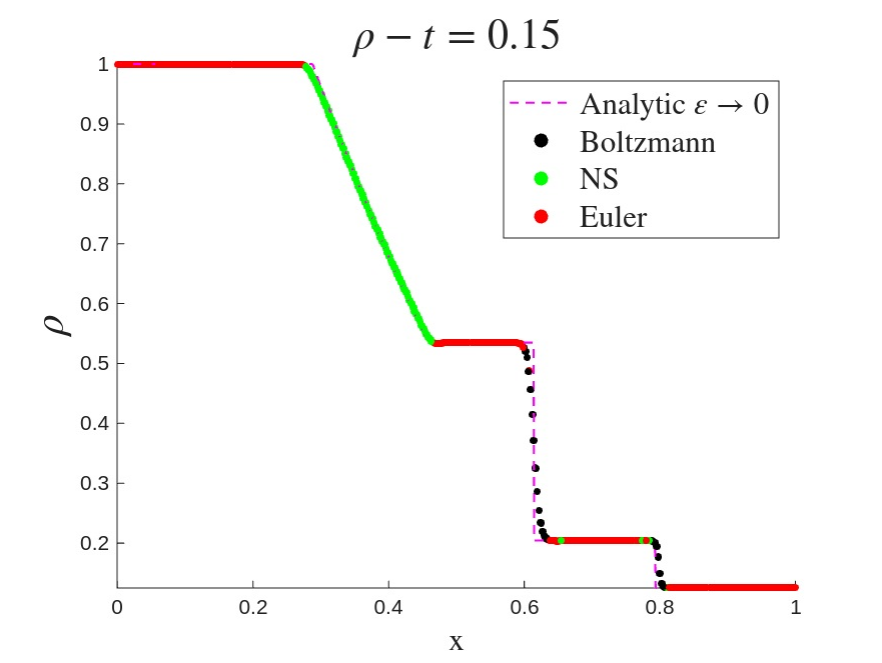}
        
        \smallskip

        \includegraphics[width=0.32\linewidth]{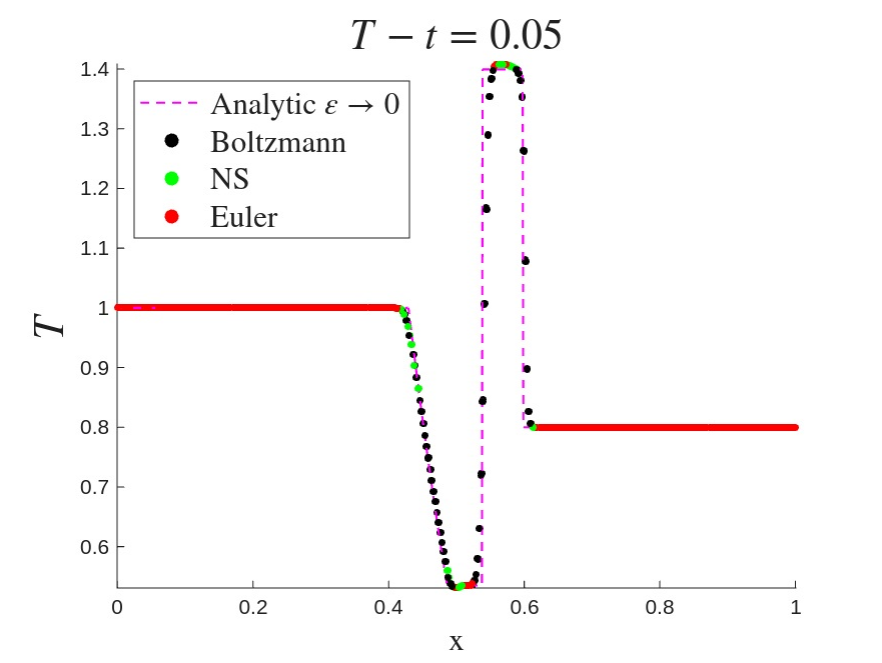}
        \includegraphics[width=0.32\linewidth]{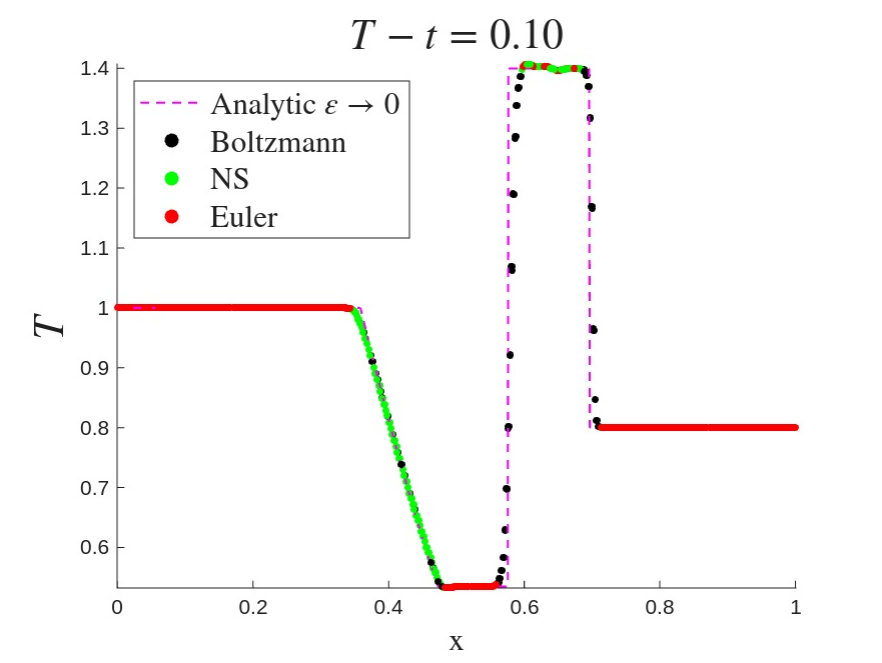}
        \includegraphics[width=0.32\linewidth]{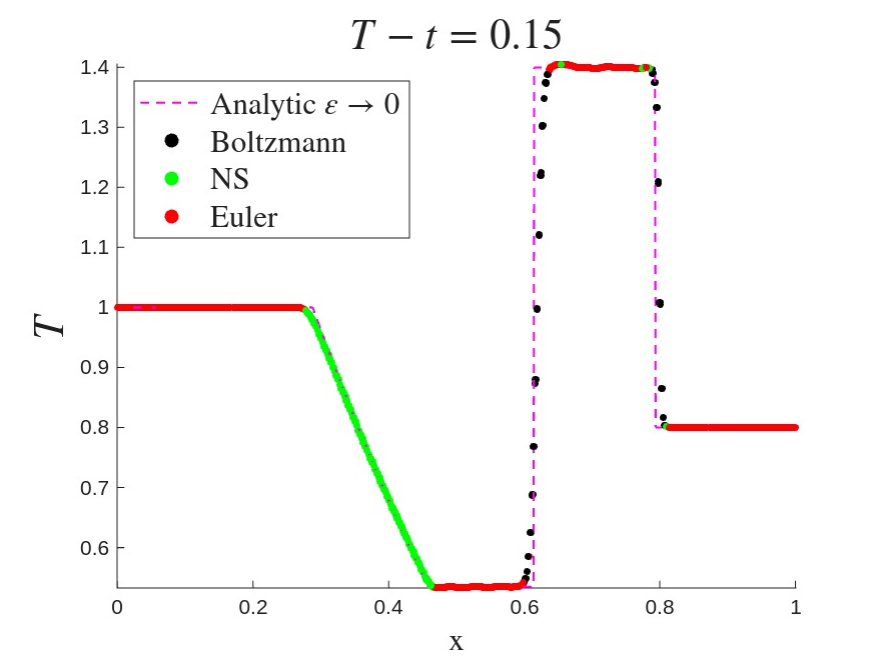}
        
        \smallskip

        \includegraphics[width=0.32\linewidth]{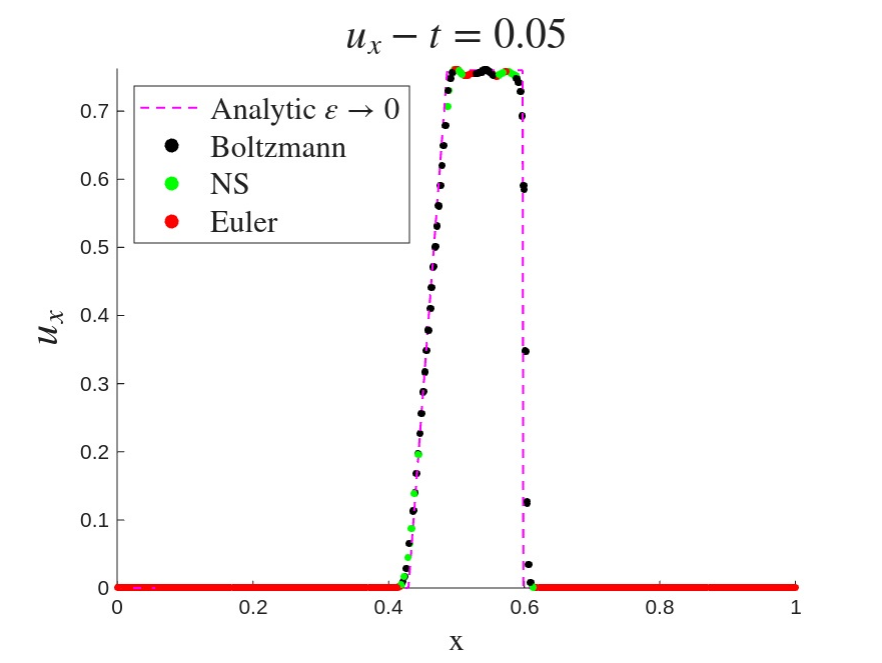}
        \includegraphics[width=0.32\linewidth]{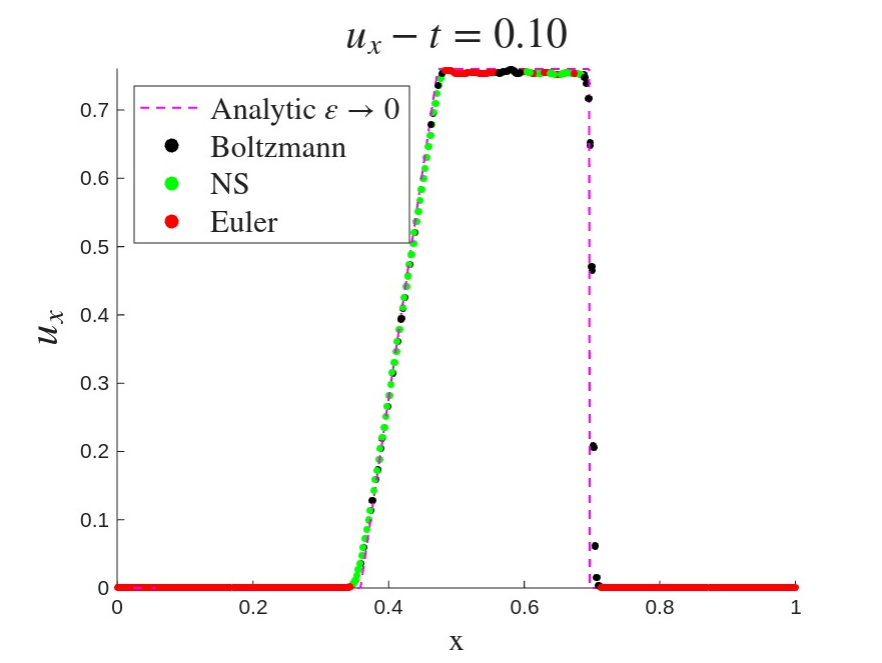}
        \includegraphics[width=0.32\linewidth]{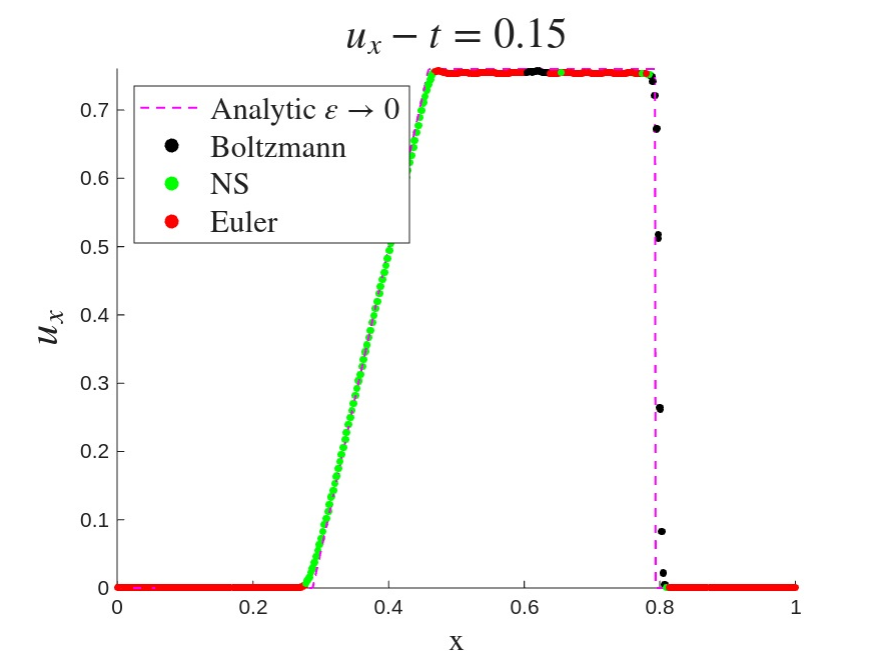}

        \smallskip

        \includegraphics[width=\linewidth]{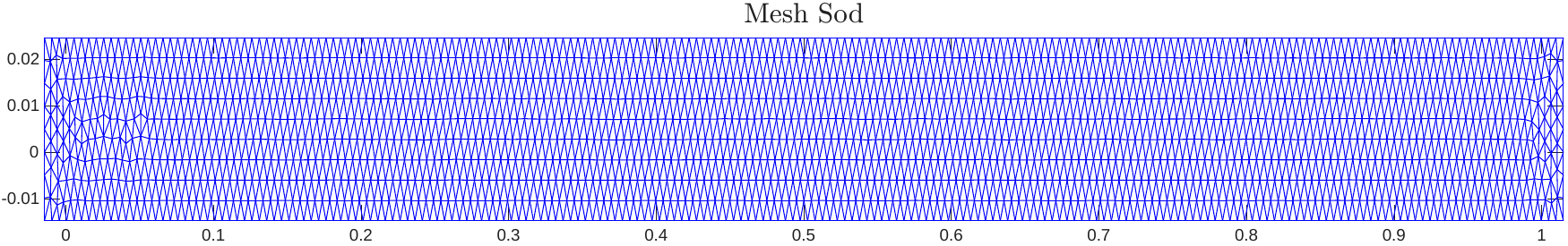}
        
        \caption{\hyperref[sec:test4]{Test 4}: Temporal evolution of the macroscopic variables for the Sod shock tube problem in the stiff fluid limit. The color coding indicates the dynamically assigned physical model: Boltzmann (black), Navier--Stokes (green), and Euler (red), while the exact analytical solution of the Euler Riemann problem is superimposed in magenta. Because of the small Knudsen number, the global solution is physically macroscopic; however, the gradient-based sensors correctly trigger the kinetic and NS solvers across steep discontinuities and mild gradients, respectively, identifying regions where non-equilibrium would typically emerge in rarefied conditions. The hybrid solution agrees with the exact solution of the Euler equations, which confirms the AP property and the correctness of the coupling. Over time, as waves propagate and new flat regions form, cells switch back to the Euler and Navier--Stokes models. The grid displays density ($\rho$, top row), temperature ($T$, second row), and $x$-velocity ($u_x$, third row) at $t = 0.05$ (left), $0.10$ (center), and $0.15$ (right). Finally, the bottom row displays the computational mesh, which extends slightly beyond the physical domain to accommodate the ghost cells used for the boundary conditions.}
        \label{Fig5}
    \end{figure}

    Figure \ref{Fig5} presents the temporal evolution of the macroscopic variables for the Sod shock tube problem. The graphical layout is structured as a grid: the top row displays the density ($\rho$), the second row shows the temperature ($T$), the third row reports the horizontal velocity ($u_x$), and the bottom row illustrates the computational mesh employed for the simulation. It is worth noting that the plotted grid extends slightly beyond the physical domain, as it explicitly includes the ghost cells required to enforce the boundary conditions. Employing an unstructured mesh for a fundamentally one-dimensional problem is a deliberate validation choice: it shows that the scheme preserves the one-dimensional structure of the solution, without spurious multidimensional artifacts, even when the element edges are not aligned with the main direction of the flow. For each macroscopic quantity, the panels from left to right correspond to the successive time instants $t = 0.05$, $0.10$, and $0.15$. The color coding applied to the spatial profiles explicitly identifies the physical model dynamically assigned to each computational cell: Boltzmann cells are depicted in black, Navier--Stokes (NS) cells in green, and Euler cells in red. Additionally, the exact analytical solution of the Riemann problem for the Euler equations is superimposed in magenta as a reference. Note that, in this configuration, the Knudsen number is small enough for the flow to remain in the continuum limit, so that the solution agrees with the macroscopic analytical solution. Consequently, the activation of the Boltzmann and NS solvers is not dictated by a severe physical departure from equilibrium, but rather by the gradient-based nature of the domain sensors. The kinetic solver is automatically triggered across the sharp discontinuities (i.e., the shock front and the contact wave), while the NS equations are activated across the mild gradients of the rarefaction fan. Essentially, the sensors identify those regions where non-equilibrium effects would naturally emerge if the gas were more rarefied. Testing the framework in this stiff limit is important to validate the AP property of the numerical scheme. 
    The results show that the adaptive strategy recovers the correct macroscopic limit, in close agreement with the analytical solution, and that the coupling between the kinetic and macroscopic models does not introduce visible artifacts. Moreover, the model selection works in both directions: as the waves propagate and leave behind flat equilibrium plateaus, the local gradients relax and the cells switch back to the Euler model.
    
    In the present numerical simulation, the threshold values for the switching criteria (Equations \eqref{crit:CompEuler_Main} and \eqref{criterioSecondoSwitch}) are set to $\eta_1 = \ve/2$ and $\eta_2 = \ve/100$.\\
    In terms of computational efficiency, the execution times on GPU are: $\text{Duration}_{\text{Euler}}=191.9\,\text{s}$, $\text{Duration}_{\text{NS}}=186.3\,\text{s}$, $\text{Duration}_{\text{Hybrid}}=370.3\,\text{s}$, and $\text{Duration}_{\text{Boltzmann}}=1053\,\text{s}$. The adaptive framework provides a computational speedup of approximately $2.8$ compared to the full Boltzmann solver. While the speedup indicator is inherently problem-dependent (being influenced by threshold selection, spatial and velocity grid refinement, and initial conditions), it shows that the proposed dynamic domain decomposition reduces the computational cost. A comprehensive discussion on computational efficiency, including the corresponding CPU and GPU execution times and a detailed analysis of the hybrid speedup, is provided in Section \ref{secNumEff}.

    \subsubsection{Test 5: Flow over a disk}\label{sec:test5}
    The second test case considers a rarefied flow past a stationary circular disk, a configuration characterized by strong non-equilibrium effects both near the solid boundary and in the wake region. This setup provides a challenging benchmark to assess the accuracy and efficiency of the proposed adaptive multiscale framework against a full Boltzmann reference solution. In the adaptive strategy, the Euler model is employed in regions close to thermodynamic equilibrium, the Navier--Stokes (NS) equations are activated in moderately rarefied zones, and the full Boltzmann description is retained only in strongly non-equilibrium regions where continuum assumptions fail.

\begin{figure}[t!]
        \centering
        \includegraphics[width=0.48\linewidth]{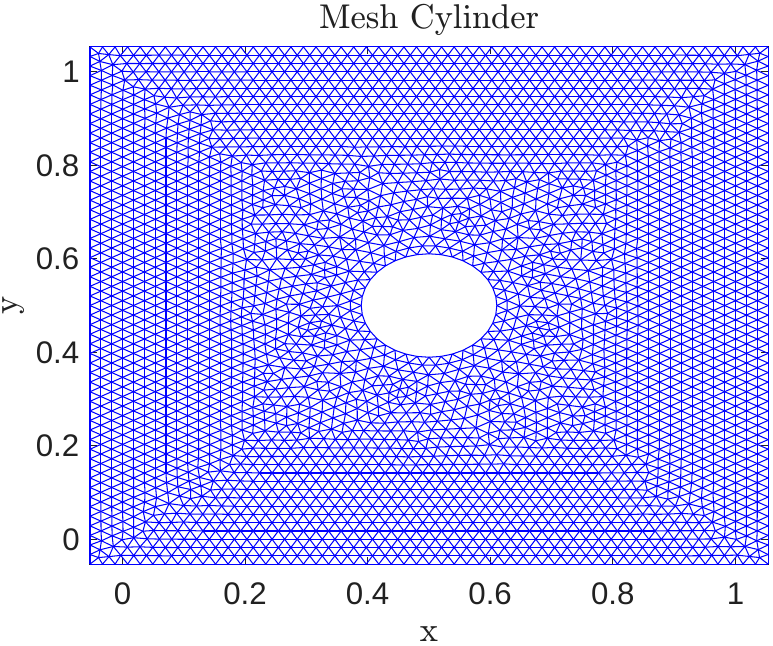}
        \caption{\hyperref[sec:test5]{Test 5}: Computational mesh layout used for the spatial discretization, including the ghost cells at the boundaries.}
        \label{Fig7}
    \end{figure}

    The computational domain is $\mathcal{T} = [0, 1]\times[0, 1]$, discretized with approximately $5100$ cells. The corresponding mesh layout, which includes the ghost cells used for the spatial discretization, is depicted in Figure~\ref{Fig7}. Zero-gradient boundary conditions are applied at the external boundaries of the domain, while slip conditions are enforced on the disk surface, which has a radius of $R = 0.11$ and is centered at $(x_c, y_c)=(0.5, 0.5)$. The velocity space $\mathcal{V}=[-7,7]^2$ is discretized using a regular $24 \times 24$ Cartesian grid. The macroscopic fields are initialized with the uniform state $W=(\rho, u_x, u_y, T) = (1,\,0.8,\,0,\,1)$, and the initial kinetic distribution function is set to the corresponding local Maxwellian state. 

    To avoid any a priori bias on the spatial distribution of the flow regimes, all cells are initially assigned to the kinetic description. From the first time step onward, the dynamic domain decomposition algorithm automatically partitions the domain into kinetic, Navier--Stokes, and Euler regions according to the local flow conditions and the regime-detection criterion.

    The numerical simulation is executed in the rarefied regime with $\varepsilon=0.1$. To evaluate the accuracy of the coupling strategy, we track the space-time averaged $L^1$ error for a representative macroscopic variable, namely the temperature $T$. 
    This metric is computed by integrating the absolute temperature difference relative to the full Boltzmann reference solution over the fluid cells, normalized by their total area and by the number of stored time instants. Mathematically, the error for a given solver is defined as:
    \begin{equation}\label{eqL1_1}
        \text{Error}_{b} = \frac{1}{N_t \sum_{i\in\mathcal{F}} |P_i|}
        \sum_{n=1}^{N_t} \sum_{i\in\mathcal{F}} \left| T_{b,i}^n - T_{\text{ref}, i}^n \right| |P_i|,
    \end{equation}
    where $b \in \{\text{Euler}, \text{NS}, \text{Hybrid}\}$ denotes the specific solver under consideration, $\mathcal{F}$ is the set of fluid cells, $N_t$ is the number of stored time instants, and $|P_i|$ is the area of the cell $P_i$. The term $T_{b, i}^n$ represents the temperature computed by solver $b$ at cell $i$ and time level $n$, while $T_{\text{ref}, i}^n$ is the corresponding reference temperature obtained from the full Boltzmann simulation.
    
    The resulting average temperature errors are $\text{Error}_{\text{Hybrid}} = 4.6\times 10^{-3}$, $\text{Error}_{\text{NS}} = 9.3\times 10^{-3}$, and $\text{Error}_{\text{Euler}} = 13\times 10^{-3}$. These results confirm that the hybrid solver significantly improves upon the macroscopic models, being about twice as accurate as the full Navier--Stokes model. \\    
    Regarding computational efficiency, the execution times on the GPU are: $\text{Duration}_{\text{Euler}} = 122.7\,\text{s}$, $\text{Duration}_{\text{NS}} = 121.0\,\text{s}$, $\text{Duration}_{\text{Hybrid}} = 247.3\,\text{s}$, and $\text{Duration}_{\text{Boltzmann}} = 621.0\,\text{s}$. The adaptive framework provides a computational speedup of approximately $2.5$ compared to the full Boltzmann solver.\\
    As for the Sod test, these indicators depend on the problem and on the numerical parameters; the CPU and GPU execution times and a detailed analysis of the hybrid speedup are given in Section~\ref{secNumEff}.
    In the present numerical simulation, the threshold values for the switching criteria (Equations \eqref{crit:CompEuler_Main} and \eqref{criterioSecondoSwitch}) are set to $\eta_1 = \ve/4$ and $\eta_2 = 10\ve$.

\begin{figure}
        \centering
        \includegraphics[width=0.32\linewidth]{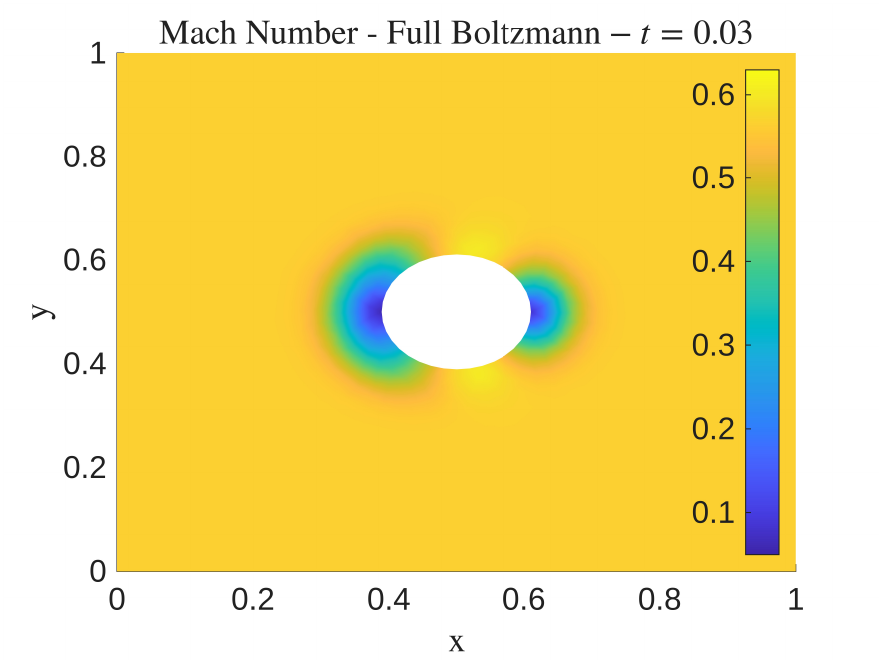}
        \includegraphics[width=0.32\linewidth]{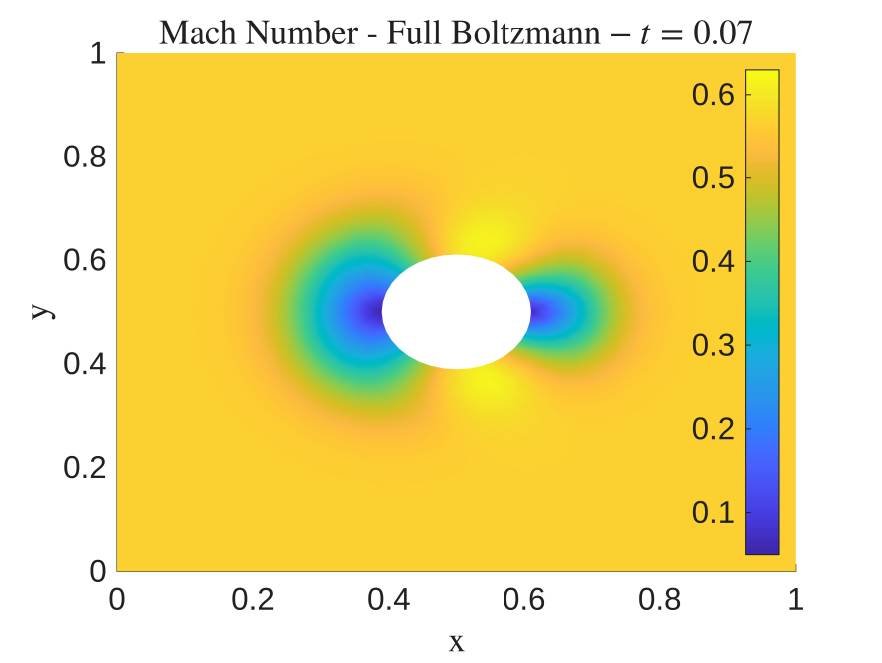}
        \includegraphics[width=0.32\linewidth]{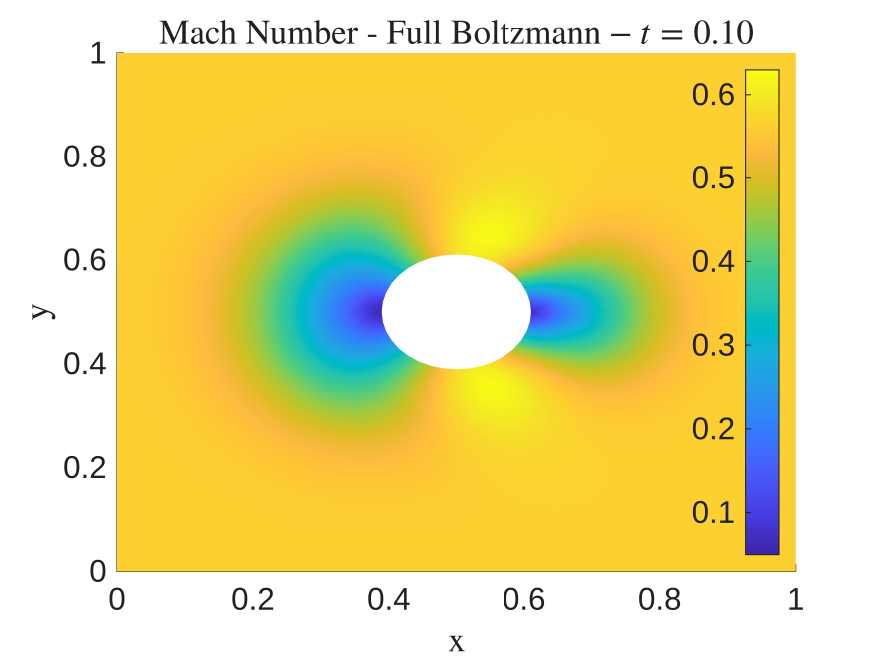}
        
        \smallskip
        
        \includegraphics[width=0.32\linewidth]{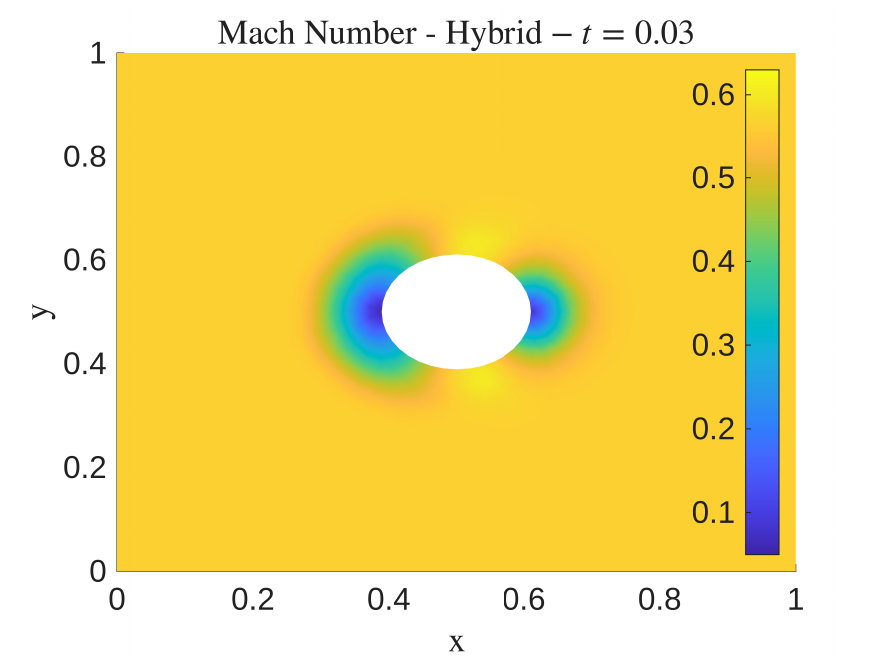}
        \includegraphics[width=0.32\linewidth]{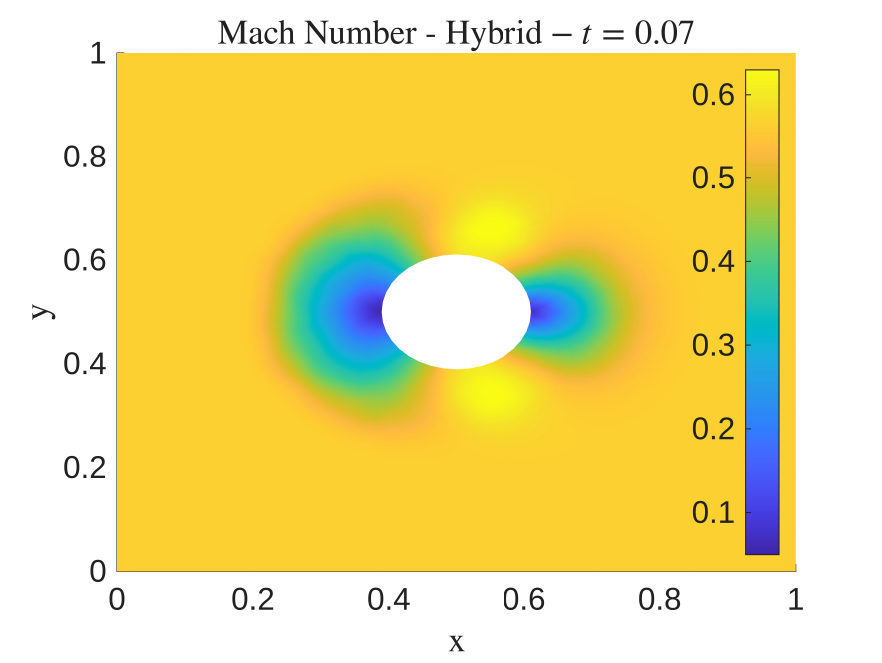}
        \includegraphics[width=0.32\linewidth]{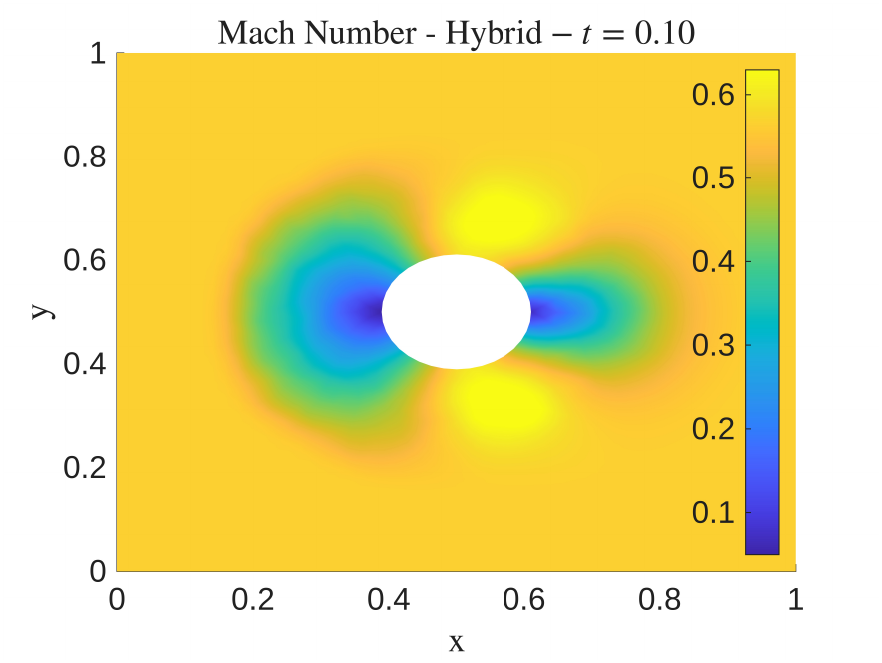}
        
        \smallskip
    
        \includegraphics[width=0.32\linewidth]{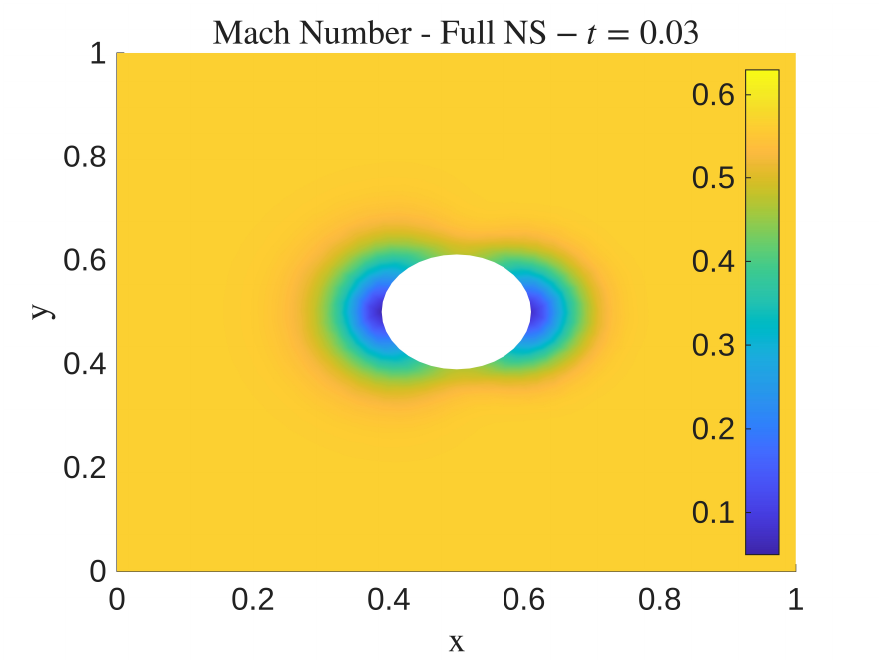}
        \includegraphics[width=0.32\linewidth]{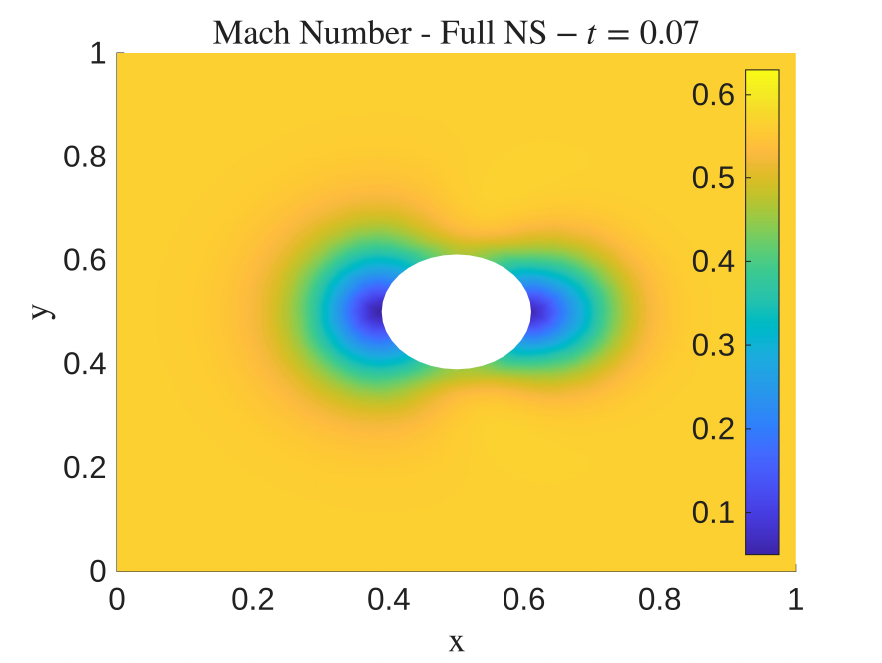}
        \includegraphics[width=0.32\linewidth]{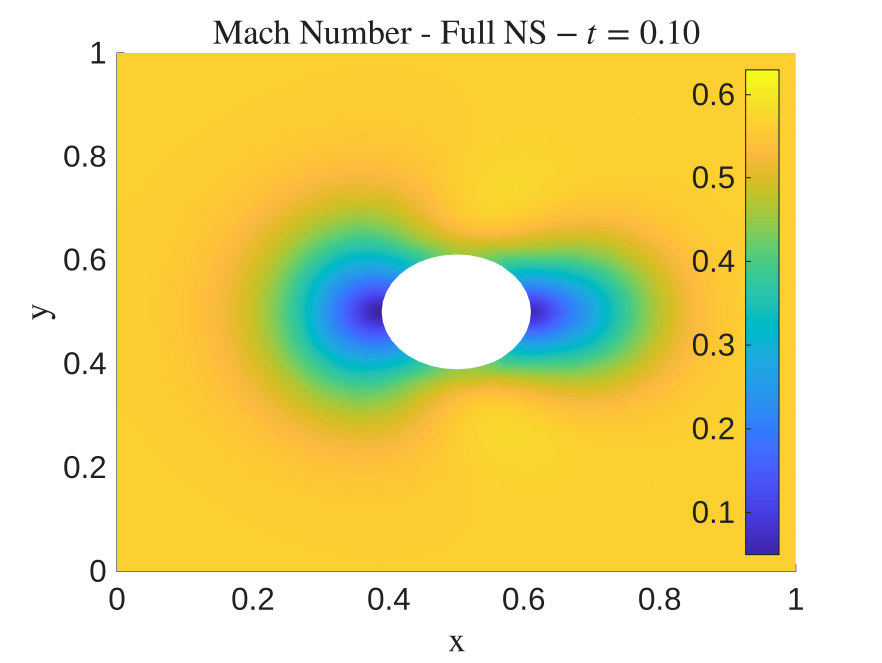}
    
        \smallskip
    
        \includegraphics[width=0.32\linewidth]{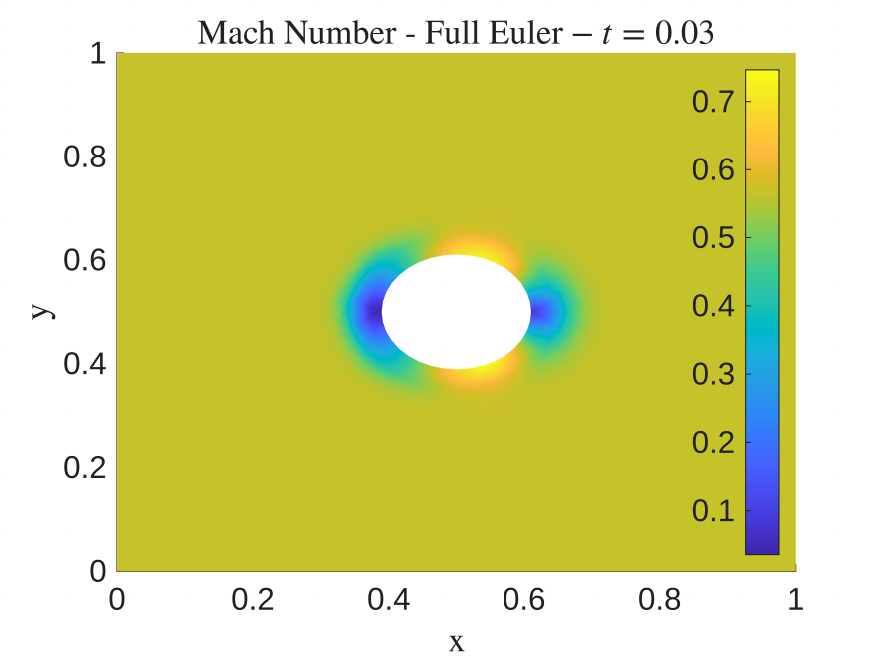}
        \includegraphics[width=0.32\linewidth]{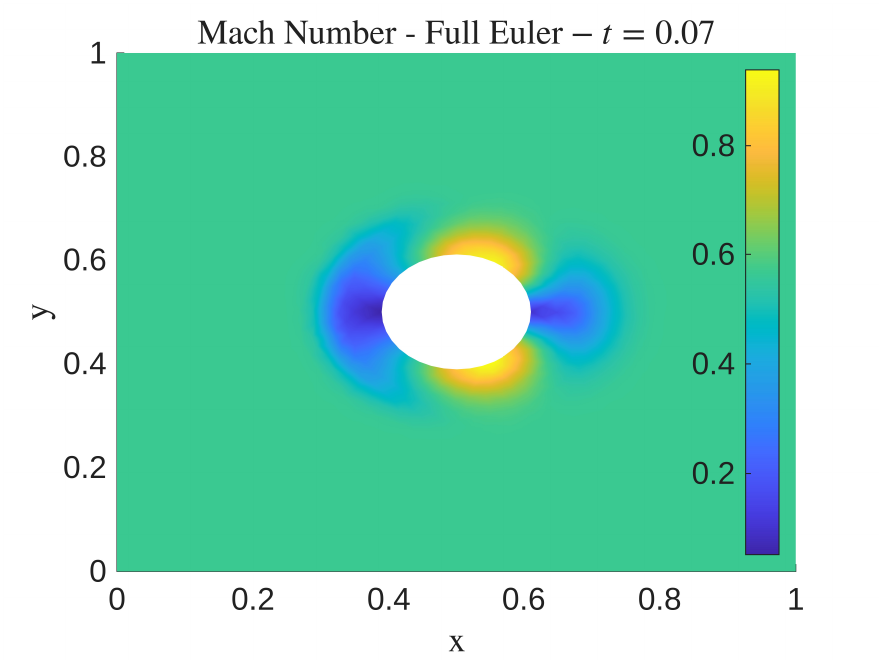}
        \includegraphics[width=0.32\linewidth]{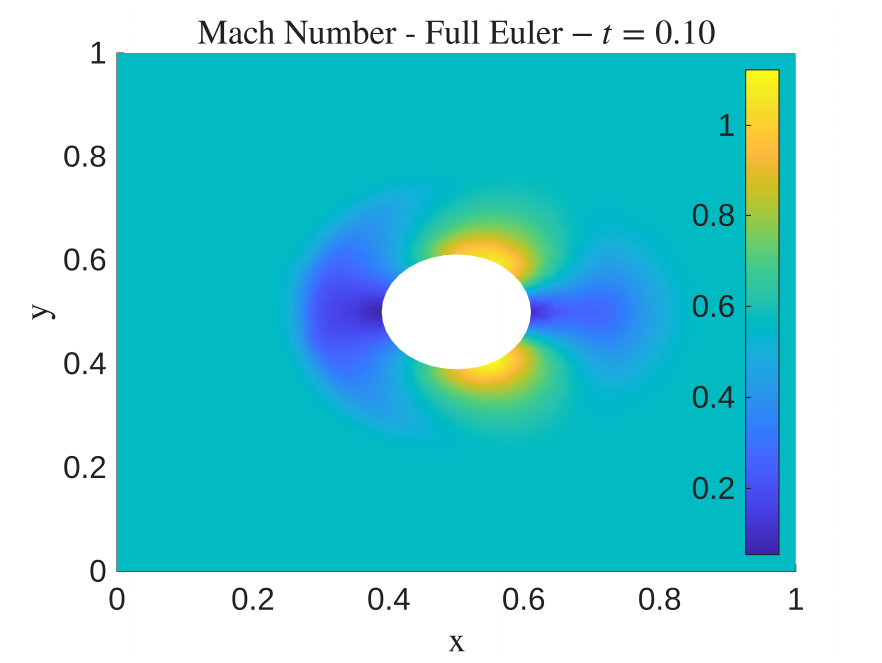}
    
        \smallskip
    
        \includegraphics[width=0.32\linewidth]{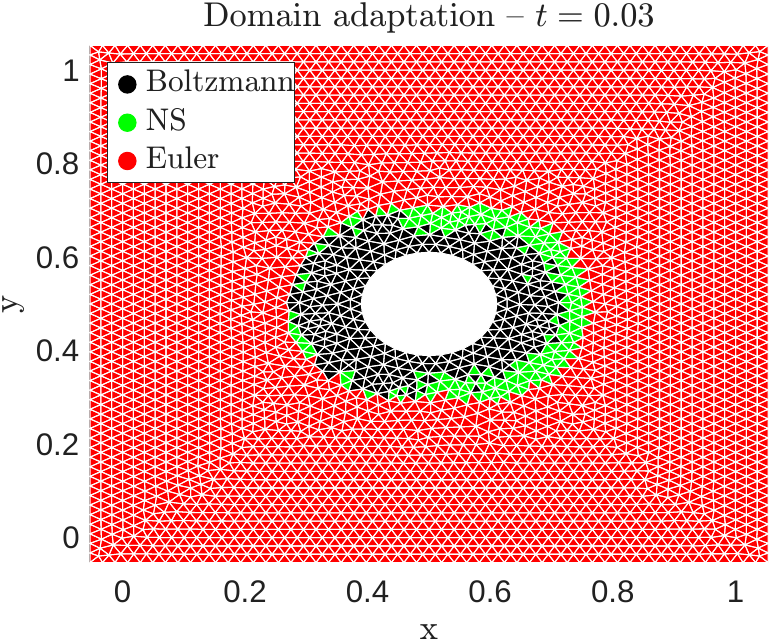}
        \includegraphics[width=0.32\linewidth]{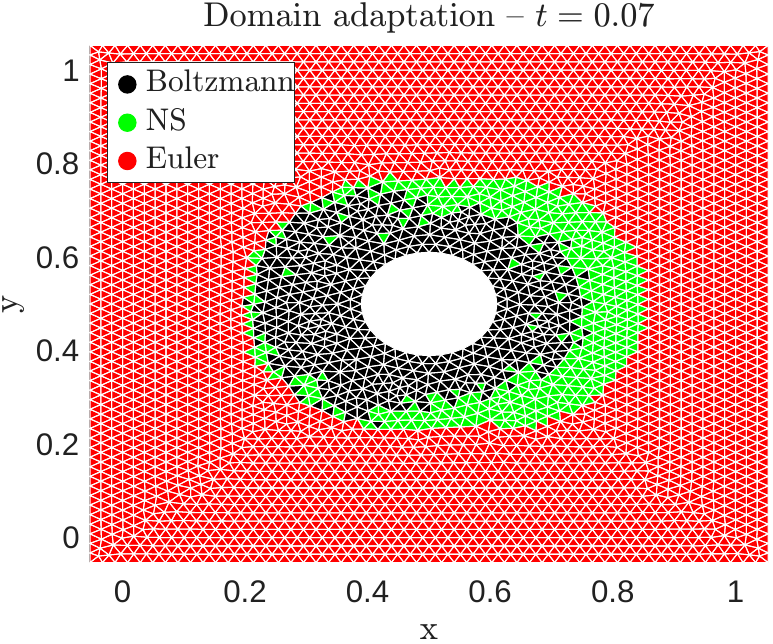}
        \includegraphics[width=0.32\linewidth]{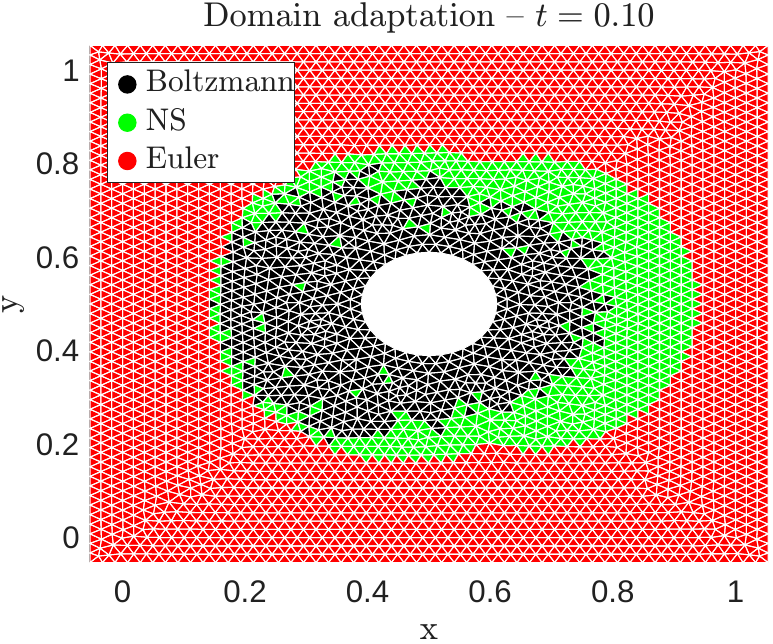}
        
        \caption{\hyperref[sec:test5]{Test 5}: Time evolution of the Mach number $\text{Ma}$ at $t=0.03$ (left), $t=0.07$ (center), and $t=0.10$ (right). Rows 1 to 4 show the solutions for the full Boltzmann, Hybrid, full Navier--Stokes, and full Euler solvers, respectively. The bottom row illustrates the dynamic domain partition (red: Euler, green: Navier--Stokes, black: Boltzmann).}
        \label{Fig6}
    \end{figure}

    Figure~\ref{Fig6} displays the time evolution of the Mach number profiles obtained with the different solvers, alongside the corresponding dynamic domain partitions generated by the adaptive framework. For a better qualitative comparison, the same colorbar is employed across all solvers, with the exception of the Euler model to better highlight the localized discrepancies. In this context, the Mach number is defined as $\text{Ma}=u/\sqrt{\gamma T}$, where $u=\sqrt{u_x^2+u_y^2}$ and $\gamma=1+2/d_v$ is the specific heat ratio. 

    \subsubsection{Test 6: Flow over a NACA 0012 airfoil}\label{sec:test6}
     The third test case considers a rarefied flow past a NACA 0012 airfoil at a high angle of attack, $\alpha = 30^\circ$, inspired by \cite{bos-dim-2020, bos-dim-2021}. As for the disk, strong non-equilibrium effects are expected near the solid boundary and in the wake region. The adaptive strategy should therefore assign the Euler model to the free stream, where the gas remains close to equilibrium, the Navier--Stokes equations to the regions with a moderate departure from equilibrium, and the full Boltzmann solver to the vicinity of the airfoil and to the wake, where the continuum description is no longer reliable.

\begin{figure}[t!]
        \centering
        \includegraphics[width=0.48\linewidth]{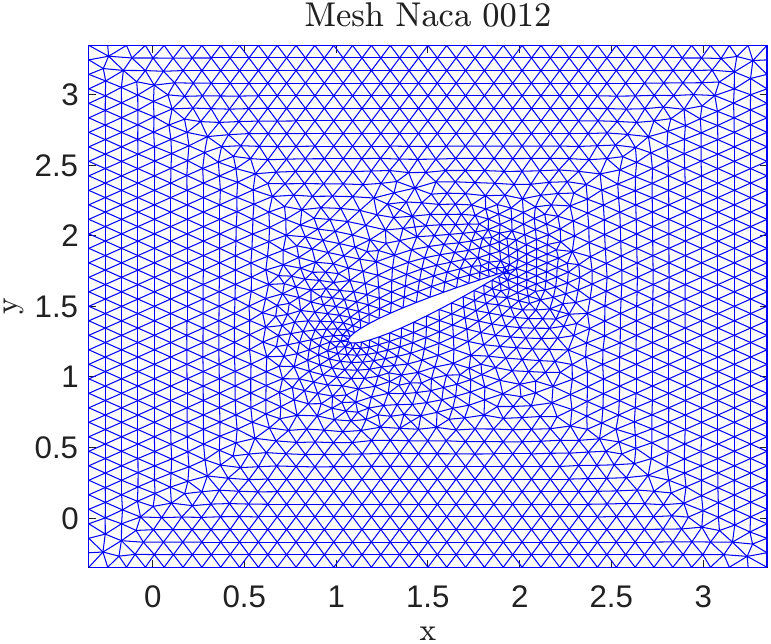}
        \caption{\hyperref[sec:test6]{Test 6}: Computational mesh layout used for the spatial discretization, including the ghost cells at the boundaries.}
        \label{Fig9}
    \end{figure}

    The computational domain is $\mathcal{T} = [0, 3]\times[0, 3]$, discretized with approximately $2400$ cells. The corresponding mesh layout, which includes the ghost cells used for the spatial discretization, is depicted in Figure~\ref{Fig9}. Zero-gradient boundary conditions are applied at the external boundaries, while slip conditions are enforced on the NACA 0012 airfoil. The velocity space $\mathcal{V}=[-6,6]^2$ is discretized using a regular $32 \times 32$ Cartesian grid. The macroscopic fields are initialized with the uniform state $W=(\rho, u_x, u_y, T) = (1,\,0.8,\,0,\,1)$, and the initial kinetic distribution function is set to the corresponding local Maxwellian state. 

    As in \hyperref[sec:test5]{Test 5}, all cells are initially assigned to the kinetic description.

    The numerical simulation is executed in the rarefied regime with $\varepsilon=0.1$. 
    To evaluate the accuracy of the coupling strategy, we track the space-time averaged $L^1$ error for a representative macroscopic variable, namely the temperature $T$. This metric is computed by integrating the absolute temperature difference relative to the full Boltzmann reference solution over the fluid cells, normalized by their total area and by the number of stored time instants. Mathematically, the error for a given solver is defined as:
    \begin{equation}
        \text{Error}_{b} = \frac{1}{N_t \sum_{i\in\mathcal{F}} |P_i|}
        \sum_{n=1}^{N_t} \sum_{i\in\mathcal{F}} \left| T_{b,i}^n - T_{\text{ref}, i}^n \right| |P_i|,
    \end{equation}
    where $b \in \{\text{Euler}, \text{NS}, \text{Hybrid}\}$ denotes the specific solver under consideration, $\mathcal{F}$ is the set of fluid cells, $N_t$ is the number of stored time instants, and $|P_i|$ is the area of the cell $P_i$. The term $T_{b, i}^n$ represents the temperature computed by solver $b$ at cell $i$ and time level $n$, while $T_{\text{ref}, i}^n$ is the corresponding reference temperature obtained from the full Boltzmann simulation.
   
    The resulting average temperature errors are $\text{Error}_{\text{Hybrid}}=1.7 \times 10^{-3}$, $\text{Error}_{\text{NS}}=3.2 \times 10^{-3}$, and $\text{Error}_{\text{Euler}}=6.3 \times 10^{-3}$, showing that the hybrid solver is almost twice as accurate as the full Navier--Stokes model. 

    In terms of computational efficiency, the execution times on the GPU are: $\text{Duration}_{\text{Euler}} = 9.297\,\text{min}$, $\text{Duration}_{\text{NS}} = 9.202\,\text{min}$, $\text{Duration}_{\text{Hybrid}} = 27.33\,\text{min}$, and $\text{Duration}_{\text{Boltzmann}} = 65.71\,\text{min}$. The adaptive framework provides a computational speedup of approximately $2.4$ compared to the full Boltzmann solver.\\
    The CPU and GPU execution times and a detailed analysis of the hybrid speedup are given in Section~\ref{secNumEff}.
    In the present numerical simulation, the threshold values for the switching criteria (Equations \eqref{crit:CompEuler_Main} and \eqref{criterioSecondoSwitch}) are set to $\eta_1 = \ve/10$ and $\eta_2 = \ve/2$.

\begin{figure}
        \centering
        \includegraphics[width=0.32\linewidth]{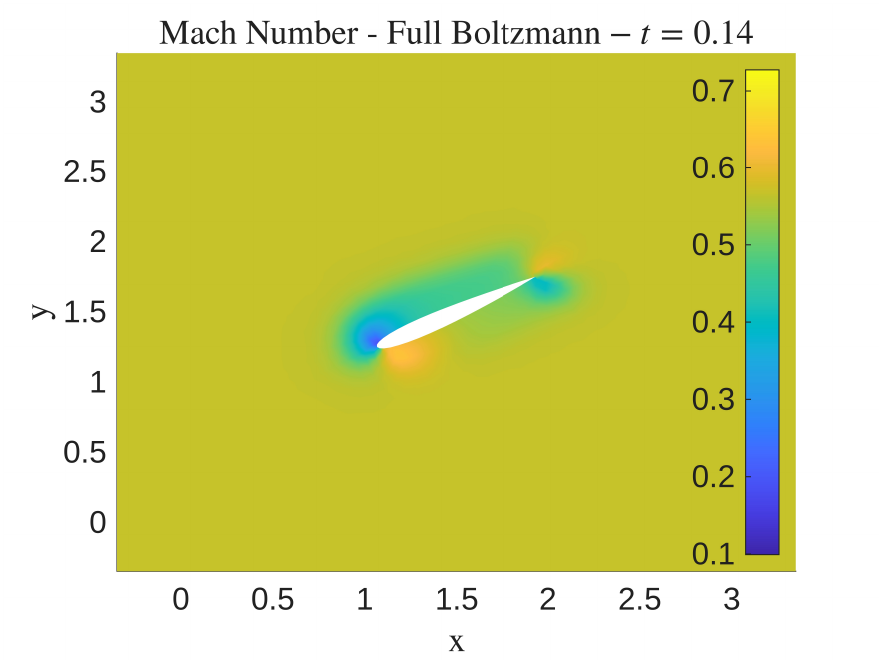}
        \includegraphics[width=0.32\linewidth]{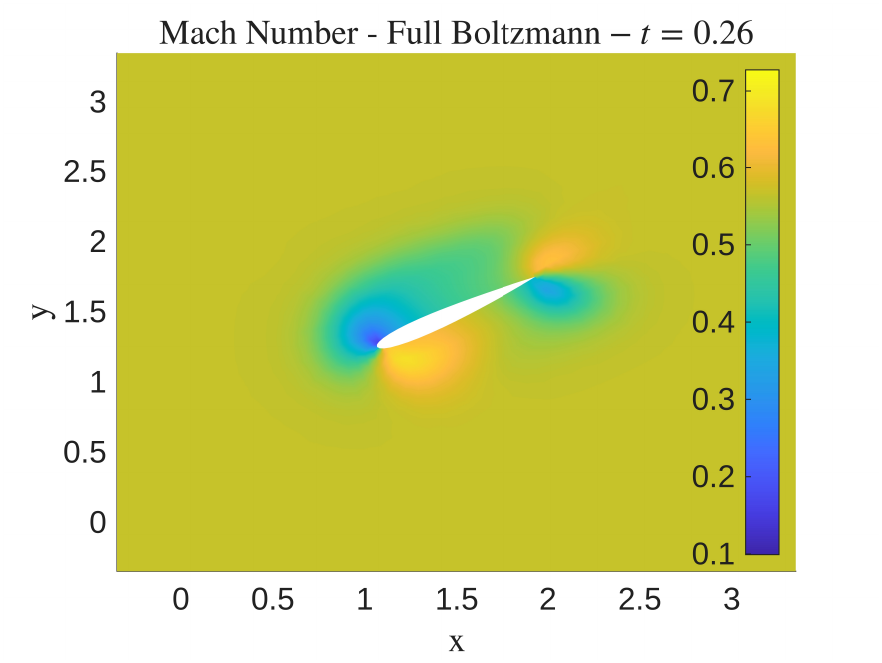}
        \includegraphics[width=0.32\linewidth]{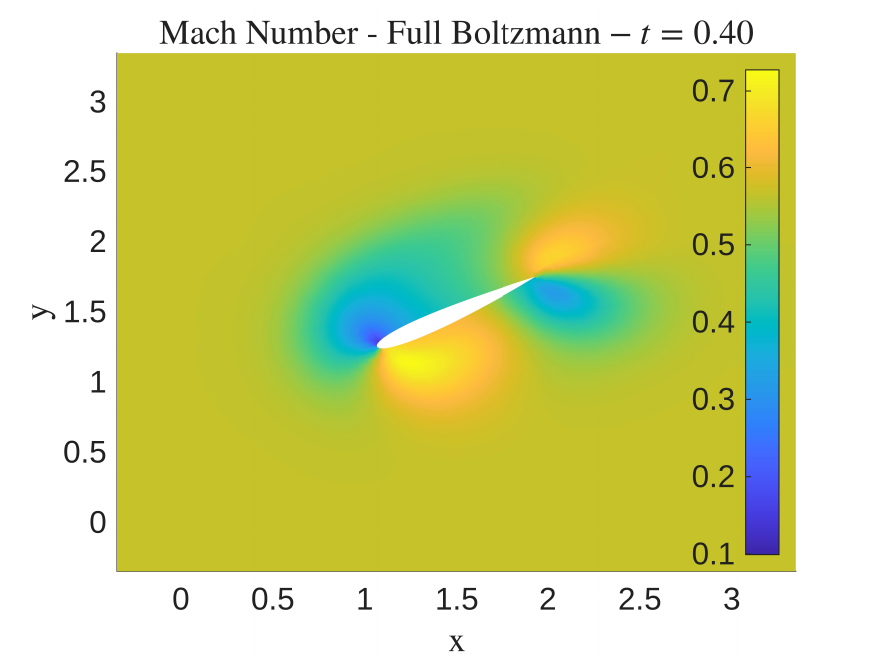}
        
        \smallskip
        
        \includegraphics[width=0.32\linewidth]{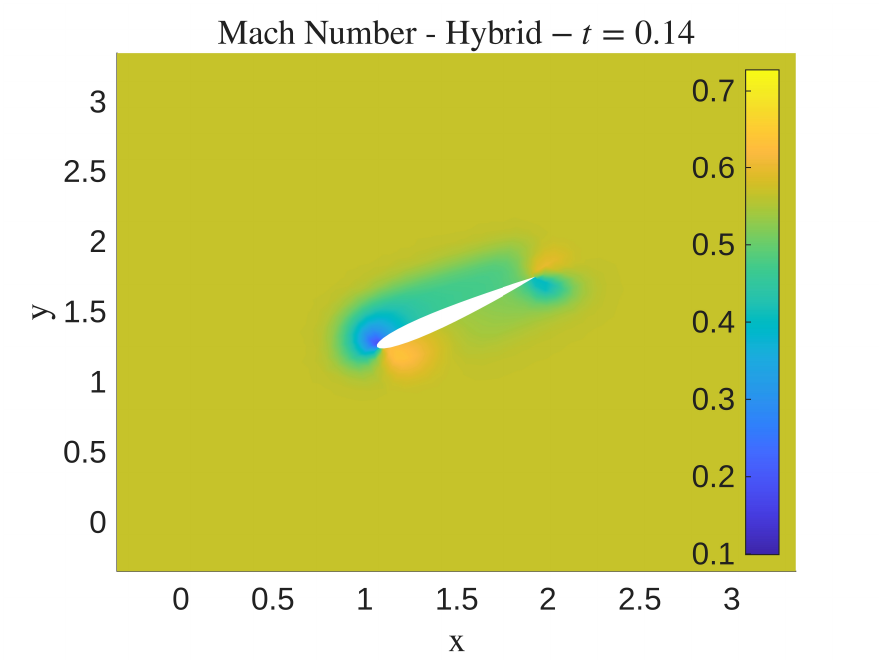}
        \includegraphics[width=0.32\linewidth]{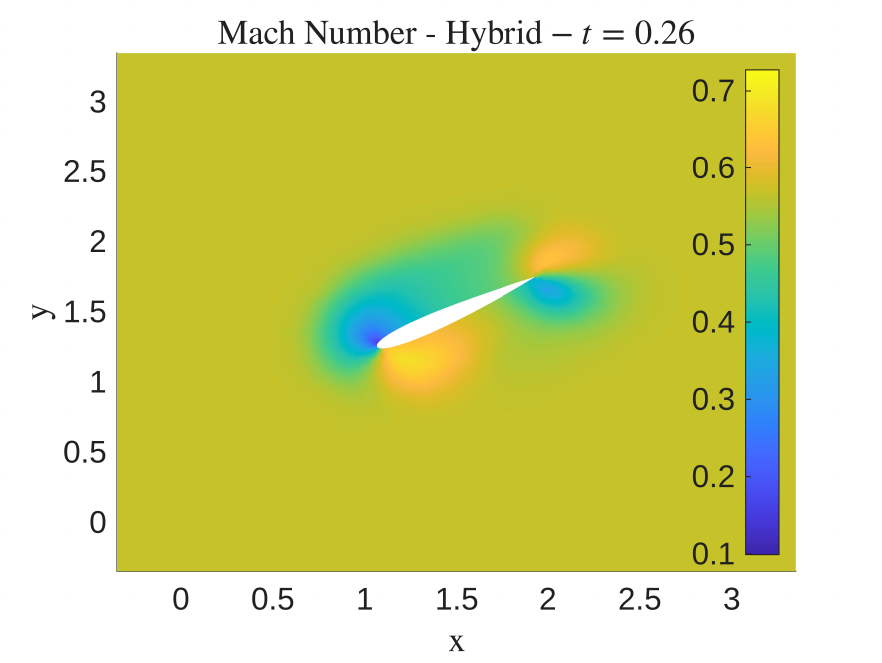}
        \includegraphics[width=0.32\linewidth]{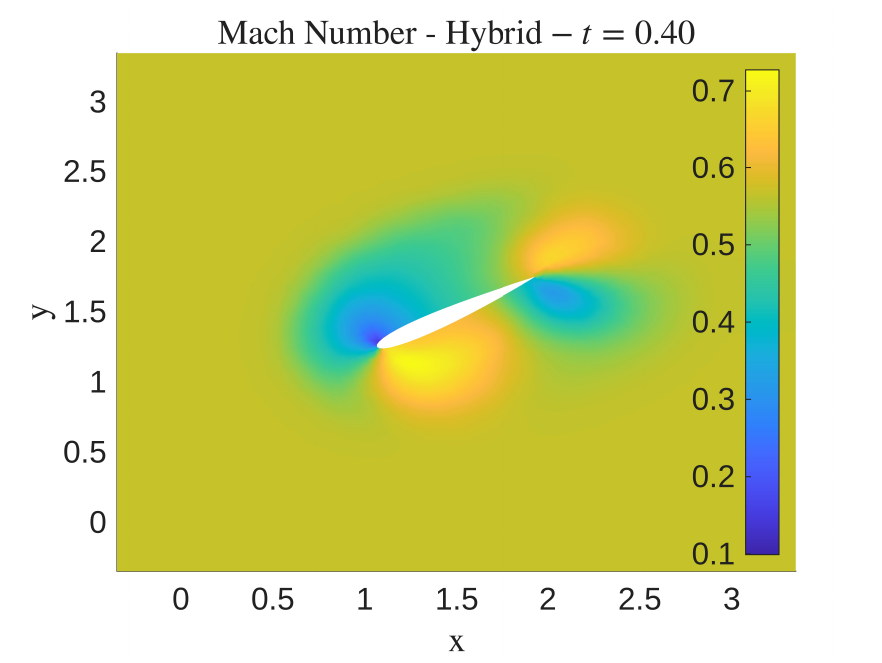}
        
        \smallskip
    
        \includegraphics[width=0.32\linewidth]{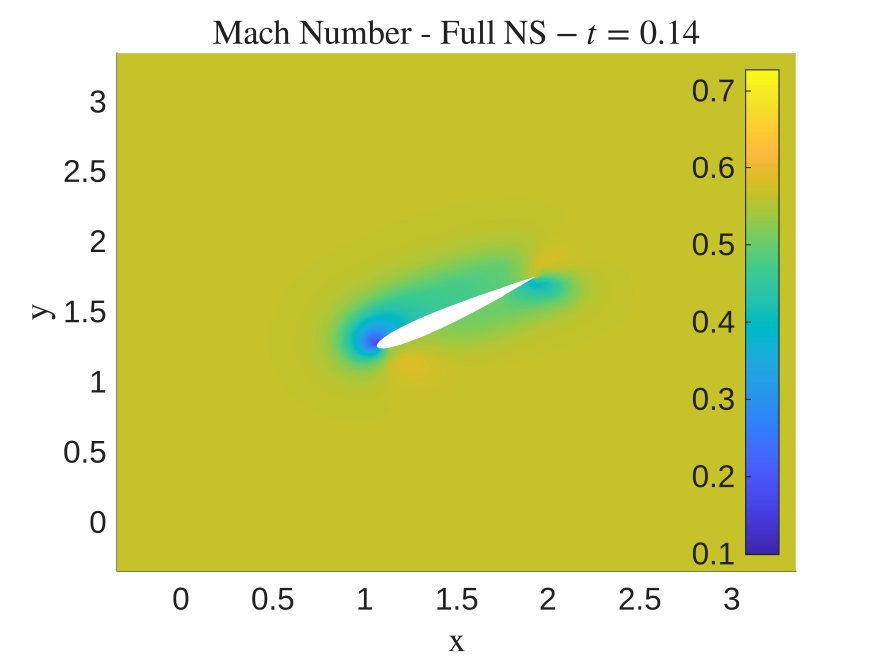}
        \includegraphics[width=0.32\linewidth]{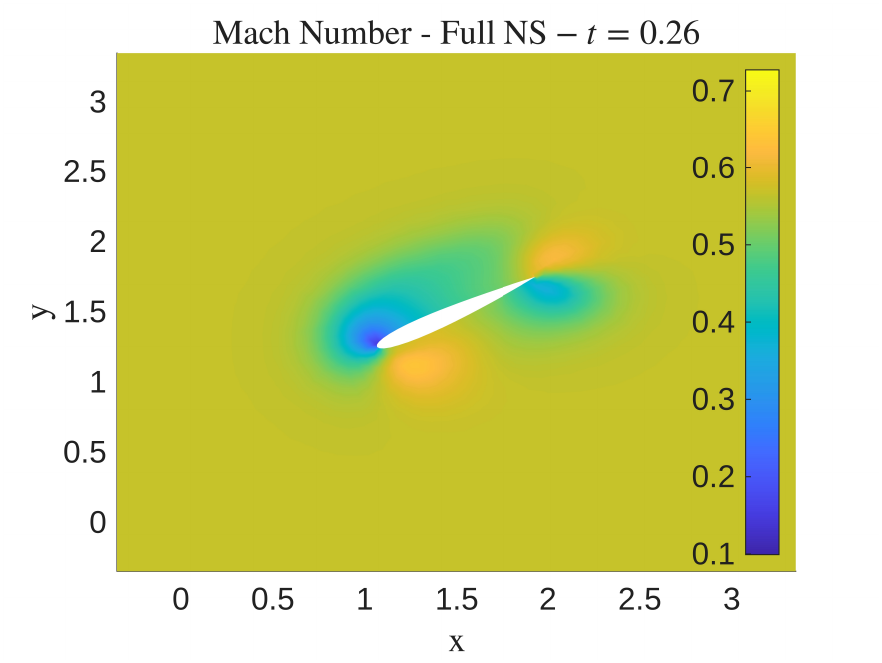}
        \includegraphics[width=0.32\linewidth]{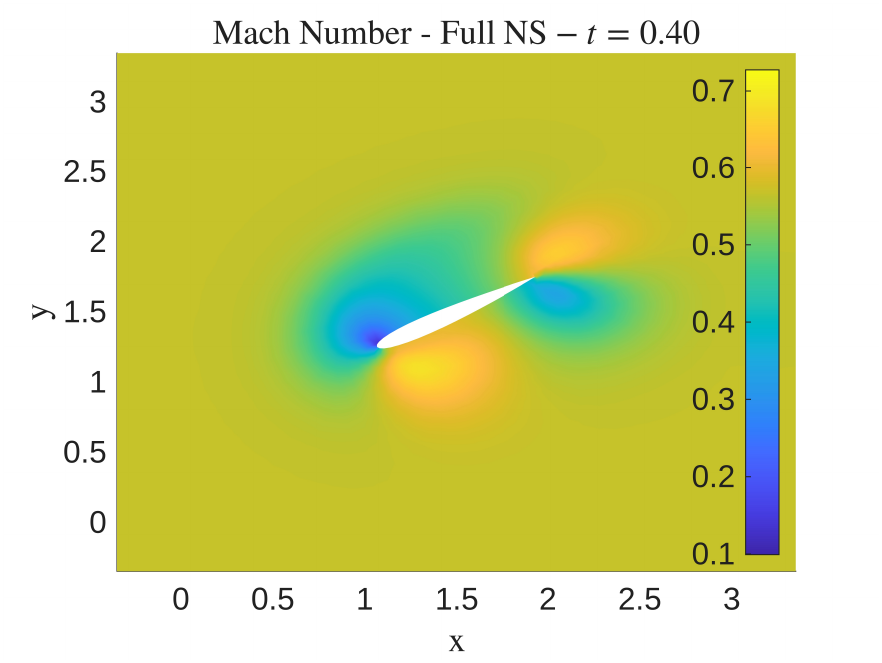}
    
        \smallskip
    
        \includegraphics[width=0.32\linewidth]{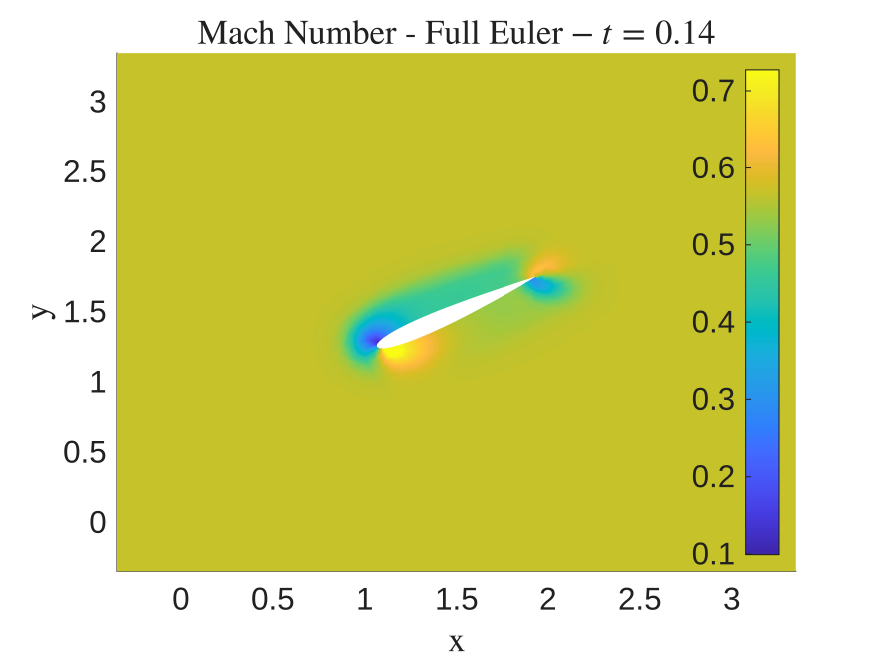}
        \includegraphics[width=0.32\linewidth]{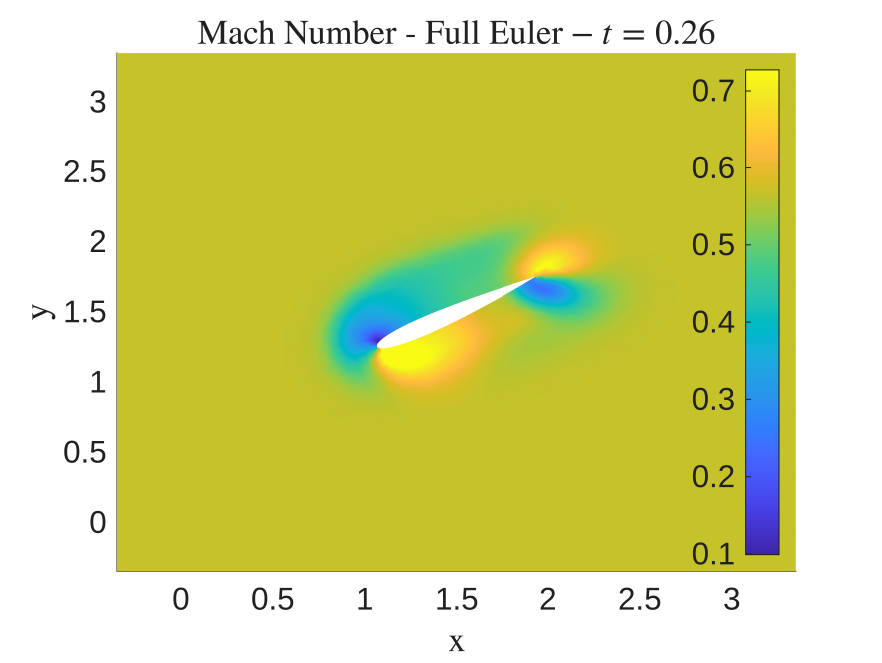}
        \includegraphics[width=0.32\linewidth]{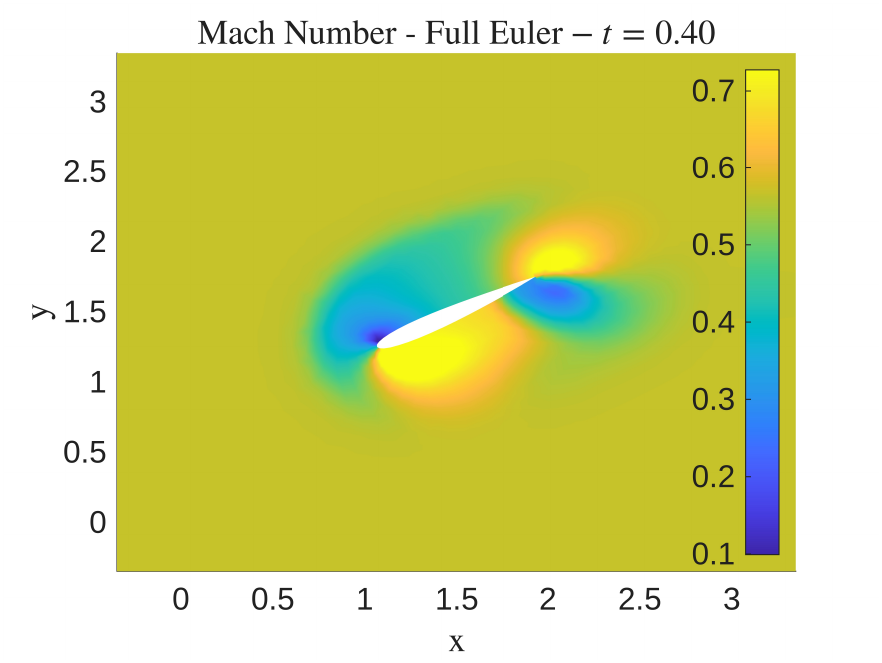}
    
        \smallskip
    
        \includegraphics[width=0.32\linewidth]{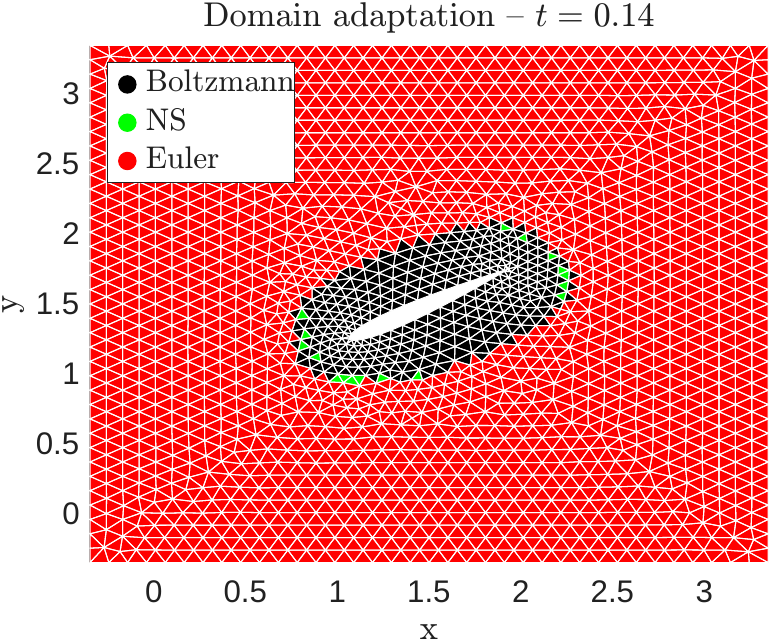}
        \includegraphics[width=0.32\linewidth]{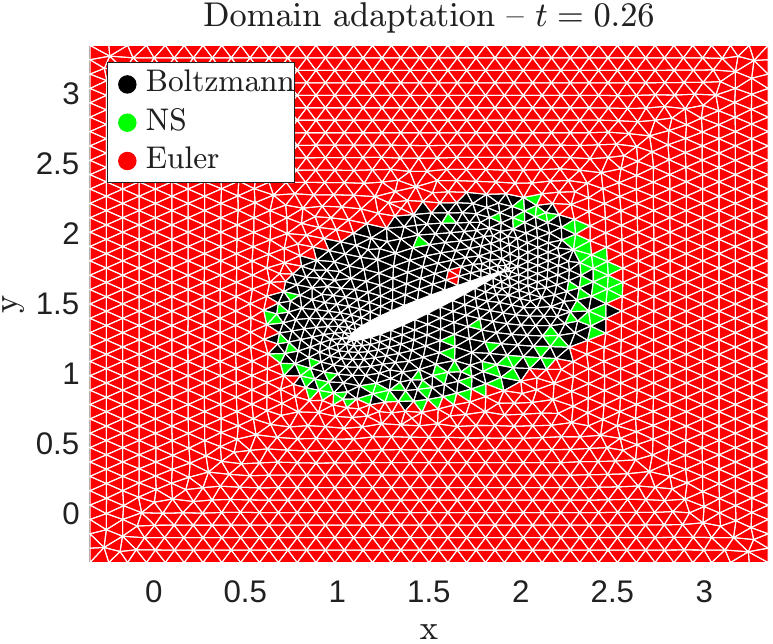}
        \includegraphics[width=0.32\linewidth]{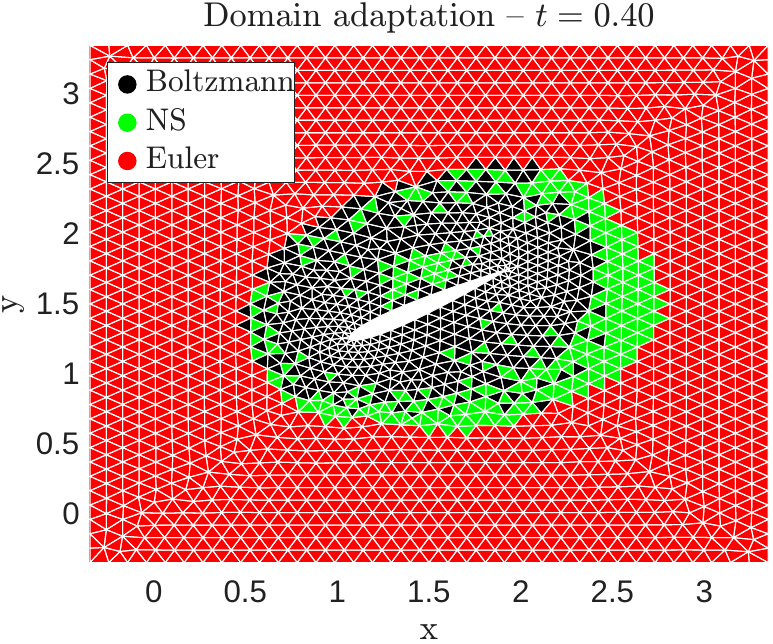}
        
        \caption{\hyperref[sec:test6]{Test 6}: Time evolution of the Mach number $\text{Ma}$ at $t=0.14$ (left), $t=0.26$ (center), and $t=0.40$ (right). Rows 1 to 4 show the solutions for the full Boltzmann, Hybrid, full Navier--Stokes, and full Euler solvers, respectively. The bottom row illustrates the dynamic domain partition (red: Euler, green: Navier--Stokes, black: Boltzmann).}
        \label{Fig8}
    \end{figure}

    Figure~\ref{Fig8} displays the time evolution of the Mach number profiles obtained with the different solvers, alongside the corresponding dynamic domain partitions generated by the adaptive framework. For a better qualitative comparison, the same colorbar, based on the range of the full Boltzmann solution, is employed across all solvers. In this context the Mach number is defined as $\text{Ma}=u/\sqrt{\gamma T}$, where $u=\sqrt{u_x^2+u_y^2}$ and $\gamma=1+2/d_v$ is the specific heat ratio. 
    As shown in Figure~\ref{Fig9}, the mesh is refined near the airfoil, where the elements are significantly smaller.

    \subsection{Numerical efficiency}\label{secNumEff}
    In this section, we discuss the numerical efficiency of the proposed domain decomposition strategy.
    Tables \ref{table_gpu} and \ref{table_cpu} summarize the execution times and the space-time averaged $L^1$ temperature errors for the Sod shock tube (\hyperref[sec:test4]{Test 4}), the flow over a disk (\hyperref[sec:test5]{Test 5}), and the flow over a NACA 0012 airfoil (\hyperref[sec:test6]{Test 6}) on GPU and CPU architectures, respectively. As previously defined in Equation~\eqref{eqL1_1}, this error metric evaluates the space-time averaged $L^1$ difference of the temperature $T$ relative to the full Boltzmann reference solution, computed over the fluid cells and normalized by their total area and by the number of stored time instants. This metric is omitted for the Sod test: since this test is performed in a highly stiff regime, all the models naturally converge toward the correct hydrodynamic limit, and the Boltzmann solution cannot be treated as the absolute reference in the same manner as for the other tests, which operate in the rarefied regime.

    As observed in the data, the three-level hybrid domain decomposition strategy significantly outperforms the full Navier--Stokes and Euler solvers in terms of accuracy. Furthermore, compared to the full Boltzmann solver, it yields a computational speedup between $2.4$ and $2.8$ on the GPU, and of up to $7.8$ on the CPU.

    This difference between the two architectures becomes clear when the execution times of each solver on the CPU and on the GPU are compared (Tables~\ref{table_gpu} and~\ref{table_cpu}). In all three test cases, the full Euler and Navier--Stokes simulations are faster on the CPU than on the GPU, by a factor of about $1.3$--$1.8$, whereas the full Boltzmann solver is about $3.0$--$3.3$ times faster on the GPU. The macroscopic models evolve only four moments per cell: their workload is too light to saturate the GPU in the present implementation. The kinetic solver, instead, evolves a large number of discrete velocities in each cell and fully exploits the massive parallelism of the GPU. Consequently, the relative gain obtained by replacing the Boltzmann solver with the macroscopic models in most of the domain is larger on the CPU, where the macroscopic models are comparatively cheaper.

    For the hybrid simulations, the benefit of the GPU therefore depends on the weight of the kinetic part. The hybrid simulation is about $1.9$ times faster on the GPU for the NACA 0012 airfoil ($27.3$~min versus $51.4$~min), where the kinetic workload per time step is the largest, but only about $1.2$ and $1.1$ times faster for the disk and the Sod shock tube, respectively; in the latter case, the kinetic region is confined to a few cells. This indicates that, for strongly localized kinetic regions, the macroscopic part of the computation does not benefit from the GPU.

    Finally, the hybrid simulations are always faster on the GPU, and, when this hardware acceleration is combined with the algorithmic efficiency of the hybrid scheme, they are about $8$ times faster ($7.9$--$8.6$) than the corresponding full Boltzmann simulations on the CPU.

    To provide a deeper insight into how computational resources are allocated within the three-level hybrid solver, Tables \ref{table_breakdown_gpu} and \ref{table_breakdown_cpu} present a detailed cost breakdown of the distinct algorithmic phases (Euler, Navier--Stokes, Boltzmann collision and transport, indicators calculation, and solver coupling) for GPU and CPU computations, respectively. A key feature of these breakdowns is that, within the kinetic region, the computational cost is dominated by the transport step rather than by the collision operator. Kinetic transport accounts for $21$--$36\,\%$ of the total cost on the GPU and $48$--$74\,\%$ on the CPU, whereas the collision operator accounts for only $6$--$10\,\%$ and $8$--$9\,\%$, respectively. In other words, the cost of the kinetic solver is largely dominated by the transport step, especially on the CPU. This is characteristic of kinetic solvers on unstructured meshes, where the transport step requires the evaluation of numerical fluxes across all cell faces for every discrete velocity.
    This observation has a direct consequence on the design of the model hierarchy. A natural candidate for the intermediate level would be a simplified kinetic model such as the ES-BGK equation, which replaces the Boltzmann collision operator with a cheaper relaxation operator. However, the ES-BGK and Boltzmann models share exactly the same transport operator and the same phase-space discretization, so that the transport cost, which dominates the kinetic workload, would remain unchanged. Replacing the Boltzmann collision operator with an ES-BGK relaxation would therefore only reduce the minor part of the cost, yielding a limited gain. For this reason, we adopted the Navier--Stokes equations as the intermediate level of the hierarchy: they remove the velocity-space discretization altogether, thus eliminating both the transport and the collision costs, while still accounting for the first non-equilibrium corrections (viscous stresses and heat flux) beyond the Euler model. Moreover, as shown in Tables \ref{table_gpu} and \ref{table_cpu}, the full Navier--Stokes simulations run at essentially the same cost as the full Euler simulations. The intermediate level is thus more accurate than Euler at a comparable cost, which explains the favourable cost-to-accuracy ratio of the proposed adaptive approach.
    
    \begin{table}[htpb]
    \centering
    \begin{tabular}{lcccc}
        \toprule
        \textbf{Simulation} & \textbf{Solver} & \textbf{Error} & \makecell[c]{\textbf{Execution} \\ \textbf{Time} \\ \textbf{GPU}} & \makecell[c]{\textbf{Hybrid} \\ \textbf{Speedup} \\ \textbf{GPU}}\\
        \midrule
        \multirow{4}{*}{Sod}      & Boltzmann              &                           &     1053 s  & 2.8 \\
                                  & Hybrid                 &                           &     370.3 s & 1.0 \\
                                  & NS                     &                           &     186.3 s & 0.50 \\
                                  & Euler                  &                           &     191.9 s & 0.52 \\
        \midrule
        \multirow{4}{*}{Disk}     & Boltzmann              &                           &     621.0 s    & 2.5 \\
                                  & Hybrid                 &       $4.6\times 10^{-3}$ &     247.3 s    & 1.0 \\
                                  & NS                     &       $9.3\times 10^{-3}$ &     121.0 s    & 0.49 \\
                                  & Euler                  &       $13 \times 10^{-3}$ &     122.7 s    & 0.50 \\
        \midrule
        \multirow{4}{*}{NACA 0012} & Boltzmann              &                           &     65.71 min    & 2.4 \\
                                  & Hybrid                 &       $1.7\times 10^{-3}$ &     27.33 min    & 1.0 \\
                                  & NS                     &       $3.2\times 10^{-3}$ &     9.202 min   & 0.34 \\
                                  & Euler                  &       $6.3\times 10^{-3}$ &     9.297 min   & 0.34 \\
        \bottomrule
    \end{tabular}
    \caption{GPU Performance: Comparison of execution times, hybrid scheme speedups, and space-time averaged $L^1$ temperature errors for the Sod shock tube (\hyperref[sec:test4]{Test 4}), flow over a disk (\hyperref[sec:test5]{Test 5}), and flow over a NACA 0012 airfoil (\hyperref[sec:test6]{Test 6}).}
    \label{table_gpu}
    \end{table}

    \begin{table}[htpb]
    \centering
    \begin{tabular}{lcccc}
        \toprule
        \textbf{Simulation} & \textbf{Solver} & \textbf{Error} & \makecell[c]{\textbf{Execution} \\ \textbf{Time} \\ \textbf{CPU}} & \makecell[c]{\textbf{Hybrid} \\ \textbf{Speedup} \\ \textbf{CPU}}\\
        \midrule
        \multirow{4}{*}{Sod}      & Boltzmann              &                           &     3171 s & 7.8 \\
                                  & Hybrid                 &                           &     407.7 s & 1.0 \\
                                  & NS                     &                           &     141.6 s & 0.35 \\
                                  & Euler                  &                           &     151.1 s & 0.37 \\
        \midrule
        \multirow{4}{*}{Disk}     & Boltzmann              &                           &     2014 s   & 6.8 \\
                                  & Hybrid                 &       $4.6\times 10^{-3}$ &     295.1 s  & 1.0 \\
                                  & NS                     &       $9.3\times 10^{-3}$ &     83.98 s  & 0.28 \\
                                  & Euler                  &       $13 \times 10^{-3}$ &     85.94 s  & 0.29 \\
        \midrule
        \multirow{4}{*}{NACA 0012} & Boltzmann              &                           &     216.4 min    & 4.2 \\
                                  & Hybrid                 &       $1.7\times 10^{-3}$ &     51.37 min   & 1.0 \\
                                  & NS                     &       $3.2\times 10^{-3}$ &     5.236 min   & 0.10 \\
                                  & Euler                  &       $6.3\times 10^{-3}$ &     5.828 min   & 0.11 \\
        \bottomrule
    \end{tabular}
    \caption{CPU Performance: Comparison of execution times, hybrid scheme speedups, and space-time averaged $L^1$ temperature errors for the Sod shock tube (\hyperref[sec:test4]{Test 4}), flow over a disk (\hyperref[sec:test5]{Test 5}), and flow over a NACA 0012 airfoil (\hyperref[sec:test6]{Test 6}).}
    \label{table_cpu}
    \end{table}

    \begin{table}[htpb]
    \centering
    \begin{tabular}{lccc}
        \toprule
        \textbf{Computational Phase GPU} & \textbf{Sod} & \textbf{Disk} & \textbf{NACA 0012} \\
        \midrule
        Euler Solver                 & 23 \%    & 25 \%        & 18 \% \\
        Navier--Stokes Solver        & 27 \%    & 26 \%        & 19 \% \\
        Kinetic Transport            & 23 \%    & 21 \%        & 36 \% \\
        Boltzmann Collision          & 7  \%     & 6 \%         & 10 \% \\
        Indicators Calculation       & 3 \%     & 4 \%         & 2 \% \\
        Solvers Coupling             & 10 \%     & 9 \%         & 7 \% \\
        Other operations             & 7 \%     & 9 \%         & 8 \% \\
        \bottomrule
    \end{tabular}
    \caption{GPU Breakdown: Relative computational cost of the different phases within the three-level hybrid solver.}
    \label{table_breakdown_gpu}
    \end{table}

    \begin{table}[htpb]
    \centering
    \begin{tabular}{lccc}
        \toprule
        \textbf{Computational Phase CPU} & \textbf{Sod} & \textbf{Disk} & \textbf{NACA 0012} \\
        \midrule
        Euler Solver                 & 13 \%    & 13 \%    & 5 \% \\
        Navier--Stokes Solver        & 14 \%    & 10 \%     & 4 \% \\
        Kinetic Transport            & 48 \%   & 57 \%    & 74 \% \\
        Boltzmann Collision          & 9 \%   & 8 \%     & 9 \% \\
        Indicators Calculation       & 3 \%    & 2 \%     & 1 \% \\
        Solvers Coupling             & 9 \%   & 5 \%     & 3 \% \\
        Other operations             & 4 \%   & 5 \%     & 4 \% \\
        \bottomrule
    \end{tabular}
    \caption{CPU Breakdown: Relative computational cost of the different phases within the three-level hybrid solver.}
    \label{table_breakdown_cpu}
    \end{table}

    \section{Conclusion}\label{secCon}
    In this work, we have proposed a hierarchical domain decomposition strategy for solving the full Boltzmann equation and, more generally, kinetic models in multiscale regimes. The switching criterion, inspired by \cite{lev-mor-nad-1998,fil-rey-2014,fil-xio-2018,tiw-2000,li-son-wan-2021,xio-qiu-2017,cap-par-rey-2025,cap-ten-2025}, uses truncated Chapman--Enskog expansions to determine when to switch locally between successive levels of description. This criterion relies solely on macroscopic quantities derived from the closure of the kinetic 
system, avoiding the need to evaluate the distribution function in the macroscopic regions.

A key novelty of the present work lies in the extension of this adaptive strategy to a 2D-in-space and 2D-in-velocity framework operating on complex unstructured meshes. More importantly, the Euler, Navier--Stokes, and Boltzmann models are not advanced as independent solvers coupled only through interface conditions. Instead, they are integrated within a unified AP-IMEX Runge--Kutta framework, share the same temporal stages, and exchange kinetic and macroscopic representations at every stage. This stage-wise coupling is used to maintain temporal consistency across model interfaces and compatibility with the fluid limit.

The framework also introduces a refined three-level model hierarchy. In the previous work \cite{cap-par-rey-2025}, the intermediate fidelity level between the Euler equations and the full Boltzmann model was represented by the ES-BGK relaxation model. However, on unstructured spatial meshes, the overall computational cost of kinetic models is largely dominated by the transport 
phase rather than the collision step. Consequently, evolving a microscopic distribution function with the ES-BGK model offers only limited performance benefits compared with the full Boltzmann equation. To address this computational bottleneck, the present approach replaces the intermediate ES-BGK model with its macroscopic Navier--Stokes counterpart. The resulting hierarchy employs the Euler equations in regions close to local thermodynamic equilibrium, the Navier--Stokes equations in regions of mild non-equilibrium, and the full Boltzmann equation in the strongly non-equilibrium kinetic regime. These three models are dynamically coupled throughout the computational domain.

In the kinetic regions, the framework employs a spectral discretization of the full Boltzmann collision operator rather than a simplified BGK-type relaxation model. The method combines an asymptotic-preserving IMEX time discretization, high-order finite-volume methods on triangular meshes, and a fast spectral solver for the Boltzmann collision integral. The numerical tests confirm third-order temporal accuracy in the rarefied and fluid limits and third-order spatial accuracy for the full kinetic discretization. As expected for the selected IMEX scheme, order reduction is observed in the intermediate Knudsen regime, and the method is therefore not uniformly accurate with respect to the Knudsen number.

The numerical experiments show that the computationally expensive kinetic regions are dynamically confined in space and time. In the configurations considered here, the hybrid solutions remain close to the corresponding full Boltzmann references while improving the accuracy obtained with the standalone Euler and Navier--Stokes solvers. Relative to the full Boltzmann 
computation, the observed speedup ranges from 2.4 to 2.8 on GPU and from 4.2 to 7.8 on CPU. These gains depend on the problem, the switching thresholds, the spatial and velocity discretizations, and the hardware architecture.

Overall, the results indicate that deterministic full Boltzmann solvers can be incorporated efficiently into multiscale simulations on nontrivial geometries by means of a tightly coupled hierarchy of kinetic and fluid models. This opens new perspectives for applications of kinetic theory to physical and engineering problems, including aerodynamics in transitional regimes. Future work will address more accurate non-equilibrium lifting at kinetic--fluid interfaces, uniformly accurate time discretizations, multi-species mixtures, and collisional plasma models, together with high-performance computing strategies tailored to the dynamic model hierarchy.

\section*{Acknowledgments}
    The authors received funding from the European Union's Horizon Europe research and innovation program under the Marie Skłodowska-Curie Doctoral Network DataHyking (Grant No. 101072546). 
    DC and TR are also supported by the French government, through the UniCA$_{JEDI}$ Investments in the Future project managed by the National Research Agency (ANR) with the reference number ANR-15-IDEX-01.
    DC acknowledges the support of the Department of Mathematics and Computer Science of the University of Ferrara through access to their computational resources. The research of LP has been supported by the Royal Society under the Wolfson Fellowship “Uncertainty quantification, data-driven simulations and learning of multiscale complex systems governed by PDEs”. LP also acknowledges the  partial support by Fondo Italiano per la Scienza (FIS2023-01334) advanced grant ADAMUS. This work has been written within the activities of GNCS group of INdAM (Italian National Institute of High Mathematics).

\appendix

\section{Chapman--Enskog expansion and hydrodynamic limits}
\label{app:ChapmanEnskog}

In this appendix, we consider a three-dimensional velocity space. Under the assumption that the distribution function $f$ is close to the local equilibrium due to the entropy relaxation properties, we can formally introduce the Chapman--Enskog expansion:
\begin{equation}
  \label{devChapEnsk_App}
     f = \, \mathcal M_{\rho,\bm{u},T}\left[1 + \ve \, g^{(1)} + \ve^2 \, g^{(2)} + \ldots \right ],
\end{equation}
where the fluctuations $g^{(i)}$ for $i \geq 1$ smoothly depend on the macroscopic moments of the distribution:
\begin{equation*}
    (\rho,\rho\bm{u},E) \, = \, \int_{\RR^{3}} f(v) \, \left(1, v, \frac{|v|^2}{2}\right) \, dv, \quad T \, = \, \frac{1}{3 \, \rho }\int_{\RR^{3}} f(v) \, |v-\bm{u}|^2 \, dv,
\end{equation*}
and satisfy the compatibility conditions:
\begin{equation*}
    \int_{\RR^{3}} \mathcal M_{\rho,\bm{u},T} \, \, g^{(i)}(v) \, \left(1, v, \frac{|v|^2}{2}\right) \, dv \, = \, \bm{0}_{\RR^{5}}^\intercal.
\end{equation*}
From the conservative properties of the collision operator, the macroscopic moments satisfy the unclosed system:
\begin{equation}
    \label{eqHydroNonClosed_App}
    \left\{ \begin{aligned}
    & \partial_t \rho  + \diverg_x (\rho \, \bm{u} ) \,= \, 0, \\
    & \partial_t(\rho \, \bm{u} ) + \diverg_x \left (\rho \bm{u} \otimes \bm{u} + \rho T \left ( \bm{I} + \bm{\bar{A}}\right )\right ) \,=\, \bm{0}_{\RR^{3}}, \\
    & \partial_t E + \diverg_x \left (\frac 12 \rho |\bm{u}|^2 \bm{u} + \rho T \left (  \frac{5}{2}\bm{I} + \bm{\bar{A}} \right ) \bm{u} + \rho (T)^{3/2} \bm{\bar{B}} \right ) \,=\, 0,
    \end{aligned} \right.
\end{equation}
where the traceless matrix $\bm{\bar{A}} \in M_{3}$ and the vector $\bm{\bar{B}} \in \RR^{3}$ are defined as:
\begin{equation}
    \label{eqDefABV_App}
    \left\{\begin{aligned}
    & \bm{\bar{A}} := \frac{1}{\rho} \int_{\RR^{3}} \bm{A}(\bm{V}) f(v) \, dv, && \bm{A}(\bm{V}) = \bm{V} \otimes \bm{V} - \frac{|\bm{V}|^2}{3} \bm{I}, \\
    & \bm{\bar{B}} := \frac{1}{\rho}\int_{\RR^{3}} \bm{B}(\bm{V}) f(v) \, dv, && \bm{B}(\bm{V}) = \frac12 \left [ |\bm{V}|^2 - 5\right ] \bm{V}, 
    \end{aligned}\right.
\end{equation}
with the normalized velocity variable $\bm{V}(v) = (v-\bm{u})/\sqrt{T}$.

\subsection{Zeroth-order approximation (Euler)}
At zeroth order in $\ve$, substituting the Maxwellian distribution $f = \mathcal{M}_{\rho, \bm{u}, T}$ directly into \eqref{eqDefABV_App} yields:
\begin{equation}
    \bm{\bar{A}}_{Euler} = \bm{0}_{M_{3}}, \quad \bm{\bar{B}}_{Euler} = \bm{0}_{\RR^{3}},
\end{equation}
recovering the compressible Euler equations.

\subsection{First-order approximation (Navier--Stokes)}
At first order, we plug $f = \mathcal M_{\rho,\bm{u},T}\left[1 + \ve \, g^{(1)}\right]$ into the Boltzmann equation. Since the Maxwellian distribution is an equilibrium point, the fluctuation $g^{(1)}$ satisfies:
\begin{equation}
    \label{eqCollisionExp_App}
    \partial_t \mathcal{M} + v \cdot \nabla_x \mathcal{M} = \LL_{\mathcal{M}} g^{(1)} + \mathcal{O}(\ve),
\end{equation}
where we write $\mathcal{M}$ as a shorthand for $\mathcal{M}_{\rho,\bm{u},T}$ and $\LL_{\mathcal M}$ is the linearized collision operator. Using the identities \cite{lev-1996, Golse:2005, sai-2009}:
\begin{equation}
    \label{eq::evolution_Maxwellian_App}
     \frac{\partial \mathcal{M}}{\partial t} + v \cdot \nabla_x \mathcal{M} = \mathcal{M}\left[\frac{1}{T}\left(c\otimes c - \frac{1}{3}|c|^2\bm{I}\right):\nabla_x \bm u + \left(\frac{|c|^2}{2T}-\frac{5}{2}\right)c\cdot \frac{\nabla_x T}{T}\right],
\end{equation}
with $c=v-\bm u$, and substituting the definitions \eqref{eqDefABV_App}, we obtain \cite{lev-1996, lev-mor-nad-1998}:
\begin{equation}
    \label{eqLinearG1_App}
    \LL_{\mathcal M} g^{(1)} = \mathcal M\left [ \frac{1}{2} \bm{A}(\bm{V}):\bm{D}(\bm{u}) + 2 \bm{B}(\bm{V}) \cdot \nabla_x \sqrt{T}\right ],
\end{equation}
where $\bm{D}(\bm{u}) = \nabla_x \bm{u} + \left (\nabla_x \bm{u}\right )^\intercal - \frac{2}{3} \left (\diverg_x \bm{u}\right ) \bm{I}$ is the deformation tensor. Inverting the linearized collision operator yields:
\begin{equation}
    \label{eqG1_App}
    g^{(1)} = \frac{1}{2} \LL_{\mathcal M}^{-1}\left ( \mathcal M \, \bm{A} \right ): \bm{D}(\bm{u}) + 2 \,\LL_{\mathcal M}^{-1}\left ( \mathcal M \, \bm{B} \right ) \cdot \nabla_x \sqrt{T}. 
\end{equation}
Plugging this fluctuation back into the definitions of $\bm{\bar{A}}$ and $\bm{\bar{B}}$ in \eqref{eqDefABV_App} leads to the standard Navier--Stokes transport relations:
\begin{equation}
    \label{eqABNS_App}
    \left\{\begin{aligned}
      & \bm{\bar{A}}_{NS} = - \ve \frac{\mu}{\rho \, T} \bm{D}(\bm{u}), \\ 
      & \bm{\bar{B}}_{NS} = - \ve \frac{\kappa}{\rho \,T^{3/2}} \nabla_x T,
    \end{aligned}\right.
\end{equation}
where the dynamic viscosity $\mu$ and thermal conductivity $\kappa$ are given by:
\begin{gather*}
    \mu := -\frac{T}{10} \int_{\RR^{3}} \bm{A}(\bm{V}): \mathcal M(v) \LL_{\mathcal M}^{-1}\left ( \mathcal M \, \bm{A} \right )(v) \, dv, \\
    \kappa := -\frac{T}{3} \int_{\RR^{3}} \bm{B}(\bm{V}) \cdot \mathcal M(v) \LL_{\mathcal M}^{-1}\left ( \mathcal M \, \bm{B} \right )(v) \, dv.
\end{gather*}

\subsection{Second-order approximation (Burnett Equations)}
By extending the expansion to second order $f=\mathcal{M}[1+\ve g^{(1)}+\ve^2 g^{(2)}]$, we recover the Burnett coefficients. For the ES-BGK model, these expressions are analytically given by \cite{zhe-str-2004}:
\begin{equation}\label{eqAepsBurnett_App}
\begin{split}
    \bm{\bar{A}}_{Burnett} &=-\frac{\ve}{\rho T}\mu \bm D(\bm u) + \frac{\ve^2}{\rho T}\frac{1}{(1-\beta)^2\nu^2}
    \Bigg[
    -2T^2\left({\rm Hess}_x(\rho)-\frac{1}{3}(\nabla_x^2 \rho)\bm I\right)\\
    &-2\beta\rho T \left({\rm Hess}_x(T)-\frac{1}{3}(\nabla_x^2 T)\bm I\right)
    +\frac{2T^2}{\rho}\left(\nabla_x\rho\otimes \nabla_x\rho - \frac{1}{3}|\nabla_x \rho|^2 \bm I\right)\\ 
    &-2T\left( \frac{1}{2}(\nabla_x T\otimes\nabla_x\rho + \nabla_x\rho\otimes\nabla_x T)-\frac{1}{3}(\nabla_x T)\cdot(\nabla_x \rho)\bm I\right)\\
    &+2P\left( (\nabla_x \bm u)(\nabla_x \bm u)^\intercal-\frac{1}{3}(\nabla_x \bm u : \nabla_x \bm u)\bm I\right)\\
    &+\frac{2-4\alpha}{3}P \frac{1}{2} \bm D(\bm u) (\nabla_x \cdot \bm u)
    +2\alpha (1-\beta)\rho \left(\nabla_x T \otimes \nabla_x T - \frac{1}{3}|\nabla_x T|^2 \bm I\right)
    \Bigg],
\end{split}
\end{equation}

\begin{equation}\label{eqBepsBurnett_App}
\begin{split}
    \bm{\bar{B}}_{Burnett} &=-\frac{\ve}{\rho T^{\frac{3}{2}}}
    \left(\kappa \nabla_x T\right)
    + \frac{\ve^2}{\rho T^{\frac{3}{2}}}\frac{1}{(1-\beta)\nu^2}
    \Bigg[
    \frac{5\beta-4}{3}P T \nabla_x(\nabla_x \cdot \bm u)
    +PT(\nabla_x^2 \bm u)\\
    &-T^2\bm D(\bm u)\nabla_x \rho
    +\left(6+\alpha-\frac{7}{2}\beta\right)P(\nabla_x \bm u)(\nabla_x T)\\
    &+\left(1+\alpha+\frac{3}{2}\beta\right)P(\nabla_x \bm u)^\intercal (\nabla_x T)\\
    &+\frac{(-13-\beta)+(14-10\beta)(1-\alpha)}{6}P(\nabla_x \cdot \bm u)(\nabla_x T)
    \Bigg],
\end{split}
\end{equation}
where $P=\rho T$, $\beta$ is the free parameter of the ES-BGK operator, and $\alpha$ is the exponent of the temperature-dependent viscosity law $\mu(T)=\mu_0 (T/T_0)^\alpha$.

\section{Moment realizability matrix and Schur complement}
\label{app:Realizability}

In a three-dimensional velocity space, we define the vector of collisional invariants and fluctuations:
\begin{equation}
    \bm{m} =
    \begin{pmatrix}
     1\\
     \frac{v- \bm u}{\sqrt{T}}\\
     \frac{|v-\bm u|^2}{2T}-\frac{5}{2}
    \end{pmatrix},
\end{equation}
and define the moment realizability matrix $\bm{M}$ as:
\begin{equation}
\label{eq:MRM_App}
    \bm{M} = \frac{1}{\rho}\int (\bm{m} \otimes \bm{m}) f dv =
    \begin{pmatrix}
     1  & \bm{0}^\intercal            & -1\\
     \bm{0}  & \bm{I}+\bm{\bar{A}}    & \bm{\bar{B}}\\
     -1 & \bm{\bar{B}}^\intercal & \bar{C}
    \end{pmatrix},
\end{equation}
where $\bm{\bar{A}}$ and $\bm{\bar{B}}$ are given in \eqref{eqDefABV_App}, and the scalar $\bar{C}$ is:
\begin{equation}
    \bar{C} = \frac{1}{\rho}\int\left(\frac{|v-\bm u|^2}{2T}-\frac{5}{2}\right)^2 f(v) \, dv = \frac{1}{4\rho T^2}\int|v-\bm u|^4 f \, dv - \frac{5}{4}.
\end{equation}

\subsection{Proof of the lower bound for $\bar{C}$}
To show that $\bar{C} > 1$, we employ Hölder's inequality:
\begin{equation}
\label{eq::Holder_App}
   \int\left|g(v)h(v)\right|dv \leq \left(\int|g(v)|^{2} dv\right)^{1/2} \left(\int|h(v)|^{2}dv\right)^{1/2},  
\end{equation}
setting $h(v)=\sqrt{f}$ and $g(v)=\sqrt{f}|v-\bm u|^2$. Assuming a standard non-degeneracy condition on $f$ (in both continuous and discrete settings), $g$ and $h$ are linearly independent, ensuring that the inequality holds strictly \cite{mit-1970}. This yields:
\begin{equation}
    \left( \int f |v-\bm{u}|^2 \, dv \right)^2 < \left( \int f |v-\bm{u}|^4 \, dv \right) \left( \int f \, dv \right).
\end{equation}
Recalling that $\int f \, dv = \rho$ and $\int f |v-\bm{u}|^2 \, dv = 3\rho T$, we get:
\begin{equation}
    (3\rho T)^2 < \rho \int f |v-\bm{u}|^4 \, dv \implies \int f |v-\bm{u}|^4 \, dv > 9\rho T^2.
\end{equation}
Substituting this back into the expression for $\bar{C}$:
\begin{equation}
    \bar{C} > \frac{9 \rho T^2}{4 \rho T^2} - \frac{5}{4} = \frac{9}{4} - \frac{5}{4} = 1.
\end{equation}
This confirms that $\bar{C} > 1$ for any non-negative distribution function $f$ satisfying the non-degeneracy assumption above.

\subsection{Two-stage Schur complement derivation}
Since $\bar{C}>1$, the positive semi-definiteness of $\bm{M}$ can be characterized through successive Schur complements. First, taking the Schur complement of $\bm{M}$ with respect to its last diagonal entry $\bar C>0$ yields the Schur complement $\tilde{\mathcal{S}}$:
\begin{equation}
    \tilde{\mathcal{S}} =
    \begin{pmatrix}
     1-\displaystyle \frac{1}{\bar{C}}     & \displaystyle \frac{\bm{\bar{B}}^\intercal}{\bar{C}} \\[8pt]
     \displaystyle \frac{\bm{\bar{B}}}{\bar{C}} & \bm{I}+\bm{\bar{A}}- \displaystyle\frac{\bm{\bar{B}}\bm{\bar{B}}^\intercal}{\bar{C}}
    \end{pmatrix}.
\end{equation}
Since $1 - 1/\bar{C} > 0$, we can compute the Schur complement of $\tilde{\mathcal{S}}$ with respect to its top-left scalar element. This yields the final reduced matrix $\mathcal{S}$:
\begin{equation}
    \mathcal{S} = \left(\bm{I}+\bm{\bar{A}}- \frac{\bm{\bar{B}}\bm{\bar{B}}^\intercal}{\bar{C}}\right) - \left(\frac{\bm{\bar{B}}}{\bar{C}}\right)\left(1 - \frac{1}{\bar{C}}\right)^{-1}\left(\frac{\bm{\bar{B}}^\intercal}{\bar{C}}\right),
\end{equation}
which simplifies to:
\begin{equation}
    \mathcal{S} = \bm{I}+\bm{\bar{A}} - \displaystyle\frac{\bm{\bar{B}}\bm{\bar{B}}^\intercal}{\bar{C}-1}.
\end{equation}
Thus, the positive semi-definiteness of $\bm{M}$ is mathematically equivalent to the positive semi-definiteness of $\mathcal{S}$.

\printbibliography

@article{bos-par-rus-2017,
  author  = {Sebastiano Boscarino and Lorenzo Pareschi and Giovanni Russo},
  title   = {A unified {IMEX} Runge--Kutta approach for hyperbolic systems
             with multiscale relaxation},
  journal = {SIAM Journal on Numerical Analysis},
  volume  = {55},
  number  = {4},
  pages   = {2085--2109},
  year    = {2017},
  doi     = {10.1137/M1111449}
}

@book{bos-par-rus-2024,
  author    = {Sebastiano Boscarino and Lorenzo Pareschi and Giovanni Russo},
  title     = {Implicit-Explicit Methods for Evolutionary Partial Differential
               Equations},
  series    = {Mathematical Modeling and Computation},
  publisher = {Society for Industrial and Applied Mathematics},
  address   = {Philadelphia, PA},
  year      = {2024},
  isbn      = {978-1-61197-819-3},
  doi       = {10.1137/1.9781611978209}
}

@article {Klar98,
    AUTHOR = {Tiwari, S. and Klar, A.},
     TITLE = {{An adaptive domain decomposition procedure for {B}oltzmann and
              {E}uler equations}},
   JOURNAL = {J. Comput. Appl. Math.},
  FJOURNAL = {Journal of Computational and Applied Mathematics},
    VOLUME = {90},
      YEAR = {1998},
    NUMBER = {2},
     PAGES = {223--237},
      ISSN = {0377-0427,1879-1778},
   MRCLASS = {76M25 (65M55 76P05)},
  MRNUMBER = {1624342},
       DOI = {10.1016/S0377-0427(98)00027-2},
       URL = {https://doi.org/10.1016/S0377-0427(98)00027-2},
}

@article {DiPa13,
    AUTHOR = {Dimarco, Giacomo and Pareschi, Lorenzo},
     TITLE = {{Asymptotic preserving implicit-explicit {R}unge-{K}utta methods for nonlinear kinetic equations}},
   JOURNAL = {SIAM J. Numer. Anal.},
  FJOURNAL = {SIAM Journal on Numerical Analysis},
    VOLUME = {51},
      YEAR = {2013},
    NUMBER = {2},
     PAGES = {1064--1087},
      ISSN = {0036-1429,1095-7170},
   MRCLASS = {65L06 (65L04)},
  MRNUMBER = {3036998},
MRREVIEWER = {Amit\ Kumar\ Verma},
       DOI = {10.1137/12087606X},
       URL = {https://doi.org/10.1137/12087606X},
}

@InCollection{DiPa08,
  author    = {Dimarco, G. and Pareschi, L.},
  booktitle = {Hyperbolic Problems: Theory, Numerics, Applications},
  publisher = {Springer, Berlin},
  title     = {{Domain decomposition techniques and hybrid multiscale methods for kinetic equations}},
  year      = {2008},
  isbn      = {978-3-540-75711-5},
  pages     = {457--464},
  doi       = {10.1007/978-3-540-75712-2\_43},
  mrclass   = {65M55 (35B25 35L60 82C80)},
  mrnumber  = {2549178},
  url       = {https://doi.org/10.1007/978-3-540-75712-2_43},
}

@book{can-hus-qua-zan-1991,
  Title                    = {{Spectral Methods in Fluid Dynamics}},
  Author                   = {C. Canuto and M.Y. Hussaini and A. Quarteroni and T.A. Zang},
  Publisher                = {Springer-Verlag},
  Year                     = {1988},

  Address                  = {New York},
  Series                   = {Springer Series in Computational Physics},

  ISBN                     = {0-387-17371-4},
  Pages                    = {xiv+557}
}

@InCollection{Golse:2005,
  Title                    = {{The Boltzmann equation and its hydrodynamic limits}},
  Author                   = {Golse, Fran\c{c}ois},
  Booktitle                = {Handbook of Differential Equations: Evolutionary Equations Vol. 2},
  Publisher                = {North-Holland},
  Year                     = {2005},
  Editor                   = {Dafermos, Constantin and Feireisl, Eduard},
  Pages                    = {159--303},

  ISBN                     = {0444520481}
}

@phdthesis{rey-2023,
  TITLE = {{Contributions to the Mathematical and Numerical Analysis of Multiscale Kinetic Equations}},
  AUTHOR = {Rey, Thomas},
  URL = {https://theses.hal.science/tel-04032898},
  NOTE = {128 pages},
  SCHOOL = {{Universit{\'e} de Lille}},
  YEAR = {2023},
  MONTH = Feb,
  TYPE = {Habilitation {\`a} diriger des recherches},
  HAL_ID = {tel-04032898},
  HAL_VERSION = {v1},
}

@article{lai-2023,
title = {Hybrid kinetic/fluid numerical method for the Vlasov-BGK equation in the diffusive scaling},
journal = {Kinetic and Related Models},
volume = {16},
number = {6},
pages = {913-947},
year = {2023},
issn = {1937-5093},
doi = {10.3934/krm.2023013},
url = {https://www.aimsciences.org/article/id/641ada1d6565307f274b58e4},
author = {Tino Laidin}
}

@InProceedings{lai-rey-2023,
  author    = {Laidin, Tino and Rey, Thomas},
  booktitle = {Finite Volumes for Complex Applications X---Volume 2, Hyperbolic and Related Problems},
  title     = {{Hybrid Kinetic/Fluid Numerical Method for the Vlasov-Poisson-BGK Equation in the Diffusive Scaling}},
  year      = {2023},
  address   = {Cham},
  editor    = {Franck, Emmanuel and Fuhrmann, J{\"u}rgen and Michel-Dansac, Victor and Navoret, Laurent},
  pages     = {229--237},
  publisher = {Springer Nature Switzerland},
  isbn      = {978-3-031-40860-1},
}

@article{fil-mou-par-2006,
  Title                    = {{Solving the Boltzmann Equation in N log2 N}},
  Author                   = {F. Filbet and C. Mouhot and L. Pareschi},
  Journal                  = {SIAM J. Sci. Comput.},
  Year                     = {2006},
  Number                   = {3},
  Pages                    = {1029--1053},
  Volume                   = {28},

  Publisher                = {Society for Industrial and Applied Mathematics}
}

@Article{bha-gro-kro-1954,
  author    = {Bhatnagar, P. L. and Gross, E. P. and Krook, M.},
  journal   = {Phys. Rev.},
  title     = {{A Model for Collision Processes in Gases. I. Small Amplitude Processes in Charged and Neutral One-Component Systems}},
  year      = {1954},
  month     = {May},
  pages     = {511--525},
  volume    = {94},
  doi       = {10.1103/PhysRev.94.511},
  number     = {3},
  publisher = {American Physical Society},
  url       = {https://link.aps.org/doi/10.1103/PhysRev.94.511},
}

@article{boy-che-can-1995,
    author = {Boyd, Iain D. and Chen, Gang and Candler, Graham V.},
    title = {{Predicting failure of the continuum fluid equations in transitional hypersonic flows}},
    journal = {Phys. Fluids},
    volume = {7},
    number = {1},
    pages = {210-219},
    year = {1995},
    month = {01},
    issn = {1070-6631},
    doi = {10.1063/1.868720},
    url = {https://doi.org/10.1063/1.868720},
    eprint = {https://pubs.aip.org/aip/pof/article-pdf/7/1/210/19222572/210\_1\_online.pdf},
}

@article{and-let-per-per-2000,
  Title                    = {{The Gaussian-BGK model of Boltzmann equation with small Prandtl number}},
  Author                   = {Andries, Pierre and {Le Tallec}, Patrick and Perlat, Jean-Philippe and Perthame, Benoit},
  Journal                  = {Eur. J. Mech. B Fluids},
  Year                     = {2000},

  Month                    = nov,
  Number                   = {6},
  Pages                    = {813--830},
  Volume                   = {19},

  Doi                      = {10.1016/S0997-7546(00)01103-1},
  ISSN                     = {09977546},
  Url                      = {http://linkinghub.elsevier.com/retrieve/pii/S0997754600011031}
}

@article{tiw-1998,
  author    = {Tiwari, S.},
  journal   = {J. Comput. Phys.},
  title     = {{Coupling of the Boltzmann and Euler equations with automatic domain decomposition}},
  year      = {1998},
  issn      = {00219991},
  month     = aug,
  number    = {2},
  pages     = {710--726},
  volume    = {144},
  doi       = {10.1006/jcph.1998.6011},
  url       = {http://linkinghub.elsevier.com/retrieve/pii/S0021999198960119 http://www.sciencedirect.com/science/article/pii/S0021999198960119},
}

@article{lev-mor-nad-1998,
  Title                    = {{Moment realizability and the validity of the {N}avier-{S}tokes equations for rarefied gas dynamics}},
  Author                   = {Levermore, C. David and Morokoff, William J. and Nadiga, B. T.},
  Journal                  = {Phys. Fluids},
  Year                     = {1998},
  Number                   = {12},
  Pages                    = {3214--3226},
  Volume                   = {10},

  Coden                    = {PHFLE6},
  Doi                      = {10.1063/1.869849},
  Fjournal                 = {Physics of Fluids},
  ISSN                     = {1070-6631},
  Mrclass                  = {76N15 (76M25 76P05)},
  Mrnumber                 = {1657273 (2000b:76088)},
  Mrreviewer               = {Hans Babovsky},
  Url                      = {http://dx.doi.org/10.1063/1.869849}
}

@Article{tiw-2000,
  author  = {{Tiwari}, S.},
  journal = {Transport Theor. Stat.},
  title   = {{Application of Moment Realizability Criteria for the Coupling of the Boltzmann and Euler Equations}},
  year    = {2000},
  month   = nov,
  number  = {7},
  pages   = {759-783},
  volume  = {29},
  adsurl  = {https://ui.adsabs.harvard.edu/abs/2000TTSP...29..759T},
  doi     = {10.1080/00411450008200001},
}

@article{fil-rey-2014,
  author    = {Filbet, Francis and Rey, Thomas},
  journal   = {SIAM J. Sci. Comput.},
  title     = {{A Hierarchy of Hybrid Numerical Methods for Multi-Scale Kinetic Equations}},
  year      = {2015},
  number    = {3},
  pages     = {A1218-A1247},
  volume    = {37},
  doi       = {10.1137/140958773},
}

@Article{PaRu:SINUM:2000,
  Title                    = {{Numerical Solution of the Boltzmann Equation I : Spectrally Accurate Approximation of the Collision Operator}},
  Author                   = {L. Pareschi and G. Russo},
  Journal                  = {SIAM J. Numer. Anal.},
  Year                     = {2000},
  Number                   = {4},
  Pages                    = {1217--1245},
  Volume                   = {37}
}

@Article{MoPa:2006,
  author     = {Mouhot, Cl{\'e}ment and Pareschi, Lorenzo},
  journal    = {Math. Comput.},
  title      = {{Fast algorithms for computing the {B}oltzmann collision operator}},
  year       = {2006},
  issn       = {0025-5718},
  number     = {256},
  pages      = {1833--1852 (electronic)},
  volume     = {75},
  coden      = {MCMPAF},
  doi        = {10.1090/S0025-5718-06-01874-6},
  fjournal   = {Mathematics of Computation},
  mrclass    = {65M70 (65T50 76M22 76P05 82C40 82C70 82C80)},
  mrnumber   = {2240637 (2007d:65095)},
  mrreviewer = {Stephen Wollman},
  url        = {http://dx.doi.org/10.1090/S0025-5718-06-01874-6},
}

@article{fil-xio-2018,
  author    = {Filbet, Francis and Xiong, Tao},
  journal   = {J. Comput. Phys.},
  title     = {{A hybrid discontinuous Galerkin scheme for multi-scale kinetic equations}},
  year      = {2018},
  pages     = {841--863},
  volume    = {372},
  publisher = {Elsevier},
}

@book{sai-2009,
  author    = {Saint-Raymond, L.},
  publisher = {Springer},
  title     = {{Hydrodynamic Limits of the Boltzmann Equation}},
  year      = {2009},
  isbn      = {9783540928461},
  lccn      = {2008943981},
  url       = {https://books.google.fr/books?id=ROUILXXb7UUC},
}

@article{dim-par-2014,
  author    = {Dimarco, G. and Pareschi, L.},
  journal   = {Acta Numer.},
  title     = {{Numerical methods for kinetic equations}},
  year      = {2014},
  pages     = {369-520},
  volume    = {23},
}

@Article{lev-1996,
  author  = {Levermore, C.},
  journal = {J. Stat. Phys.},
  title   = {{Moment closure hierarchies for kinetic theories}},
  year    = {1996},
  month   = {06},
  pages   = {1021-1065},
  volume  = {83},
  doi     = {10.1007/BF02179552},
}

@InBook{par-rus-2001,
  author  = {Pareschi, Lorenzo and Russo, Giovanni},
  pages   = {269-288},
  title   = {{Implicit–Explicit Runge–Kutta schemes for stiff systems of differential equations}},
  publisher = {Nova Sci. Publ.},
  year    = {2001},
  month   = {01},
  volume  = {3},
  journal = {Recent Trends in Numerical Analysis},
}

@Article{fil-jin-2011,
  author  = {Filbet, Francis and Jin, Shi},
  journal = {J. Sci. Comput.},
  title   = {{An Asymptotic Preserving Scheme for the ES-BGK Model of the Boltzmann Equation}},
  year    = {2011},
  month   = {02},
  pages   = {204-224},
  volume  = {46},
  doi     = {10.1007/s10915-010-9394-x},
}

@Article{fil-jin-2010,
  author   = {Francis Filbet and Shi Jin},
  journal  = {J. Comput. Phys.},
  title    = {{A class of asymptotic-preserving schemes for kinetic equations and related problems with stiff sources}},
  year     = {2010},
  issn     = {0021-9991},
  number   = {20},
  pages    = {7625-7648},
  volume   = {229},
  doi      = {https://doi.org/10.1016/j.jcp.2010.06.017},
  url      = {https://www.sciencedirect.com/science/article/pii/S0021999110003323},
}

@Article{cor-per-1991,
  author  = {François Coron and Beno{\^i}t Perthame},
  journal = {SIAM J. Numer. Anal.},
  title   = {{Numerical passage from kinetic to fluid equations}},
  year    = {1991},
  pages   = {26-42},
  volume  = {28},
  url     = {https://api.semanticscholar.org/CorpusID:121757327},
}

@Article{jin-1995,
  author  = {Shi Jin},
  journal = {J. Comput. Phys.},
  title   = {{Runge-Kutta Methods for Hyperbolic Conservation Laws with Stiff Relaxation Terms}},
  year    = {1995},
  issn    = {0021-9991},
  number  = {1},
  pages   = {51-67},
  volume  = {122},
  doi     = {https://doi.org/10.1006/jcph.1995.1196},
  url     = {https://www.sciencedirect.com/science/article/pii/S0021999185711965},
}

@Article{caf-shi-rus-1995,
  author  = {Caflisch, Russel and Jin, Shi and Russo, Giovanni},
  journal = {SIAM J. Numer. Anal.},
  title   = {{Uniformly Accurate Schemes For Hyperbolic Systems With Relaxation}},
  year    = {1997},
  month   = {2},
  volume  = {34},
  doi     = {10.1137/S0036142994268090},
}

@Article{jin-par-tos-1998,
  author  = {Jin, Shi and Pareschi, Lorenzo and Toscani, Giuseppe},
  journal = {SIAM J. Numer. Anal.},
  title   = {{Diffusive Relaxation Schemes For Multiscale Discrete-Velocity Kinetic Equations}},
  year    = {1998},
  month   = {12},
  pages   = {2405-2439},
  volume  = {35},
  doi     = {10.1137/S0036142997315962},
}

@Article{jin-1999,
  author     = {Jin, Shi},
  journal    = {SIAM J. Sci. Comput.},
  title      = {{Efficient Asymptotic-Preserving (AP) Schemes For Some Multiscale Kinetic Equations}},
  year       = {1999},
  issn       = {1064-8275},
  month      = jan,
  number     = {2},
  pages      = {441–454},
  volume     = {21},
  address    = {USA},
  doi        = {10.1137/S1064827598334599},
  issue_date = {1999},
  numpages   = {14},
  publisher  = {Society for Industrial and Applied Mathematics},
  url        = {https://doi.org/10.1137/S1064827598334599},
}

@InCollection{jin-par-2001,
  author    = {Jin, Shi and Pareschi, Lorenzo},
  booktitle = {Hyperbolic Problems: Theory, Numerics, Applications},
  publisher = {Birkh{\"a}user Basel},
  title     = {{Asymptotic-Preserving (Ap) Schemes for Multiscale Kinetic Equations: a Unified Approach}},
  year      = {2001},
  address   = {Basel},
  editor    = {Freist{\"u}hler, Heinrich and Warnecke, Gerald},
  isbn      = {978-3-0348-8372-6},
  pages     = {573--582},
}

@Article{klar-1998,
  author  = {Klar, Axel},
  journal = {SIAM J. Numer. Anal.},
  title   = {{An Asymptotic-Induced Scheme for Nonstationary Transport Equations in the Diffusive Limit}},
  year    = {1998},
  number  = {3},
  pages   = {1073-1094},
  volume  = {35},
  doi     = {10.1137/S0036142996305558},
  eprint  = {https://doi.org/10.1137/S0036142996305558},
  url     = {https://doi.org/10.1137/S0036142996305558},
}

@article{gos-tos-2002,
title = {An asymptotic-preserving well-balanced scheme for the hyperbolic heat equations},
journal = {Comptes Rendus Mathematique},
volume = {334},
number = {4},
pages = {337-342},
year = {2002},
issn = {1631-073X},
doi = {https://doi.org/10.1016/S1631-073X(02)02257-4},
url = {https://www.sciencedirect.com/science/article/pii/S1631073X02022574},
author = {Laurent Gosse and Giuseppe Toscani},
}

@Article{kol-ars-ari-fro-zab-2007,
  author   = {V.I. Kolobov and R.R. Arslanbekov and V.V. Aristov and A.A. Frolova and S.A. Zabelok},
  journal  = {J. Comput. Phys.},
  title    = {{Unified solver for rarefied and continuum flows with adaptive mesh and algorithm refinement}},
  year     = {2007},
  issn     = {0021-9991},
  number   = {2},
  pages    = {589-608},
  volume   = {223},
  doi      = {https://doi.org/10.1016/j.jcp.2006.09.021},
  url      = {https://www.sciencedirect.com/science/article/pii/S0021999106004633},
}

@Article{deg-dim-2012,
  author   = {Pierre Degond and Giacomo Dimarco},
  journal  = {J. Comput. Phys.},
  title    = {{Fluid simulations with localized boltzmann upscaling by direct simulation Monte-Carlo}},
  year     = {2012},
  issn     = {0021-9991},
  number   = {6},
  pages    = {2414-2437},
  volume   = {231},
  doi      = {https://doi.org/10.1016/j.jcp.2011.11.030},
  url      = {https://www.sciencedirect.com/science/article/pii/S0021999111006851},
}

@Article{tiw-kla-har-2009,
  author   = {Sudarshan Tiwari and Axel Klar and Steffen Hardt},
  journal  = {J. Comput. Phys.},
  title    = {{A particle–particle hybrid method for kinetic and continuum equations}},
  year     = {2009},
  issn     = {0021-9991},
  number   = {18},
  pages    = {7109-7124},
  volume   = {228},
  doi      = {https://doi.org/10.1016/j.jcp.2009.06.019},
  url      = {https://www.sciencedirect.com/science/article/pii/S0021999109003428},
}

@Article{ala-pup-2012,
  author   = {Alessandro Alaia and Gabriella Puppo},
  journal  = {J. Comput. Phys.},
  title    = {{A hybrid method for hydrodynamic-kinetic flow – Part II – Coupling of hydrodynamic and kinetic models}},
  year     = {2012},
  issn     = {0021-9991},
  number   = {16},
  pages    = {5217-5242},
  volume   = {231},
  doi      = {https://doi.org/10.1016/j.jcp.2012.02.022},
  url      = {https://www.sciencedirect.com/science/article/pii/S0021999112001234},
}

@Article{ala-pup-2011,
  author   = {Alessandro Alaia and Gabriella Puppo},
  journal  = {J. Comput. Phys.},
  title    = {{A hybrid method for hydrodynamic-kinetic flow Part I: A particle-grid method for reducing stochastic noise in kinetic regimes}},
  year     = {2011},
  issn     = {0021-9991},
  number   = {14},
  pages    = {5660-5683},
  volume   = {230},
  doi      = {https://doi.org/10.1016/j.jcp.2011.03.047},
  url      = {https://www.sciencedirect.com/science/article/pii/S0021999111002087},
}

@Article{dim-mie-ris-2014,
  author   = {Giacomo Dimarco and Luc Mieussens and Vittorio Rispoli},
  journal  = {J. Comput. Phys.},
  title    = {{An asymptotic preserving automatic domain decomposition method for the Vlasov–Poisson–BGK system with applications to plasmas}},
  year     = {2014},
  issn     = {0021-9991},
  pages    = {122-139},
  volume   = {274},
  doi      = {https://doi.org/10.1016/j.jcp.2014.06.002},
  url      = {https://www.sciencedirect.com/science/article/pii/S0021999114004069},
}

@article{li-son-wan-2021,
author = {Li, Weiming and Song, Peng and Wang, Yanli},
title = {{A Hybrid Moment Method for Multi-scale Kinetic Equations Based on Maximum Entropy Principle}},
year = {2021},
volume = {89},
number = {2},
issn = {0885-7474},
url = {https://doi.org/10.1007/s10915-021-01639-0},
doi = {10.1007/s10915-021-01639-0},
journal = {J. Sci. Comput.}
}

@Article{xio-qiu-2017,
  author   = {Tao Xiong and Jing-Mei Qiu},
  journal  = {{J. Comput. Phys.}},
  title    = {{A hierarchical uniformly high order DG-IMEX scheme for the 1D BGK equation}},
  year     = {2017},
  issn     = {0021-9991},
  pages    = {164-191},
  volume   = {336},
  doi      = {https://doi.org/10.1016/j.jcp.2017.01.032},
  url      = {https://www.sciencedirect.com/science/article/pii/S0021999117300487},
}

@Article{tia-mar-2021,
  author  = {Xiao, Tianbai and Frank, Martin},
  journal = {J. Comput. Phys.},
  title   = {{Using neural networks to accelerate the solution of the Boltzmann equation}},
  year    = {2021},
  month   = {10},
  pages   = {110521},
  volume  = {443},
  doi     = {10.1016/j.jcp.2021.110521},
}

@Article{gan-tap-esm-myo-2024,
  author  = {Garg, Gagan and Mankodi, Tapan and Esmaeilifar, Esmaeil and Myong, R. S.},
  journal = {Phys. Fluids},
  title   = {{Neural network-based finite volume method and direct simulation Monte Carlo solutions of non-equilibrium shock flow guided by nonlinear coupled constitutive relations}},
  year    = {2024},
  month   = {10},
  pages   = {106113},
  volume  = {36},
  doi     = {10.1063/5.0223654},
}

@article{zha-kun-2021,
author = {Guo, Zhaoli and Xu, Kun},
year = {2021},
month = {},
pages = {},
title = {{Progress of discrete unified gas-kinetic scheme for multiscale flows}},
volume = {3},
journal = {{Advances in Aerodynamics}},
doi = {10.1186/s42774-020-00058-3}
}

@Article{zha-wan-xu-2014,
  author  = {Guo, Zhaoli and Wang, Ruijie and Xu, Kun},
  journal = {Phys. Rev. E},
  title   = {{Discrete unified gas kinetic scheme for all Knudsen number flows: II. Thermal compressible case}},
  year    = {2015},
  volume  = {91},
  doi     = {10.1103/PhysRevE.91.033313},
}

@Article{guo-xu-wan-2013,
  author  = {Guo, Zhaoli and Xu, Kun and Wang, Ruijie},
  journal = {Phys. Rev. E},
  title   = {{Discrete unified gas kinetic scheme for all Knudsen number flows: Low-speed isothermal case.}},
  year    = {2013},
  month   = {09},
  pages   = {033305},
  volume  = {88},
}

@article{wu-shi-cha-wan-2016,
title = {{Discrete unified gas kinetic scheme with a force term for incompressible fluid flows}},
journal = {{Computers \& Mathematics with Applications}},
volume = {71},
number = {12},
pages = {2608-2629},
year = {2016},
issn = {0898-1221},
doi = {https://doi.org/10.1016/j.camwa.2016.04.025},
url = {https://www.sciencedirect.com/science/article/pii/S0898122116302061},
author = {Chen Wu and Baochang Shi and Zhenhua Chai and Peng Wang},
}

@Unpublished{lee-sch-xia-kru-fra-2024,
  author        = {Jae Yong Lee and Steffen Schotthöfer and Tianbai Xiao and Sebastian Krumscheid and Martin Frank},
  note          = {arXiv preprint arXiv:2402.16613},
  title         = {{Structure-Preserving Operator Learning: Modeling the Collision Operator of Kinetic Equations}},
  year          = {2024},
  archiveprefix = {{arXiv}},
  eprint        = {2402.16613},
  primaryclass  = {math.NA},
  url           = {https://arxiv.org/abs/2402.16613},
}

@article{cap-par-rey-2025,
  title={Hierarchical dynamic domain decomposition for the multiscale Boltzmann equation},
  author={Caparello, Domenico and Pareschi, Lorenzo and Rey, Thomas},
  journal={arXiv preprint arXiv:2505.03360},
  year={2025}
}

@article{bos-dim-2021,
 author = {Boscheri, Walter and Dimarco, Giacomo},
 title = {High order finite volume schemes with {IMEX} time stepping for the {Boltzmann} model on unstructured meshes},
 fjournal = {Computer Methods in Applied Mechanics and Engineering},
 journal = {Comput. Methods Appl. Mech. Eng.},
 issn = {0045-7825},
 volume = {387},
 year = {2021},
 language = {English},
 doi = {10.1016/j.cma.2021.114180},
 zbMATH = {7427413},
 Zbl = {1507.65147}
}

@article{bos-dim-2020,
 author = {Boscheri, Walter and Dimarco, Giacomo},
 title = {High order central {WENO}-implicit-explicit {Runge} {Kutta} schemes for the {BGK} model on general polygonal meshes},
 fjournal = {Journal of Computational Physics},
 journal = {J. Comput. Phys.},
 issn = {0021-9991},
 volume = {422},
 note = {Id/No 109766},
 year = {2020},
 language = {English},
 doi = {10.1016/j.jcp.2020.109766},
 zbMATH = {7508389},
 Zbl = {1559.76336}
}

@book{cer-1988,
 author = {Cercignani, Carlo},
 title = {The {Boltzmann} equation and its applications},
 fseries = {Applied Mathematical Sciences},
 series = {Appl. Math. Sci.},
 issn = {0066-5452},
 volume = {67},
 isbn = {0-387-96637-4},
 year = {1988},
 publisher = {New York etc.: Springer-Verlag},
 language = {English},
 zbMATH = {192917},
 Zbl = {0646.76001}
}

@INPROCEEDINGS{hol-1965,
       author = {{Holway}, Jr., Lowell H.},
        title = "{Kinetic Theory of Shock Structure Using an Ellipsoidal Distribution Function}",
    booktitle = {Rarefied Gas Dynamics, Volume 1},
         year = 1965,
       editor = {{de Leeuw}, J.~H.},
       volume = {1},
        month = jan,
        pages = {193},
       adsurl = {https://ui.adsabs.harvard.edu/abs/1965rgd1.conf..193H}
}

@article{bar-gol-lev-1991,
author = {Bardos, Claude and Golse, Francois and Levermore, David},
year = {1991},
month = {04},
pages = {323-344},
title = {Fluid dynamic limits of kinetic equations. I. Formal derivations},
volume = {63},
journal = {Journal of Statistical Physics},
doi = {10.1007/BF01026608}
}

@book{bou-gol-pul-2000,
author = {Bouchut, François and Golse, Francois and Pulvirenti, Mario},
year = {2000},
month = {03},
pages = {},
title = {Kinetic equations and asymptotic theory},
volume = {4},
isbn = {2-84299-110-9}
}

@article{cap-ten-2025,
  title={A Coupled IMEX Domain Decomposition Method for High-Order Time Integration of the ES-BGK Model of the Boltzmann Equation},
  author={Caparello, Domenico and Tenna, Tommaso},
  journal={arXiv preprint arXiv:2512.03586},
  year={2025}
}

@article{zhe-str-2004,
    author = {Zheng, Y. and Struchtrup, Henning},
    year = {2004},
    month = {02},
    pages = {97-108},
    title = {Burnett equations for the ellipsoidal statistical BGK model},
    volume = {16},
    journal = {Continuum Mechanics and Thermodynamics},
    doi = {10.1007/s00161-003-0143-3}
}

@article{cap-par-rey-2025-2,
    author = {Caparello, Domenico and Pareschi, Lorenzo and Rey, Thomas},
    year = {2025},
    month = {11},
    pages = {},
    title = {High-Order Asymptotic-Preserving IMEX schemes for an ES-BGK model for Gas Mixtures},
    doi = {10.48550/arXiv.2511.22304},
    journal = {arXiv preprint arXiv:2511.22304}
}

@article{jin-li-2013,
    author = {Jin, Shi and Li, Qin},
    title = {A BGK-penalization-based asymptotic-preserving scheme for the multispecies Boltzmann equation},
    journal = {Numerical Methods for Partial Differential Equations},
    volume = {29},
    number = {3},
    pages = {1056-1080},
    doi = {https://doi.org/10.1002/num.21746},
    url = {https://onlinelibrary.wiley.com/doi/abs/10.1002/num.21746},
    eprint = {https://onlinelibrary.wiley.com/doi/pdf/10.1002/num.21746},
    year = {2013}
}

@article{deg-dim-mie-2010,
    title = {A multiscale kinetic–fluid solver with dynamic localization of kinetic effects},
    journal = {Journal of Computational Physics},
    volume = {229},
    number = {13},
    pages = {4907-4933},
    year = {2010},
    issn = {0021-9991},
    doi = {https://doi.org/10.1016/j.jcp.2010.03.009},
    url = {https://www.sciencedirect.com/science/article/pii/S0021999110001191},
    author = {Pierre Degond and Giacomo Dimarco and Luc Mieussens}
}

@article{deg-jin-mie-2005,
    title = {A smooth transition model between kinetic and hydrodynamic equations},
    journal = {Journal of Computational Physics},
    volume = {209},
    number = {2},
    pages = {665-694},
    year = {2005},
    issn = {0021-9991},
    doi = {https://doi.org/10.1016/j.jcp.2005.03.025},
    url = {https://www.sciencedirect.com/science/article/pii/S0021999105001932},
    author = {Pierre Degond and Shi Jin and Luc Mieussens}
}

@article{deg-liu-mie-2006,
    author = { Degond, Pierre and  Liu, Jian‐Guo and  Mieussens, Luc},
    title = {Macroscopic Fluid Models with Localized Kinetic Upscaling Effects},
    journal = {Multiscale Modeling \& Simulation},
    volume = {5},
    number = {3},
    pages = {940-979},
    year = {2006},
    doi = {10.1137/060651574},
    URL = {https://doi.org/10.1137/060651574},
    eprint = {https://doi.org/10.1137/060651574}
}

@article{deg-dim-mie-2007,
    title = {A moving interface method for dynamic kinetic–fluid coupling},
    journal = {Journal of Computational Physics},
    volume = {227},
    number = {2},
    pages = {1176-1208},
    year = {2007},
    issn = {0021-9991},
    doi = {https://doi.org/10.1016/j.jcp.2007.08.027},
    url = {https://www.sciencedirect.com/science/article/pii/S0021999107003890},
    author = {Pierre Degond and Giacomo Dimarco and Luc Mieussens},
}

@article{bal-gar-flo-bos-2020,
 author = {Balsara, Dinshaw S. and Garain, Sudip and Florinski, Vladimir and Boscheri, Walter},
 title = {An efficient class of {WENO} schemes with adaptive order for unstructured meshes},
 fjournal = {Journal of Computational Physics},
 journal = {J. Comput. Phys.},
 issn = {0021-9991},
 volume = {404},
 year = {2020},
 language = {English},
 doi = {10.1016/j.jcp.2019.109062},
 zbMATH = {7302779},
 Zbl = {1453.65208}
}

@article{hu-shu-1999,
 author = {Hu, Changqing and Shu, Chi-Wang},
 title = {Weighted essentially non-oscillatory schemes on triangular meshes},
 fjournal = {Journal of Computational Physics},
 journal = {J. Comput. Phys.},
 issn = {0021-9991},
 volume = {150},
 number = {1},
 pages = {97--127},
 year = {1999},
 language = {English},
 doi = {10.1006/jcph.1998.6165},
 zbMATH = {1278951},
 Zbl = {0926.65090}
}

@article{bos-cho-rus-2022,
 author = {Boscarino, Sebastiano and Cho, Seung Yeon and Russo, Giovanni},
 title = {A local velocity grid conservative semi-{Lagrangian} schemes for {BGK} model},
 fjournal = {Journal of Computational Physics},
 journal = {J. Comput. Phys.},
 issn = {0021-9991},
 volume = {460},
 year = {2022},
 language = {English},
 doi = {10.1016/j.jcp.2022.111178},
 zbMATH = {7525158},
 Zbl = {1559.76816}
}

@article{bos-par-rus-2013,
 author = {Boscarino, S. and Pareschi, L. and Russo, G.},
 title = {Implicit-explicit {Runge}--{Kutta} schemes for hyperbolic systems and kinetic equations in the diffusion limit},
 fjournal = {SIAM Journal on Scientific Computing},
 journal = {SIAM J. Sci. Comput.},
 issn = {1064-8275},
 volume = {35},
 number = {1},
 pages = {a22--a51},
 year = {2013},
 language = {English},
 doi = {10.1137/110842855},
 zbMATH = {6159085},
 Zbl = {1264.65150}
}

@Article{mie-2000,
  author  = {Mieussens, Luc},
  journal = {Math. Models Method Appl. Sci.},
  title   = {Discrete velocity model and implicit scheme for the {BGK} equation of rarefied gas dynamic},
  year    = {2000},
  month   = {11},
  volume  = {10},
  doi     = {10.1142/S0218202500000562},
}

@Article{gam-tha-2009,
  author  = {Gamba, Irene and Tharkabhushanam, Sri},
  journal = {J. Comput. Phys.},
  title   = {Spectral-Lagrangian methods for Collisional Models of Non-Equilibrium Statistical States},
  year    = {2009},
  month   = {04},
  pages   = {2012-2036},
  volume  = {228},
  doi     = {10.1016/j.jcp.2008.09.033},
}

@article{cai-li-2010,
title = {An h-adaptive mesh method for Boltzmann-BGK/hydrodynamics coupling},
journal = {Journal of Computational Physics},
volume = {229},
number = {5},
pages = {1661-1680},
year = {2010},
issn = {0021-9991},
doi = {https://doi.org/10.1016/j.jcp.2009.10.050},
url = {https://www.sciencedirect.com/science/article/pii/S0021999109006123},
author = {Zhenning Cai and Ruo Li},
}

@article{wu-li-che-koo-tse-2006,
title = {Development and verification of a coupled DSMC–NS scheme using unstructured mesh},
journal = {Journal of Computational Physics},
volume = {219},
number = {2},
pages = {579-607},
year = {2006},
issn = {0021-9991},
doi = {https://doi.org/10.1016/j.jcp.2006.04.013},
url = {https://www.sciencedirect.com/science/article/pii/S0021999106001938},
author = {J.-S. Wu and Y.-Y. Lian and G. Cheng and R.P. Koomullil and K.-C. Tseng},
}

@book{mit-1970,
  title={{Analytic Inequalities}},
  author={Mitrinovic, Dragoslav S},
  volume={1},
  year={1970},
  publisher={Springer}
}

@article{cli-mcd-2009,
 author = {Groth, Clinton P. T. and McDonald, James G.},
 title = {Towards physically realizable and hyperbolic moment closures for kinetic theory},
 fjournal = {Continuum Mechanics and Thermodynamics},
 journal = {Contin. Mech. Thermodyn.},
 issn = {0935-1175},
 volume = {21},
 number = {6},
 pages = {467--493},
 year = {2009},
 language = {English},
 doi = {10.1007/s00161-009-0125-1},
 zbMATH = {5679076},
 Zbl = {1184.76839}
}

\end{document}